\documentclass[12pt]{amsart}
\usepackage{epsf}
\usepackage{psfrag}
\usepackage{fullpage}
\usepackage{float}
\usepackage{mathrsfs}
\usepackage{amsfonts}
\usepackage[centertags]{amsmath}
\usepackage{amssymb}
\usepackage{amsthm}
\usepackage{graphicx}
\usepackage{float}
\usepackage[all]{xy}
\usepackage{tikz-cd}
\usepackage{tikz}
\usetikzlibrary[shapes]
\usepackage{multirow}
\usepackage{caption}
\usepackage{xparse}
\usepackage{lscape}
\usepackage{dsfont}
\usepackage{inputenc}
\usepackage{subcaption} 
\usepackage{enumerate,enumitem}
\usepackage[margin=1in]{geometry} 
\usetikzlibrary{matrix,arrows,decorations.pathmorphing}
\usepackage{listing} 
\usepackage{hyperref}
\hypersetup{
bookmarksnumbered,
pdfstartview={FitH},
breaklinks=true,
linkcolor=blue,
urlcolor=blue,
citecolor=blue,
bookmarksdepth=2
}
\usepackage{cleveref}
\usepackage{ytableau}
\usepackage{adjustbox}
\usepackage{booktabs}

\usepackage{mathtools}
 
\theoremstyle{plain}
   
   \newtheorem{theorem}{Theorem}[section]

   \newtheorem{lemma}[theorem]{Lemma}
   \newtheorem{corollary}[theorem]{Corollary}
   \newtheorem{conjecture}[theorem]{Conjecture}
   
\theoremstyle{definition}
   \newtheorem{definition}[theorem]{Definition}
   
   \newtheorem{example}[theorem]{Example}

   \newtheorem{remark}[theorem]{Remark}

\numberwithin{equation}{section}

\newcommand{\CC}{{\mathbb {C}}}
\newcommand{\QQ}{{\mathbb {Q}}}
\newcommand{\RR}{{\mathbb {R}}}
\newcommand{\ZZ}{{\mathbb {Z}}}

\newcommand{\PP}{{\mathbb {P}}}

\newcommand{\TT}{{\mathbb {T}}}

\newcommand{\ch}{{\operatorname{ch}}}

\newcommand{\SSYT}{{\operatorname{SSYT}}}

\DeclareMathOperator*{\diag}{diag}

\DeclareMathOperator*{\wt}{wt}

\DeclareMathOperator{\Br}{Br}

\DeclareMathOperator{\sgn}{sgn}

\DeclareMathOperator{\Gr}{Gr} 

\newcommand\scalemath[2]{\scalebox{#1}{\mbox{\ensuremath{\displaystyle #2}}}}

\newcommand\mutline[1]{-^{#1} \!\!\!\!\!\! - \!\!\!\!\!- \!\!\!-\,}

\newcommand{\mathsm}[1]{\text{\small$\displaystyle #1$}}

\begin{document}

\title{Tropical symmetries of cluster algebras}

\author{James Drummond, \"{O}mer G\"{u}rdo\u{g}an, and Jian-Rong Li}
 
\address{James Drummond, School of Physics and Astronomy, University of Southampton, Highfield, Southampton, SO17 1BJ, United Kingdom.}
\email{j.m.drummond@soton.ac.uk}
\address{\"{O}mer G\"{u}rdo\u{g}an, School of Physics and Astronomy, University of Southampton, Highfield, Southampton, SO17 1BJ, United Kingdom.}
\email{o.c.gurdogan@soton.ac.uk}
\address{Jian-Rong Li, Faculty of Mathematics, University of Vienna, Oskar-Morgenstern-Platz 1, 1090 Vienna, Austria}
\email{lijr07@gmail.com}

\begin{abstract}
We study tropicalisations of quasi-automorphisms of cluster algebras and show that their induced action on the g-vectors can be realized by tropicalising their action on the homogeneous $\hat{y}$ (or $\mathcal{X}$) variables of a chosen initial cluster. This perspective allows us to interpret the action on g-vectors as a change of coordinates in the tropical setting. Focusing on Grassmannian cluster algebras, we analyse tropicalisations of quasi-automorphisms in detail. We derive tropical analogues of the braid group action and the twist map on both g-vectors and tableaux. We introduce the notions of unstable and stable fixed points for quasi-automorphisms. 

As an application, we demonstrate that the counts of prime non-real tableaux with a fixed number of columns in $\mathrm{SSYT}(3, [9])$ and $\mathrm{SSYT}(4, [8])$, arising from the braid group action on stable fixed points, are governed by Euler's totient function. Furthermore, we apply our findings to scattering amplitudes in physics, providing a novel interpretation of the square root associated with the four-mass box integral via stable fixed points of quasi-automorphisms of the Grassmannian cluster algebra $\CC[\Gr(4,8)]$.
\end{abstract}

\maketitle



\setcounter{tocdepth}{1}
\tableofcontents

\section{Introduction}

Cluster algebras, introduced by Fomin and Zelevinsky in \cite{FZ02}, are commutative rings generated from an initial combinatorial object called a seed, together with an iterative procedure known as mutation. This framework has since revealed deep connections across diverse areas of mathematics. Notably, cluster algebras play a central role in the representation theory of finite-dimensional algebras \cite{BBM06, Kel10} and KLR algebras \cite{KKKO}, as well as in Lie theory \cite{BFZ05, GLS13, HL10}, mirror symmetry \cite{GHKK18}, Poisson geometry \cite{GSV03}, quantum affine algebras \cite{HL10}, stability conditions in the sense of Bridgeland \cite{Bri07}, symplectic geometry \cite{CG22}, and Teichmüller theory \cite{FG06}. Beyond pure mathematics, cluster structures have found significant applications in physics, particularly in the study of scattering amplitudes and the amplituhedron \cite{ABCGPT16, ALS21, PSW22, DFGK20, HP21}. This wide-ranging impact highlights the broad and unifying influence of cluster algebras.

Fraser \cite{Fra16} introduced the notion of quasi-homomorphisms of cluster algebras, generalizing the concept of cluster automorphisms originally defined in \cite{ASS12}. These quasi-homomorphisms, also referred to as twist automorphisms in \cite{KQW22}, map clusters to clusters and send cluster variables to cluster variables, modulo factors of frozen variables. Quasi-homomorphisms play a fundamental role as a tool for identifying cluster variables within cluster algebras.

The goal of this paper is to study a tropical version of quasi-automorphisms of cluster algebras, focusing in particular on the tropicalisation of braid group actions and twist maps within Grassmannian cluster algebras.

Every quasi-automorphism $\mathcal{A} \to \mathcal{A}$ of a cluster algebra induces a map on the set of g-vectors of $\mathcal{A}$. We show that this induced map on g-vectors can be obtained by tropicalising the corresponding map on the homogeneous $\hat{y}$-variables (or $\mathcal{X}$-variables) associated with a fixed initial cluster of $\mathcal{A}$. This perspective allows us to interpret the action on g-vectors as arising from a change of coordinates. We first establish this result in the more general setting of a quasi-isomorphism between two cluster algebras (see Theorem \ref{thm:quasi automorphism sends g vectors of cluster variables to g vectors of cluster variables}), and then specialise to quasi-automorphisms of a single cluster algebra (see Corollary \ref{cor-quasiauto}).

The concept of g-vectors is generalised to tropical points in \cite{Qin17}, see Remark \ref{remark:tropical points and g-vectors}. Quasi-homomorphisms send tropical points to tropical points. Note that not every tropical point is the g-vector of a cluster variable, and it is important to classify which tropical points arise as g-vectors of cluster variables.

Tropicalisations of quasi-automorphisms of cluster algebras provide a more effective method to compute cluster variables. While using quasi-automorphisms directly to compute cluster variables can be complicated, applying tropical quasi-automorphisms to compute the g-vectors of cluster variables is considerably simpler. Knowing the g-vectors of cluster variables allows, in principle, the recovery of the cluster variables themselves. For example, consider the cluster algebra $K_0(\mathcal{C}_{\ell})$ introduced by Hernandez and Leclerc \cite{HL10}, where $\mathcal{C}_{\ell}$ is a certain subcategory of the category of finite-dimensional modules over a quantum affine algebra, and $K_0(\mathcal{C}_{\ell})$ denotes its Grothendieck ring. In this setting, cluster variables can be recovered using the $q$-characters of simple modules corresponding to the g-vectors. For $k \le n$, the Grassmannian $\Gr(k,n)$ is the variety of $k$-dimensional subspaces in an $n$-dimensional vector space. Scott \cite{Sco} proved that the coordinate ring $\CC[\Gr(k,n)]$ carries a cluster algebra structure. This coordinate ring is known as the Grassmannian cluster algebra. For $\CC[\Gr(k,n)]$, cluster variables can be recovered via the formula for $\ch(T)$ in Theorem 5.8 of \cite{CDFL}, where $T$ is the tableau associated with a g-vector obtained through tropical quasi-automorphisms.

We apply our general approach to the tropicalisations of quasi-automorphisms of Grassmannian cluster algebras. In particular, we define a braid group action on the set of g-vectors of the Grassmannian cluster algebra $\CC[\Gr(k,n)]$ by tropicalising Fraser's braid group action on $\CC[\Gr(k,n)]$. Using the correspondence between g-vectors and semistandard Young tableaux of rectangular shape \cite{CDFL}, we thereby obtain a braid group action on semistandard tableaux, see Sections \ref{sec:g vectors and tableaux} and \ref{sec:braid group actions on tableaux}.

Berenstein, Fomin, and Zelevinsky introduced a twist map on unipotent cells of a semisimple group of simply-laced type \cite{BFZ96}. Marsh and Scott \cite{MS16} subsequently defined a twist map on $\Gr(k,n)$. This twist map on $\Gr(k,n)$ is closely related to Berenstein, Fomin, and Zelevinsky's twist map on a unipotent cell of $SL_n(\CC)$; see Section 5 of \cite{MS16} for details. In this paper, we define a tropical version of the twist map on the Grassmannian cluster algebra. Additionally, we introduce tropical versions of the cyclic rotation and reflection maps on the Grassmannian cluster algebra.

Along the way, we discover an interesting conjectural relations (Conjecture \ref{conj:theta g is g of theta of tau inverse of b}) 
involving the reflection map. For example, for any cluster variable $b$ in $\CC[\Gr(k,n)]$, we have that 
\begin{align*}
\theta\tau({\bf g}(b)) = - {\bf g}(b)^\pi,
\end{align*}  
where $\pi$ is the permutation which reverses the order of the entries of a g-vector. This conjecture has an interesting corollary: there exists a seed whose set of 
g-vectors consists precisely of the negatives of the standard unit vectors (recall that the 
g-vectors of the initial cluster variables are the unit vectors). Indeed, starting from an initial seed and applying $\theta \circ \tau$ to the g-vectors of the initial cluster variables produces such a seed, as follows from the identity $\theta\tau({\bf g}(b)) = - {\bf g}(b)^\pi$.

We introduce the concepts of unstable fixed points and stable fixed points of quasi-automorphisms of cluster algebras. Furthermore, we conjecture that for any $(k,n)$, all stable fixed points are prime non-real elements in the dual canonical basis of $\CC[\Gr(k,n)]$; see Conjecture \ref{conj:stable fixed points and unstable fixed points prime elements and cluster variables}.

Recently, Sakamoto \cite{Sak25} classified simple modules that are simultaneously real and imaginary in certain module categories of quantum affine algebras arising in monoidal categorifications of cluster algebras of affine type. These include, in particular, the cluster algebras $\CC[\Gr(3,9)]$ and $\CC[\Gr(4,8)]$. To the best of our knowledge, the classification of imaginary prime modules in the categories arising from monoidal categorifications of cluster algebras, including the categories considered in \cite{Sak25}, remains open.

We explain that the numbers of tableaux with a given number of columns in $\SSYT(3, [9])$ and $\SSYT(4, [8])$, obtained from the braid group action on stable fixed points, are determined by Euler’s totient function; see Theorem \ref{thm:prime non real modules Gr39 and Gr48}. This gives a large family of imaginary prime modules in monoidal categorifications of $\CC[\Gr(3,9)]$ and $\CC[\Gr(4,8)]$ via the correspondence between simple modules and tableaux in Section \ref{subsec:correspondence between modules and tableaux}. We conjecture that all imaginary prime tableaux in $\SSYT(3, [9])$ and $\SSYT(4, [8])$ can be obtained from these imaginary prime tableaux (also all imaginary prime modules in monoidal categorifications of $\CC[\Gr(3,9)]$ and $\CC[\Gr(4,8)]$ can be obtained from these imaginary prime modules), see Conjecture \ref{conj: all imaginary prime tableaux in Gr39 and Gr48}. 

We explain a connection between braid group actions on Grassmannian cluster algebras and scattering amplitudes in physics. In particular, we provide a new interpretation of the well-known square root associated with the four-mass box using stable fixed points of quasi-automorphisms of cluster algebras, see Section \ref{sec:connection to scattering amplitudes}.

In \cite{CDFL}, it is proved that the dual canonical basis of $\CC[\Gr(k,n,\sim)]$ is given by $\{\widetilde{\ch}(T) : T \in \SSYT(k,[n],\sim)\}$, where $\widetilde{\ch}$ is a map from $\SSYT(k,[n],\sim)$ to $\CC[\Gr(k,n,\sim)]$; see Section \ref{sec:Grassmannian cluster algebras}. It is conjectured in \cite{CDFL} that the dual canonical basis of $\CC[\Gr(k,n)]$ is $\{\ch(T) : T \in \SSYT(k,[n])\}$. In Section \ref{sec:Grassmannian cluster algebras}, we prove this conjecture; see Theorem \ref{thm:dual canonical basis of CGrkn}.

The paper is organised as follows. In Section \ref{sec:tropicalisation of quasi homomorphisms}, we define a tropical version of quasi-automorphisms of cluster algebras. Section \ref{sec:Grassmannian cluster algebras} recalls results on Grassmannian cluster algebras, while Section \ref{sec:quantum affine algebras} reviews representations of quantum affine algebras. In Section \ref{sec:braid group actions on g vectors}, we define braid group actions on g-vectors, and in Section \ref{sec:braid group actions on tableaux}, we define these actions on tableaux. Section \ref{sec:Tropical cyclic rotations and promotions of tableaux} studies tropical cyclic rotations and promotions of tableaux. In Section \ref{sec:tropical Marsh Scott twist map}, we study the tropical version of Marsh and Scott's twist maps on Grassmannian cluster algebras. In Section \ref{sec:tropical reflection map}, we study tropical reflection maps and evacuations of tableaux. In Section \ref{sec:fixed points of quasi-automorphisms}, we introduce unstable and stable fixed points of quasi-automorphisms of cluster algebras. Section \ref{sec:unstable fixed points and cluster monomials} investigates unstable fixed points of braid group actions on Grassmannian cluster algebras. In Section \ref{sec:stable fixed points}, we apply stable fixed points to construct prime non-real elements in the dual canonical basis of $\CC[\Gr(k,n)]$ and prove that tableaux arising from braid group actions on stable fixed points in $\CC[\Gr(3,9)]$ and $\CC[\Gr(4,8)]$ are governed by Euler’s totient function. Finally, Section \ref{sec:connection to scattering amplitudes} explains connections of quasi-automorphisms of Grassmannian cluster algebras with scattering amplitudes in physics.

\subsection*{Acknowledgements}
We are very grateful to Lara Bossinger for many helpful and insightful discussions. JMD and \"OCG, are supported by the STFC consolidated grant ST/X000583/1. \"OCG is also supported by the Royal Society University Research Fellowship URF\textbackslash R1\textbackslash221236. JRL was supported by the Austrian Science Fund (FWF): P-34602, Grant DOI: 10.55776/P34602, and PAT 9039323, Grant-DOI 10.55776/PAT9039323.

\section{Tropicalisation of quasi-homomorphisms of cluster algebras} \label{sec:tropicalisation of quasi homomorphisms}

In this section, we prove that the map on g-vectors induced by any quasi-isomorphism of cluster algebras can be obtained by using the induced map on the $\hat{y}$ variables (or cluster $\mathcal{X}$-coordinates). 

\subsection{Cluster algebras}
We begin by recalling the definition of cluster algebras and several related notions, as introduced by Fomin and Zelevinsky \cite{FZ02, FZ07}. 

Let $(\PP, \oplus, \cdot)$ be a semifield, i.e. an abelian multiplicative group endowed with a binary operation $\oplus$ which is commutative associative and distributive with respect to the multiplication in $\PP$. 

For $n \in \ZZ_{\ge 1}$, we denote $[n]=\{1, \ldots, n\}$, and for $n \le m \in \ZZ$, denote $[n,m]=[n, n+1, \ldots, m]$. Denote $[x]_+ = \max(x, 0)$ and 
\begin{align*}
\sgn(x)=\begin{cases}
    1, & x>0, \\
    0, & x=0, \\
    -1, & x<0.
\end{cases}    
\end{align*}

Choose $m \ge n \in \ZZ_{>0}$. Let $\mathcal{F}$ be an ambient field of rational functions in $n$ independent variables with coefficients in $\QQ\PP$. A labelled seed in $\mathcal{F}$ is a triple $\Sigma = ({\bf x}, {\bf y}, B)$, where 
\begin{itemize}
    \item ${\bf x}=(x_1, \ldots, x_n)$ is an ordered $n$-tuple of elements of $\mathcal{F}$ forming a free generating set,

    \item ${\bf y}=(y_1, \ldots, y_n)$ is an ordered $n$-tuple of elements of $\PP$,

    \item $B = (b_{ij})$ is an $n \times n$ integer matrix which is skew-symmetrizable\footnote{i.e. there exists a diagonal matrix $D= {\rm diag}(d_1,\ldots,d_n)$ with $d_i \in \mathbb{N}$ such that $DB$ is skew symmetric.}  
\end{itemize}
The $n$-tuple ${\bf x}$ is called the labelled cluster of $\Sigma$, ${\bf y}$ is called the labelled coefficient tuple, and $B$ is called the exchange matrix. 

For $k \in [n]$, the seed mutation $\mu_k$ sends $\Sigma$ to $\mu_k(\Sigma) = ({\bf x}', {\bf y}', B')$, where 
\begin{itemize}
    \item ${\bf x}' = (x_1', \ldots, x_n')$, $x_j'=x_j$ if $j \ne k$, and 
\begin{align}
x_k' = \frac{y_k \prod_{i} x_i^{[b_{ik}]_+} + \prod_{i} x_i^{[-b_{ik}]_+} }{(y_k \oplus 1)x_k},         
\end{align}

\item ${\bf y}' = (y_1', \ldots, y_n')$, and
\begin{align}
y_j' = \begin{cases}
y_k^{-1}, & j=k, \\
y_jy_k^{[b_{kj}]_+}(y_k \oplus 1)^{-b_{kj}}, & j \ne k,
\end{cases}
\end{align}

\item $B' = (b_{ij}')$, and
\begin{align}
b_{ij}' = \begin{cases}
-b_{ij}, & i=k \text{ or } j=k, \\
b_{ij}+\sgn(b_{ik})[b_{ik}b_{kj}]_+, & \text{otherwise}. 
\end{cases}
\label{Bmutation}
\end{align}
\end{itemize}
The operations of flipping the sign and taking the transpose of $B$ both commute with the mutation rule (\ref{Bmutation}). 
%That is we have $(-B)' = - B'$ and $(B^T)' = (B')^T$, with the latter relying on $\sgn (b_{ij}) = - \sgn (b_{ji})$, following from skew-symmetrisability. 
We define the opposite exchange matrix $B^{\rm op} = -B$ and the dual exchange matrix $B^\vee = - B^T$ obeying $(B^{\rm op})' = (B')^{\rm op}$ and $(B^\vee)' = (B')^\vee$.

Consider the $n$-regular tree (called labelled tree) $\TT_n$ whose edges are labelled by $1, \ldots, n$, so that the $n$ edges emanating from each vertex receive different labels. A cluster pattern is an assignment of a labelled seed $\Sigma_t = ({\bf x}_t, {\bf y}_t, B_t)$ to every vertex $t \in \TT_n$ such that the labelled seeds assigned to the end-points of any edge
$t \mutline{k} t'$
are each obtained from the other by the mutation $\mu_k$. We usually write 
\begin{equation}
{\bf x}_t=(x_{1;t}, \ldots, x_{n;t})\,, \qquad  {\bf y}_t=(y_{1;t}, \ldots, y_{n;t})\,, \qquad B_t = (b_{ij}^t)_{n \times n}\,
\end{equation}
to denote the elements of the labelled seed $\Sigma_t$.
%and we sometimes distinguish the initial seed $\Sigma_{t_0}=({\bf x}_{t_0}, {\bf y}_{t_0}, B_{t_0})$ with the notation
%\begin{equation}
%{\bf x}_{t_0}=(x_1, \ldots, x_n)\,, \qquad {\bf y}_{t_0}=(y_1, \ldots, y_n)\,, \qquad B_{t_0}=B^0=(b_{ij}^0)\,.
%\label{initialnotation}
%\end{equation}
Given a cluster pattern, we often use the shorthand notation $t'=\mu_k t$ whenever $t \mutline{k} t'$ so that $\mu_k \Sigma_t = \Sigma_{\mu_k t}$. Note that every pair of vertices $(t,t')$ induces a unique sequence of mutations $\mu = \mu_{k_r} \ldots \mu_{k_1}$ such that $t' = \mu t$.

We can obtain another cluster pattern by applying a permutation $\pi$ to the $n$ elements of each labelled seed, as well as each edge of the labelled tree $\mathbb{T}_n$. We define $\Sigma^\pi = ({\bf x}^\pi, {\bf y}^\pi, B^\pi)$ such that
\begin{equation}
    {\bf x}^\pi = (x_{\pi(1}\\\ldots, x_{\pi(n)})\,, \quad {\bf y}^\pi = (y_{\pi(1)},\ldots,y_{\pi(n)})\,, \quad B^\pi = (b^\pi_{ij}) = (b_{\pi(i) \pi(j)})\,. \label{permseed}
\end{equation}
Then replacing each $\Sigma_t$ with $\Sigma_t^\pi$ and each edge label $k$ with $\pi(k)$ in $\mathbb{T}_n$ gives another cluster pattern.

Two labelled seeds $\Sigma=({\bf x},{\bf y},B)$ and $\Sigma' = ({\bf x}', {\bf y}', B')$ are said to define the same (unlabelled) seed if there exists a permutation $\pi$ of the $n$ labels such that $\Sigma = \Sigma'^\pi$.
%\begin{equation}
%    x_i = x'_{\pi(i)}\,, \qquad y_i = y_{\pi(i)}'\,, \qquad b_{ij} = b'_{\pi(i),\pi(j)}\,.
%\end{equation}
In such a case we write $\Sigma \sim \Sigma'$ and denote the equivalence class as $[\Sigma] = [\Sigma']$. The corresponding unordered sets $[{\bf x}] = [{\bf x}'] = \{x_1,\ldots,x_n\}$ are called clusters.

Given a cluster pattern, let $X = \cup_{t \in \TT_n} {\bf x}_t = \{ x_{i,t}: t \in \TT_n, i \in [n] \}$ be the union of clusters of all the seeds in the pattern. The elements $x_{i,t} \in X$ are called cluster variables. Two cluster variables are compatible if they are in a common cluster. Cluster monomials are the products of compatible cluster variables. The cluster algebra $\mathcal{A}$ associated with a given cluster pattern is the $\ZZ\PP$-subalgebra $\ZZ\PP[X]$ of $\mathcal{F}$ generated by all cluster variables. The number $n$ of cluster variables in a cluster is called the rank of the cluster algebra $\mathcal{A}$. 

Given a cluster pattern, we have three other cluster patterns, each with their own sets of variables and coefficients associated to the opposite, dual and opposite dual exchange matrices. We denote such labelled seeds in the obvious way, e.g. $\Sigma^{\rm op}_t = ({\bf x}^{\rm op}_t , {\bf y}^{\rm op}_t , B^{\rm op}_t)$ etc.

For a labelled seed $\Sigma_t=({\bf x}_t, {\bf y}_t, B_t)$ of a cluster algebra, denote  
\begin{align} \label{eq:definition of y hat}
\hat{y}_{j;t} = y_{j;t} \prod_{i=1}^n x_{i;t}^{b_{ij}^t}\,,
%\quad \hat{y}_j = \hat{y}_{j;t_0}, 
\quad 1 \le j \le n\,.
%\hat{y}_\Sigma(x_i) = \frac{\prod_{j: i \to j} x_j}{\prod_{j:j\to i} x_j},
\end{align} 
These $\hat{y}$ variables mutate in the same way as the coefficients, except in $\mathcal{F}$ instead of $\mathbb{P}$. In other words, under mutation on node $k$ we have
\begin{align}\label{yhatmutation}
    \hat{y}_j' = 
    \begin{cases}
        \hat{y}_k^{-1}\,, & j=k\,, \\
        \hat{y}_j \hat{y}_k^{[b_{kj}]_+}(\hat{y}_k + 1)^{-b_{kj}}\,, & j \neq k\,.
    \end{cases}
\end{align}
For a labelled seed in the opposite cluster pattern we have the corresponding variables $\hat{y}^{\rm op}_{j;t}$. We sometimes make use of the notation $\check{y}_{j;t}=\hat{y}^{\rm op}_{j;t}$. Note that the $\check{y}_j$ mutate in the same way as $\hat{y}_j^{-1}$ since
\begin{align}
    \bigl[\hat{y}_j \hat{y}_k^{[b_{kj}]_+}(\hat{y}_k + 1)^{-b_{kj}}\bigr]^{-1} 
    %= \hat{y}_j^{-1} (\hat{y}_k^{-1})^{[b_{kj}]_+} (1/\hat{y}_k^{-1} + 1)^{b_{kj}} \notag 
    = \hat{y}_j^{-1} (\hat{y}_k^{-1})^{[b_{kj}]_+ - b_{kj}} (\hat{y}_k^{-1} + 1)^{b_{kj}} = \hat{y}_j^{-1} (\hat{y}_k^{-1})^{[b^{\rm op}_{kj}]_+} (\hat{y}_k^{-1} + 1)^{-b^{\rm op}_{kj}}\,. \notag
\end{align}

A tropical semifield ${\rm Trop}(u_j: j \in J)$ is an abelian group (written multiplicatively) freely generated by the elements $u_j$, $j \in J$, for some finite set $J$, and the addition $\oplus$ in  ${\rm Trop}(u_j: j \in J)$ is defined by 
\begin{align}\label{tropsemifield}
    \prod_j u_j^{a_j} \oplus \prod_j u_j^{b_j} = \prod_j u_j^{\min(a_j, b_j)}.
\end{align}

A cluster algebra is of geometric type if the coefficient semifield $\PP$ is a tropical semifield. For a cluster algebra of geometric type, denote the generators of $\PP$ by $x_{n+1}, \ldots, x_m$ for some $m \ge n$ and $\PP = {\rm Trop}(x_{n+1}, \ldots, x_m)$. The variables $x_{n+1},\ldots,x_m$ are called frozen variables (while the variables $x_1,\ldots,x_n$ are called active). The coefficients $y_{1;t}, \ldots, y_{n;t}$ in each labelled seed $\Sigma_t=({\bf x}_t, {\bf y}_t, B_t)$ can then be written as 
\begin{align*}
    y_{j;t} = \prod_{i=n+1}^m x_i^{b_{ij}^t},
\end{align*}
for some integers $b_{ij}^t$. The extended exchange matrix $\tilde{B}_t = (b_{ij}^t)_{m \times n}$ contains $B_t$ as an $n \times n$ submatrix as its first $n$ rows. In this case we refer to $B_t$ as the principal part of $\tilde{B}_t$. The exchange relation for cluster variables can be written as
\begin{align*}
x_k' = x_k^{-1}\left( \prod_{i=1}^m x_i^{[b_{ik}^t]_+} + \prod_{i=1}^m x_i^{[-b_{ik}^t]_+} \right),
\end{align*}
where $x_i = x_{i;t}$, $i \in [n]$. If we extend the matrix mutation (\ref{Bmutation}) for $B_t$ to $\tilde{B}_t$ then the mutation of coefficients is encoded in the mutation of $\tilde{B}_t$. 
We may also define $\mathcal{A}^{\rm op}$ via its extended exchange matrices $\tilde{B}_t^{\rm op} = - \tilde{B}_t$, encoding the mutation rules among its own active and frozen generators $x_i^{\rm op}$.

A cluster algebra has principal coefficients at a vertex $t_0$ if $\PP = {\rm Trop}(y_1, \ldots, y_n)$ and ${\bf y}_{t_0} = (y_1, \ldots, y_n)$. Such a cluster algebra is denoted by $\mathcal{A}_\bullet(B_{t_0})$. Note that in this case, the extended exchange matrix $\tilde{B}_t =(b^t_{ij})$ is a $2n \times n$ matrix with $B_t$ as its principal part. 

We refer to the $n \times n$ matrix $C_{t,t_0}$ with entries $c^t_{ij} = b^t_{n+i,j}$ as the coefficient matrix. The $j$th column of the coefficient matrix is referred to as the $c$-vector ${\bf c}_{j;t}$. The coefficient matrix for the labelled seed at vertex $t_0$ is the $n \times n$ identity matrix, $C_{t_0,t_0} = 1\!\!1_{n \times n}$. 

By repeated mutation, we can express any cluster variable $x_{\ell;t}$ in $\mathcal{A}_\bullet(B_{t_0})$ as a unique subtraction-free rational function $X_{\ell;t}$ of the initial seed data $(x_1,\ldots,x_n,y_1,\ldots,y_n)$,
\begin{equation}
    x_{\ell;t} = X_{\ell;t} = X_{\ell;t}^{B^0;t_0} \in \mathbb{Q}(x_1,\ldots,x_n,y_1,\ldots,y_n)\,.
\end{equation}
%The notation $X_{\ell;t}^{B^0;t_0}$ emphasises that the rational function depends not only on the choice of cluster variable $x_{\ell,t}$ but also on the choice of initial vertex $t_0$ at which $\mathcal{A}_\bullet(B_{t_0})$ has principal coefficients.

Note that the rational function $X_{\ell;t}$ above exhibits the Laurent phenomenon (c.f. Prop. 3.6 in \cite{FZ07}),
\begin{equation}
    X_{\ell;t}^{B^0;t_0} \in \mathbb{Z}[x_1^{\pm 1},\ldots,x_n^{\pm 1},y_1,\ldots,y_n]\,.
\end{equation}
If we specialise to the case where all of the $x_i$ are set to 1 we therefore obtain a polynomial (called the $F$-polynomial),
\begin{equation}
    F_{\ell;t}^{B^0,t_0}(y_1,\ldots,y_n) = X_{\ell;t}^{B^0;t_0}(1,\ldots,1,y_1,\ldots,y_n)\,.
\end{equation}
%\subsection{g-vectors}
The concept g-vector was introduced by Fomin and Zelevinsky in \cite{FZ07}. In Corollary 6.2 in \cite{FZ07}, it is shown that every cluster variable in a cluster algebra $\mathcal{A}$ can be written as
\begin{align*}
x_{l;t} = x_{1}^{g_1} \cdots x_n^{g_n} \frac{F_{\ell; t}^{B^0; t_0}|_{\mathcal{F}}( \hat{y}_1, \ldots, \hat{y}_n ) }{  F_{\ell; t}^{B^0; t_0}|_{\mathbb{P}}( y_1, \ldots, y_n ) },
\end{align*}
where $\hat{y}_j$ is defined in (\ref{eq:definition of y hat}). The vector 
\begin{equation}
{\bf g}_{l;t} = {\bf g}(x_{l,t};t_0) = (g_1, \ldots, g_n)
\end{equation} 
is called the g-vector of the cluster variable $x_{l;t}$ with respect to the vertex $t_0$. Note that the g-vector depends on both the cluster variable $x_{l;t}$ and the choice of initial vertex $t_0$, as we have emphasised by the notation ${\bf g}(x_{l;t};t_0)$ above. There will also be g-vectors associated to the opposite, dual and opposite dual exchange matrices, for which we use the notation 
\begin{equation}
{\bf g}^{\rm op}_{l;t} = {\bf g}^{\rm op}(x_{l;t} ; t_0) = {\bf g}(x^{\rm op}_{l,t};t_0) = (g^{\rm op}_1, \ldots, g^{\rm op}_n)
\end{equation}
etc.

Just as for the cluster variables, one can obtain all g-vectors from the initial vertex by mutation. The g-vectors of the initial vertex are given by unit vectors,
\begin{equation}\label{initialgs}
    {\bf g}_{l;t_0} = {\bf e}_l\,.
\end{equation}
The remaining g-vectors are then obtained by the following rule \cite{FZ07} for mutating to the vertex $t'$ from the vertex $t$ on node $k$,
\begin{align}\label{gvecmutation}
    {\bf g}_{l;t'} &= {\bf g}_{l;t}\,, \qquad l \neq k\,,  \\
    {\bf g}_{k;t'} &= - {\bf g}_{k,t} + \sum_{i=1}^n [b^t_{ik}]_+ {\bf g}_{i;t} - \sum_{j=1}^n [c^{t,t_0}_{jk}]_+ {\bf b}_j^0\,, \notag \\
    &= - {\bf g}_{k;t} + \sum_{i=1}^n [-b^t_{ik}]_+ {\bf g}_{i;t} - \sum_{j=1}^n [-c^{t,t_0}_{jk}]_+ {\bf b}_j^0\,. \notag
\end{align}
Here $c^{t,t_0}_{ij}$ are the entries of the $n \times n$ coefficient matrix $C_{t,t_0}$ while ${\bf b}_j^0$ is column $j$ of $B_{t_0}$. Note that the mutation rules make clear that, having fixed the initial vertex $t_0$ (and hence $C_{t_0,t_0} = 1\!\!1_{n \times n}$), the full set of g-vectors is then determined by the exchange matrix $B_{t_0}$ and the initial data (\ref{initialgs}). 

In analogy to the coefficient matrix we introduce the $n \times n$ matrix $G_{t,t_0}=(g^{t,t_0}_{ij})$ whose $j$th column is the g-vector ${\bf g}_{j;t}$ and we note that $G_{t_0,t_0} = 1\!\!1_{n \times n}$. Having fixed an initial vertex $t_0$ we may adjoin the matrices $C_{t,t_0}$ and $G_{t,t_0}$ to the data for each labelled seed $\Sigma_t$. If we then permute the cluster pattern to obtain a new cluster pattern in the sense described below (\ref{permseed}) then we replace each edge label $k$ with $\pi(k)$ in the labelled tree $\mathbb{T}_n$ and the data at each vertex are replaced via $(\Sigma_t,C_{t,t_0},G_{t,t_0}) \mapsto (\Sigma^\pi_t,C^\pi_{t,t_0},G^\pi_{t,t_0})$. Here $C^\pi_{t,t_0}$ and $G_{t,t_0}^\pi$ are defined similarly to $B_t^\pi$ in (\ref{permseed}).

\begin{remark} \label{remark:tropical points and g-vectors}
More generally, geometric type cluster algebras have certain bases formed by the so called pointed elements of the form ${\bf x}_t^{{\bf g}_t}F_t(\hat{y}_{1,t}, \ldots, \hat{y}_{n,t})$ (with respect to a seed $({\bf x}_t, \tilde{B}_t)$), where $F_t$ has constant term $1$, see \cite{Qin19}. Each such basis contains cluster monomials as a subset. With slight abuse of notation, we also call ${\bf g}_t$ a g-vector (it is possible that ${\bf x}_t^{{\bf g}_t}F_t(\hat{y}_{1,t}, \ldots, \hat{y}_{n,t})$ is not a cluster variable). For example, every element in the dual canonical basis of the Grassmannian cluster algebra $\CC[\Gr(k,n)]$ has a g-vector, see Section 7 of \cite{CDFL}. In particular, the elements in the dual canonical basis of $\CC[\Gr(k,n)]$ corresponding to limit rays have g-vectors. 
\end{remark}

\subsection{\texorpdfstring{Change of g-vectors under mutation of the initial seed as a tropicalised dual $\hat{y}$ mutation}{Change of g-vectors under mutation of the initial seed as a tropicalised dual y-hat mutation}}

It is natural to ask how the g-vectors transform under a mutation of the initial vertex $t_0$ to an adjacent one $t_1$. Here we prove a lemma relating such a change of g-vectors to the tropicalisation of the mutation of $\hat{y}$-coordinates in a cluster algebra with the dual exchange matrix $B_t^{\vee} = - B_t^T$ in each seed. This result will be used in the rest of this section to define the tropicalisation of quasi-automorphisms of cluster algebras. 

Here, by tropicalisation, we mean a slightly different notion from the definition of tropical semifield introduced above in (\ref{tropsemifield}). We will define the tropicalisation of the field of rational functions in the $\hat{y}$ variables as follows. Given a rational function,
\begin{equation}
R(\hat{y}_1,\ldots,\hat{y}_n) = \frac{P(\hat{y}_1,\ldots,\hat{y}_n)}{Q(\hat{y}_1,\ldots,\hat{y}_n)}
\end{equation}
for $P$ and $Q$ polynomials, we define
\begin{equation}
{\rm Trop}^+ R(\hat{y}_1,\ldots,\hat{y}_n) = {\rm Trop}^+ P(\hat{y}_1,\ldots,\hat{y}_n) -  {\rm Trop}^+ Q(\hat{y}_1,\ldots,\hat{y}_n)
\end{equation} where ${\rm Trop}^+$ applied to a polynomial results in a piecewise linear function of the tropical variables $v_i = {\rm Trop^+} \hat{y}_i$, obtained by replacing multiplication by addition and addition by taking the maximum. For example,
\begin{equation}
{\rm Trop}^+ (1 + \hat{y}_1^2 + 2 \hat{y}_1 \hat{y}_2) = {\rm max}(0,2 v_1, v_1 + v_2)\,.
\end{equation}
We similarly define ${\rm Trop}^-$ where the maximum is replaced by the minimum and we obtain a piecewise linear function in variables $w_i = {\rm Trop}^- \hat{y}_i$.

%For a cluster algebra $\mathcal{A}_{\bullet}(B_t)$ with skew-symmetrisable exchange matrix $B_t$, the dual of $\mathcal{A}_{\bullet}(B_t)$ is the cluster algebra $\mathcal{A}_{\bullet}(B_t^D)$, where $B_t^D = - B_t^T$. 

\begin{lemma} \label{lem:for every mutation, change of g-vector is the same as change of tropicalisation of x-coordinates}
For every mutation $t \mutline{k} t'$ in a cluster algebra with skew-symmetrisable exchange matrix $B_{t}$, the change of the g-vectors with respect to the two vertices $t, t'$ is the same as the `max' tropicalisation of the mutation of the $\hat{y}$-variables induced by the dual matrix $B_{t}^\vee$. The change of the g-vectors associated to the opposite exchange matrix $B_t^{\rm op} = - B_t$ is given by the `min' tropicalisation of the same mutation.
\end{lemma}

\begin{proof}

We will make significant use of Conjecture 7.12 in \cite{FZ07}. This conjecture was proven in Theorem 1.7 in \cite{DWZ10}, Theorem 3.7 in \cite{Pla11} for cluster algebras with skew-symmetric exchange matrices, and, as proven in Theorem 2.11 in \cite{Rea20}, follows more generally from Theorem 4.9, Corollary 5.9 in \cite{GHKK18} for cluster algebras with skew-symmetrisable exchange matrices. The result describes the how the g-vectors change under a mutation of the initial seed to an adjacent seed. Let $t'$ be the vertex obtained by mutating from vertex $t$ on node $k$, $t \mutline{k} t'$. For any cluster variable $x$ in $\mathcal{A}$, let ${\bf g}(x;t) = (g_1,\ldots,g_n)$ and let ${\bf g}(x;t') = (g_1',\ldots,g_n')$. Then the two vectors are related by
\begin{align}
g_j' = 
\begin{cases}
- g_k\,, & j=k\,, \\
g_j + [b_{jk}]_+ g_k - b_{jk} \,{\rm min}(g_k,0), & j\, \neq k\,. 
\end{cases}
\label{gveccoordchange}
\end{align}
Now consider a cluster algebra with the dual exchange matrix $B^\vee = (b^\vee_{ij})$ for each labelled seed and denote its $\hat{y}$ variables by $\hat{y}^\vee$.  If we consider the `max' convention tropicalised coordinates $v^\vee_j = {\rm Trop}^+ \hat{y}^\vee_j$ we then have the tropicalised mutation rule
\begin{align}
    v^{\vee \prime}_j = 
    \begin{cases}
        -v^\vee_k \,, & j=k\,, \\
        v^\vee_j + [b^\vee_{kj}]_+ v^\vee_k - b^\vee_{kj} \,{\rm max}(v^\vee_k,0)\,, & j \neq k\,.
    \end{cases}
\end{align}
Now, using the definition of the dual exchange matrix $b^\vee_{ij} = - b_{ji}$, we have
\begin{align}
    [b^\vee_{kj}]_+ v^\vee_k - b^\vee_{kj} \, {\rm max}(v^\vee_k,0) &= [-b_{jk}]_+ v^\vee_k + b_{jk} \, {\rm max}(v^\vee_k,0) \, \\
    &= ([b_{jk}]_+  - b_{jk}) v^\vee_k - b_{jk} \, {\rm min}(-v^\vee_k,0) \notag \\
    &= [b_{jk}]_+ v^\vee_k - b_{jk} \, {\rm min}(0,v^\vee_k)\,.
    \notag 
\end{align}
Therefore we can rewrite the mutation rule for the tropical coordinates $v^\vee_j$ as
\begin{align}\label{tropyhatmutation}
    v^{\vee \prime}_j = 
    \begin{cases}
        -v^\vee_k \,, & j=k\,, \\
        v^\vee_j + [b_{jk}]_+ v^\vee_k - b_{jk} \, {\rm min}(v^\vee_k,0)\,, & j \neq k\,,
    \end{cases}
\end{align}
We now observe that the piecewise linear map (\ref{tropyhatmutation}) is exactly the change of g-vectors given by (\ref{gveccoordchange}). 
Note that to make use of the tropicalised transformation of the $\hat{y}^\vee$ variables to compute the transformation of the g-vectors, we assume that there are no algebraic relations among the $\hat{y}^\vee_j$ of a given seed so that the form of the RHS of (\ref{tropyhatmutation}) is unique.

If we consider the `min' convention tropicalised coordinates $w^\vee_j = {\rm Trop}^- \hat{y}^\vee_j$ we then have the tropicalised mutation rule
\begin{align}\label{mintropyhatmutation}
    w^{\vee \prime}_j &= 
    \begin{cases}
        -w^\vee_k \,, & j=k\,, \\
        w^\vee_j + [b^\vee_{kj}]_+ w^\vee_k - b^\vee_{kj} \,{\rm min}(w^\vee_k,0)\,, & j \neq k\,.
    \end{cases}  \\
    &= \begin{cases}
        -w^\vee_k \,, & j=k\,, \\
        w^\vee_j + [b^{\rm op}_{jk}]_+ w^\vee_k - b_{jk}^{\rm op} \,{\rm min}(w^\vee_k,0)\,, & j \neq k\,.
    \end{cases}
    \notag
\end{align}
We observe that this is exactly the equivalent transformation rule between the g-vectors ${\bf g}(x^{\rm op};t) = {\bf g}^{\rm op}(x;t)$ and ${\bf g}(x^{\rm op};t') = {\bf g}^{\rm op}(x;t')$ associated to the opposite exchange matrix $B^{\rm op}_t = - B_t$. Note that the rule (\ref{mintropyhatmutation}) is the same as the `max' tropicalisation applied to the mutation rule of the $\hat{y}^{\rm op}$ variables while the rule (\ref{tropyhatmutation}) is the same as the `min' rule applied to the $\hat{y}^{\rm op}$ variables.

\end{proof}

We introduce the notation $Q^\pm_{t',t}$ for the piecewise linear maps in (\ref{tropyhatmutation}) and (\ref{mintropyhatmutation}) above,
\begin{align}\label{Qmaps}
    Q^+_{t',t} &: (v^\vee_1,\ldots,v^\vee_n) \mapsto (v_1^{\vee \prime}, \ldots, v_n^{\vee \prime})\,,  &&Q^+_{t',t} {\bf g}(x;t) = {\bf g}(x;t')\,,  \\
    Q^-_{t',t} &: (w^\vee_1,\ldots,w^\vee_n) \mapsto (w_1^{\vee \prime}, \ldots, w_n^{\vee \prime})\,,  &&Q^-_{t',t} {\bf g}^{\rm op}(x;t) = {\bf g}^{\rm op}(x;t')\,. 
    \notag
\end{align}
Since the maps arise from tropicalisation of mutation rules for the $\hat{y}$ variables, they compose in a straightforward way such that (\ref{Qmaps}) may be used when $t$ and $t'$ are related by a sequence of mutations. In particular we have
\begin{equation}
  Q^\pm_{t'',t'}   Q^\pm_{t',t} = Q^\pm_{t'',t}\,.
\end{equation}

We similarly introduce the notation $Q^{\vee \pm}_{t',t}$ for corresponding the piecewise linear maps obtained by tropicalising the $\hat{y}$ variables themselves rather than their duals. These maps act similarly on the g-vectors of cluster algebras with the dual exchange matrices,
\begin{align}
    Q^{\vee +}_{t',t} &: (v_1,\ldots,v_n) \mapsto (v_1', \ldots, v_n')\,,  &&Q^{\vee +}_{t',t} {\bf g}^\vee(x;t) = {\bf g}^\vee(x;t')\,,  \\
    Q^{\vee -}_{t',t} &: (w_1,\ldots,w_n) \mapsto (w_1', \ldots, w_n')\,,  &&Q^{\vee -}_{t',t} {\bf g}^{{\rm op}, \vee}(x;t) = {\bf g}^{{\rm op}, \vee}(x;t')\,.
    \notag
\end{align}

\begin{example}
    Let us consider the case of a cluster algebra of type $C_2$ and a cluster algebra based on the dual exchange matrix, of type $B_2$. Table \ref{C2examplegvecs} lists the six inequivalent seeds in type $C_2$ with the corresponding $C$-matrices and $G$-matrices, taking either $t_0$ or $t_1$ as the initial seed. Table \ref{B2examplegvecs} lists the same information for the dual $B_2$ type cluster algebra. In each seed we also list the $\hat{y}$-variables in terms of those of the seed $t_0$ in both cases.

    The $\hat{y}^\vee$-variables in seed $t_1$ are given in terms of those in seed $t_0$ in the second row of Table \ref{B2examplegvecs}. Tropicalising the corresponding rational functions with `max' convention we obtain
    \begin{align}
        v^\vee_{1;t_1} &= - v^\vee_1\,,  \\
        v^\vee_{2;t_1} &= v^\vee_1 + v^\vee_2 - \,{\rm max}(0,v^\vee_1)\,.
        \notag
    \end{align}
    Applying the piecewise linear map
    \begin{align}
        Q^+_{t_1,t_0}:  \begin{pmatrix} v^\vee_1 \\ v^\vee_2 \end{pmatrix} \mapsto \begin{pmatrix} v^\vee_{1;t_1} \\ v^\vee_{2;t_1}\end{pmatrix}
    \end{align}
    to the columns of the matrices $G_{t,t_0}$ in Table \ref{C2examplegvecs} we obtain the columns of the matrices $G_{t,t_1}$.
    
    Tropicalising the same rational functions with `min' convention we obtain
    \begin{align}
        w^\vee_{1;t_1} &= - w^\vee_1\,,  \\
        w^\vee_{2;t_1} &= w^\vee_1 + w^\vee_2 - \,{\rm min}(0,w^\vee_1)\,.
        \notag
    \end{align}
    Applying the piecewise linear map
    \begin{align}
       Q^-_{t_1,t_0}:    \begin{pmatrix} w^\vee_1 \\ w^\vee_2 \end{pmatrix} \mapsto \begin{pmatrix} w^\vee_{1;t_1} \\ w^\vee_{2;t_1}\end{pmatrix}
    \end{align}
    to the columns of the matrices $G^{\rm op}_{t,t_0}$ in Table \ref{C2opexamplegvecs} we obtain the columns of the matrices $G^{\rm op}_{t,t_1}$.

    Likewise, the $\hat{y}$-variables in seed $t_1$ are given in terms of those in seed $t_0$ in the second row of Table \ref{C2examplegvecs}. Tropicalising these rational functions with `max' convention we find
    \begin{align}
        v_{1;t_1} &= - v_1 \,,  \\
        v_{2;t_1} &= 2 v_1 + v_2 - 2 \,{\rm max}(0,v_1)\,.
        \notag
    \end{align}
    Applying the piecewise linear map
    \begin{align}
        Q^{\vee +}_{t_1,t_0}: \begin{pmatrix} v_1 \\ v_2 \end{pmatrix} \mapsto \begin{pmatrix} v_{1;t_1} \\ v_{2;t_1}\end{pmatrix}
    \end{align}
    to the columns of the matrices $G^\vee_{t,t_0}$ in Table \ref{B2examplegvecs} we obtain the columns of the matrices $G^\vee_{t,t_1}$.

    Tropicalising the same rational functions with `min' convention we find
    \begin{align}
        w_{1;t_1} &= - w_1 \,,  \\
        w_{2;t_1} &= 2 w_1 + w_2 - 2 \,{\rm min}(0,w_1)\,.
        \notag
    \end{align}
    Applying the piecewise linear map
    \begin{align}
       Q^{\vee -}_{t_1,t_0}: \begin{pmatrix} w_1 \\ w_2 \end{pmatrix} \mapsto \begin{pmatrix} w_{1;t_1} \\ w_{2;t_1}\end{pmatrix}
    \end{align}
    to the columns of the matrices $G^{{\rm op},\vee}_{t,t_0}$ in Table \ref{B2opexamplegvecs} we obtain the columns of the matrices $G^{{\rm op},\vee}_{t,t_1}$.
{
  \renewcommand{\arraystretch}{1.2}
\begin{table}
  \centering
  \begin{tabular}{cccccccc}
    \toprule
    $t$&$B_t$&$C_{t,t_0}$&$G_{t,t_0}$&$C_{t,t_1}$&$G_{t,t_1}$ & $\hat{y}_{1;t}$ & $\hat{y}_{2;t}$\\
    \midrule
    $t_0$&$ \biggl[ \mathsm{ \begin{matrix} 0 & 2 \\ -1 & 0 \end{matrix}}\biggr]$&$\mathsm{ \begin{bmatrix} 1 & 0 \\ 0 & 1 \end{bmatrix}}$&$\mathsm{ \begin{bmatrix} 1 & 0 \\ 0 & 1 \end{bmatrix}}$&$\mathsm{ \begin{bmatrix} -1 & 0 \\ 0 & 1 \end{bmatrix}}$&$\mathsm{ \begin{bmatrix} -1 & 0 \\ 0 & 1 \end{bmatrix}}$& $\hat{y}_1$ & $\hat{y}_2$\\
    \midrule
    $t_1$ &$\mathsm{ \begin{bmatrix} 0 & -2 \\ 1 & 0 \end{bmatrix}}$&$\mathsm{ \begin{bmatrix} -1 & 2 \\ 0 & 1 \end{bmatrix}}$&$\mathsm{ \begin{bmatrix} -1 & 0 \\ 1 & 1 \end{bmatrix}}$&$\mathsm{ \begin{bmatrix} 1 & 0 \\ 0 & 1 \end{bmatrix}}$&$\mathsm{ \begin{bmatrix} 1 & 0 \\ 0 & 1 \end{bmatrix}}$ & $\frac{1}{\hat{y}_1}$ & $\frac{\hat{y}_1^2 \hat{y}_2}{(1+\hat{y}_1)^2}$\\
    \midrule
    $t_2$ &$\mathsm{ \begin{bmatrix} 0 & 2 \\ -1 & 0 \end{bmatrix}}$&$\mathsm{ \begin{bmatrix} 1 & -2 \\ 1 & -1 \end{bmatrix}}$&$\mathsm{ \begin{bmatrix} -1 & -2 \\ 1 & 1 \end{bmatrix}}$&$\mathsm{ \begin{bmatrix} 1 & 0 \\ 1 & -1 \end{bmatrix}}$&$\mathsm{ \begin{bmatrix} 1 & 2 \\ 0 & -1 \end{bmatrix}}$ & $\frac{\hat{y}_1 \hat{y}_2}{1+2\hat{y}_1 +\hat{y}_1^2 + \hat{y}_1^2 \hat{y}_2}$ & $\frac{(1+\hat{y}_1)^2}{\hat{y}_1^2 \hat{y}_2}$\\
    \midrule
    $t_3$&$\mathsm{ \begin{bmatrix} 0 & -2 \\ 1 & 0 \end{bmatrix}}$&$\mathsm{ \begin{bmatrix} -1 & 0 \\ -1 & 1 \end{bmatrix}}$&$\mathsm{ \begin{bmatrix} -1 & -2 \\ 0 & 1 \end{bmatrix}}$&$\mathsm{ \begin{bmatrix} -1 & 2 \\ -1 & 1 \end{bmatrix}}$&$\mathsm{ \begin{bmatrix} 1 & 2 \\ -1 & -1 \end{bmatrix}}$ & $\frac{1+2\hat{y}_1 +\hat{y}_1^2 + \hat{y}_1^2 \hat{y}_2}{\hat{y}_1 \hat{y}_2}$ & $\frac{\hat{y}_2}{(1+\hat{y}_1 + \hat{y}_1 \hat{y}_2)^2}$\\
    \midrule
    $t_4$&$\mathsm{ \begin{bmatrix} 0 & 2 \\ -1 & 0 \end{bmatrix}}$&$\mathsm{ \begin{bmatrix} -1 & 0 \\ 0 & -1 \end{bmatrix}}$&$\mathsm{ \begin{bmatrix} -1 & 0 \\ 0 & -1 \end{bmatrix}}$&$\mathsm{ \begin{bmatrix} 1 & -2 \\ 0 & -1 \end{bmatrix}}$&$\mathsm{ \begin{bmatrix} 1 & 0 \\ -1 & -1 \end{bmatrix}}$ &$\frac{1}{\hat{y}_1(1+\hat{y}_2)}$ & $\frac{(1+\hat{y}_1 + \hat{y}_1 \hat{y}_2)^2}{\hat{y}_2}$ \\
    \midrule
    $t_5$&$\mathsm{ \begin{bmatrix} 0 & -2 \\ 1 & 0 \end{bmatrix}}$&$\mathsm{ \begin{bmatrix} 1 & 0 \\ 0 & -1 \end{bmatrix}}$&$\mathsm{ \begin{bmatrix} 1 & 0 \\ 0 & -1 \end{bmatrix}}$&$\mathsm{ \begin{bmatrix} -1 & 0 \\ 0 & -1 \end{bmatrix}}$&$\mathsm{ \begin{bmatrix} -1 & 0 \\ 0 & -1 \end{bmatrix}}$ & $\hat{y}_1(1+\hat{y}_2)$ & $\frac{1}{\hat{y}_2}$ \\
    \bottomrule
  \end{tabular}
  \caption{The six inequivalent seeds for a cluster algebra of type $C_2$, showing $c$-vectors and g-vectors with respect to two choices of initial seed and the $\hat{y}$-variables in terms of those of the seed at $t_0$.}
  \label{C2examplegvecs}
\end{table}
}

{
  \renewcommand{\arraystretch}{1.2}
\begin{table}
  \centering
  \begin{tabular}{cccccccc}
    \toprule
    $t$&$B^{\rm op}_t$&$C^{\rm op}_{t,t_0}$&$G^{\rm op}_{t,t_0}$&$C^{\rm op}_{t,t_1}$&$G^{\rm op}_{t,t_1}$ & $\hat{y}_{1;t}^{\rm op}$ & $\hat{y}_{2;t}^{\rm op}$\\
    \midrule
    $t_0$&$ \biggl[\mathsm{ \begin{matrix} 0 & -2 \\ 1 & 0 \end{matrix}}\biggr]$&$\mathsm{ \begin{bmatrix} 1 & 0 \\ 0 & 1 \end{bmatrix}}$&$\mathsm{ \begin{bmatrix} 1 & 0 \\ 0 & 1 \end{bmatrix}}$&$\mathsm{ \begin{bmatrix} -1 & 2 \\ 0 & 1 \end{bmatrix}}$&$\mathsm{ \begin{bmatrix} -1 & 0 \\ 1 & 1 \end{bmatrix}}$& $\check{y}_1$ & $\check{y}_2$\\
    \midrule
    $t_1$ &$\mathsm{ \begin{bmatrix} 0 & 2 \\ -1 & 0 \end{bmatrix}}$&$\mathsm{ \begin{bmatrix} -1 & 0 \\ 0 & 1 \end{bmatrix}}$&$\mathsm{ \begin{bmatrix} -1 & 0 \\ 0 & 1 \end{bmatrix}}$&$\mathsm{ \begin{bmatrix} 1 & 0 \\ 0 & 1 \end{bmatrix}}$&$\mathsm{ \begin{bmatrix} 1 & 0 \\ 0 & 1 \end{bmatrix}}$ & $\frac{1}{\check{y}_1}$ & $(1+\check{y}_1)^2\check{y}_2$\\
    \midrule
    $t_2$ &$\mathsm{ \begin{bmatrix} 0 & -2 \\ 1 & 0 \end{bmatrix}}$&$\mathsm{ \begin{bmatrix} -1 & 0 \\ 0 & -1 \end{bmatrix}}$&$\mathsm{ \begin{bmatrix} -1 & 0 \\ 0 & -1 \end{bmatrix}}$&$\mathsm{ \begin{bmatrix} 1 & 0 \\ 0 & -1 \end{bmatrix}}$&$\mathsm{ \begin{bmatrix} 1 & 0 \\ 0 & -1 \end{bmatrix}}$ & $\frac{1+\check{y}_2  + 2\check{y}_1 \check{y}_2+ \check{y}_1^2 \check{y}_2 }{\check{y}_1}$ & $\frac{1}{(1+\check{y}_1)^2\check{y}_2}$\\
    \midrule
    $t_3$&$\mathsm{ \begin{bmatrix} 0 & 2 \\ -1 & 0 \end{bmatrix}}$&$\mathsm{ \begin{bmatrix} 1 & -2 \\ 0 & -1 \end{bmatrix}}$&$\mathsm{ \begin{bmatrix} 1 & 0 \\ -1 & -1 \end{bmatrix}}$&$\mathsm{ \begin{bmatrix} -1 & 0 \\ 0 & -1 \end{bmatrix}}$&$\mathsm{ \begin{bmatrix} -1 & 0 \\ 0 & -1 \end{bmatrix}}$ & $\frac{\check{y}_1}{1+\check{y}_2  + 2\check{y}_1 \check{y}_2+ \check{y}_1^2 \check{y}_2 }$ & $\frac{(1 + \check{y}_2 + \check{y}_1 \check{y}_2)^2}{\check{y}_1^2 \check{y}_2}$\\
    \midrule
    $t_4$&$\mathsm{ \begin{bmatrix} 0 & -2 \\ 1 & 0 \end{bmatrix}}$&$\mathsm{ \begin{bmatrix} -1 & 2 \\ -1 & 1 \end{bmatrix}}$&$\mathsm{ \begin{bmatrix} 1 & 2 \\ -1 & -1 \end{bmatrix}}$&$\mathsm{ \begin{bmatrix} -1 & 0 \\ -1 & 1 \end{bmatrix}}$&$\mathsm{ \begin{bmatrix} -1 & -2 \\ 0 & 1 \end{bmatrix}}$ &$\frac{1+\check{y}_2}{\check{y}_1 \check{y}_2}$ & $\frac{\check{y}_1^2 \check{y}_2}{(1 + \check{y}_2 + \check{y}_1 \check{y}_2)^2}$ \\
    \midrule
    $t_5$&$\mathsm{ \begin{bmatrix} 0 & 2 \\ -1 & 0 \end{bmatrix}}$&$\mathsm{ \begin{bmatrix} 1 & 0 \\ 1 & -1 \end{bmatrix}}$&$\mathsm{ \begin{bmatrix} 1 & 2 \\ 0 & -1 \end{bmatrix}}$&$\mathsm{ \begin{bmatrix} 1 & -2 \\ 1 & -1 \end{bmatrix}}$&$\mathsm{ \begin{bmatrix} -1 & -2 \\ 1 & 1 \end{bmatrix}}$ & $\frac{\check{y}_1 \check{y}_2}{1+\check{y}_2}$ & $\frac{1}{\check{y}_2}$ \\
    \bottomrule
  \end{tabular}
  \caption{The same data given in Table \ref{C2examplegvecs} but for a cluster algebra with the opposite exchange matrices. We use the notation $\check{y}_i$ to denote the corresponding $\hat{y}$ variables.}
  \label{C2opexamplegvecs}
\end{table}
}

{
  \renewcommand{\arraystretch}{1.2}
\begin{table}
  \centering
  \begin{tabular}{cccccccc}
    \toprule
    $t$&$B^\vee_t$&$C^\vee_{t,t_0}$&$G^\vee_{t,t_0}$&$C^\vee_{t,t_1}$&$G^\vee_{t,t_1}$ & $\hat{y}^\vee_{1;t}$ & $\hat{y}^\vee_{2;t}$\\
    \midrule
    $t_0$&$ \biggl[\mathsm{ \begin{matrix} 0 & 1 \\ -2 & 0 \end{matrix}}\biggr]$&$\mathsm{ \begin{bmatrix} 1 & 0 \\ 0 & 1 \end{bmatrix}}$&$\mathsm{ \begin{bmatrix} 1 & 0 \\ 0 & 1 \end{bmatrix}}$&$\mathsm{ \begin{bmatrix} -1 & 0 \\ 0 & 1 \end{bmatrix}}$&$\mathsm{ \begin{bmatrix} -1 & 0 \\ 0 & 1 \end{bmatrix}}$& $\hat{y}^\vee_1$ & $\hat{y}^\vee_2$\\
    \midrule
    $t_1$ &$\mathsm{ \begin{bmatrix} 0 & -1 \\ 2 & 0 \end{bmatrix}}$&$\mathsm{ \begin{bmatrix} -1 & 1 \\ 0 & 1 \end{bmatrix}}$&$\mathsm{ \begin{bmatrix} -1 & 0 \\ 2 & 1 \end{bmatrix}}$&$\mathsm{ \begin{bmatrix} 1 & 0 \\ 0 & 1 \end{bmatrix}}$&$\mathsm{ \begin{bmatrix} 1 & 0 \\ 0 & 1 \end{bmatrix}}$ & $\frac{1}{\hat{y}^\vee_1}$ & $\frac{\hat{y}^\vee_1 \hat{y}^\vee_2}{1+\hat{y}^\vee_1}$\\
    \midrule
    $t_2$ &$\mathsm{ \begin{bmatrix} 0 & 1 \\ -2 & 0 \end{bmatrix}}$&$\mathsm{ \begin{bmatrix} 1 & -1 \\ 2 & -1 \end{bmatrix}}$&$\mathsm{ \begin{bmatrix} -1 & -1 \\ 2 & 1 \end{bmatrix}}$&$\mathsm{ \begin{bmatrix} 1 & 0 \\ 2 & -1 \end{bmatrix}}$&$\mathsm{ \begin{bmatrix} 1 & 1 \\ 0 & -1 \end{bmatrix}}$ & $\frac{\hat{y}^\vee_1 (\hat{y}^\vee_2)^2}{(1+\hat{y}^\vee_1 +\hat{y}^\vee_1 \hat{y}^\vee_2)^2}$ & $\frac{1+\hat{y}^\vee_1}{\hat{y}^\vee_1 \hat{y}^\vee_2}$\\
    \midrule
    $t_3$&$\mathsm{ \begin{bmatrix} 0 & -1 \\ 2 & 0 \end{bmatrix}}$&$\mathsm{ \begin{bmatrix} -1 & 0 \\ -2 & 1 \end{bmatrix}}$&$\mathsm{ \begin{bmatrix} -1 & -1 \\ 0 & 1 \end{bmatrix}}$&$\mathsm{ \begin{bmatrix} -1 & 1 \\ -2 & 1 \end{bmatrix}}$&$\mathsm{ \begin{bmatrix} 1 & 1 \\ -2 & -1 \end{bmatrix}}$ & $\frac{(1+\hat{y}^\vee_1 +\hat{y}^\vee_1 \hat{y}^\vee_2)^2}{\hat{y}^\vee_1 (\hat{y}^\vee_2)^2}$ & $\frac{\hat{y}^\vee_2}{1+\hat{y}^\vee_1 + 2 \hat{y}^\vee_1 \hat{y}^\vee_2 + \hat{y}^\vee_1(\hat{y}^\vee_2)^2}$\\
    \midrule
    $t_4$&$\mathsm{ \begin{bmatrix} 0 & 1 \\ -2 & 0 \end{bmatrix}}$&$\mathsm{ \begin{bmatrix} -1 & 0 \\ 0 & -1 \end{bmatrix}}$&$\mathsm{ \begin{bmatrix} -1 & 0 \\ 0 & -1 \end{bmatrix}}$&$\mathsm{ \begin{bmatrix} 1 & -1 \\ 0 & -1 \end{bmatrix}}$&$\mathsm{ \begin{bmatrix} 1 & 0 \\ -2 & -1 \end{bmatrix}}$ &$\frac{1}{\hat{y}^\vee_1(1+\hat{y}^\vee_2)^2}$ & $\frac{1+\hat{y}^\vee_1 + 2 \hat{y}^\vee_1 \hat{y}^\vee_2 + \hat{y}^\vee_1(\hat{y}^\vee_2)^2}{\hat{y}^\vee_2}$ \\
    \midrule
    $t_5$&$\mathsm{ \begin{bmatrix} 0 & -1 \\ 2 & 0 \end{bmatrix}}$&$\mathsm{ \begin{bmatrix} 1 & 0 \\ 0 & -1 \end{bmatrix}}$&$\mathsm{ \begin{bmatrix} 1 & 0 \\ 0 & -1 \end{bmatrix}}$&$\mathsm{ \begin{bmatrix} -1 & 0 \\ 0 & -1 \end{bmatrix}}$&$\mathsm{ \begin{bmatrix} -1 & 0 \\ 0 & -1 \end{bmatrix}}$ & $\hat{y}^\vee_1(1+\hat{y}^\vee_2)^2$ & $\frac{1}{\hat{y}^\vee_2}$ \\
    \bottomrule
  \end{tabular}
  \caption{The inequivalent seeds for a cluster algebra of type $B_2$, dual to that shown in Table \ref{C2examplegvecs}.}
  \label{B2examplegvecs}
\end{table}
}

{
  \renewcommand{\arraystretch}{1.2}
\begin{table}
  \centering
  \begin{tabular}{cccccccc}
    \toprule
    $t$&$B^{{\rm op}, \vee}_t$&$C^{{\rm op}, \vee}_{t,t_0}$&$G^{{\rm op}, \vee}_{t,t_0}$&$C^{{\rm op}, \vee}_{t,t_1}$&$G^{{\rm op}, \vee}_{t,t_1}$ & $\hat{y}^{{\rm op}, \vee}_{1;t}$ & $\hat{y}^{{\rm op}, \vee}_{2;t}$\\
    \midrule
    $t_0$&$ \biggl[\mathsm{ \begin{matrix} 0 & -1 \\ 2 & 0 \end{matrix}}\biggr]$&$\mathsm{ \begin{bmatrix} 1 & 0 \\ 0 & 1 \end{bmatrix}}$&$\mathsm{ \begin{bmatrix} 1 & 0 \\ 0 & 1 \end{bmatrix}}$&$\mathsm{ \begin{bmatrix} -1 & 1 \\ 0 & 1 \end{bmatrix}}$&$\mathsm{ \begin{bmatrix} -1 & 0 \\ 2 & 1 \end{bmatrix}}$& $\check{y}^\vee_1$ & $\check{y}^\vee_2$\\
    \midrule
    $t_1$ &$\mathsm{ \begin{bmatrix} 0 & 1 \\ -2 & 0 \end{bmatrix}}$&$\mathsm{ \begin{bmatrix} -1 & 0 \\ 0 & 1 \end{bmatrix}}$&$\mathsm{ \begin{bmatrix} -1 & 0 \\ 0 & 1 \end{bmatrix}}$&$\mathsm{ \begin{bmatrix} 1 & 0 \\ 0 & 1 \end{bmatrix}}$&$\mathsm{ \begin{bmatrix} 1 & 0 \\ 0 & 1 \end{bmatrix}}$ & $\frac{1}{\check{y}^\vee_1}$ & $(1+\check{y}^\vee_1) \check{y}^\vee_2$\\
    \midrule
    $t_2$ &$\mathsm{ \begin{bmatrix} 0 & -1 \\ 2 & 0 \end{bmatrix}}$&$\mathsm{ \begin{bmatrix} -1 & 0 \\ 0 & -1 \end{bmatrix}}$&$\mathsm{ \begin{bmatrix} -1 & 0 \\ 0 & -1 \end{bmatrix}}$&$\mathsm{ \begin{bmatrix} 1 & 0 \\ 0 & -1 \end{bmatrix}}$&$\mathsm{ \begin{bmatrix} 1 & 0 \\ 0 & -1 \end{bmatrix}}$ & $\frac{(1+ \check{y}^\vee_2 + \check{y}^\vee_1 \check{y}^\vee_2)^2}{\check{y}^\vee_1}$ & $\frac{1}{(1+\check{y}^\vee_1) \check{y}^\vee_2}$\\
    \midrule
    $t_3$&$\mathsm{ \begin{bmatrix} 0 & 1 \\ -2 & 0 \end{bmatrix}}$&$\mathsm{ \begin{bmatrix} 1 & -1 \\ 0 & -1 \end{bmatrix}}$&$\mathsm{ \begin{bmatrix} 1 & 0 \\ -2 & -1 \end{bmatrix}}$&$\mathsm{ \begin{bmatrix} -1 & 0 \\ 0 & -1 \end{bmatrix}}$&$\mathsm{ \begin{bmatrix} -1 & 0 \\ 0 & -1 \end{bmatrix}}$ & $\frac{\check{y}^\vee_1}{(1+ \check{y}^\vee_2 + \check{y}^\vee_1 \check{y}^\vee_2)^2}$ & $\frac{1+ 2 \check{y}^\vee_2 + (\check{y}^\vee_2)^2 + \check{y}^\vee_1 (\check{y}^\vee_2)^2}{\check{y}^\vee_1 \check{y}^\vee_2}$\\
    \midrule
    $t_4$&$\mathsm{ \begin{bmatrix} 0 & -1 \\ 2 & 0 \end{bmatrix}}$&$\mathsm{ \begin{bmatrix} -1 & 1 \\ -2 & 1 \end{bmatrix}}$&$\mathsm{ \begin{bmatrix} 1 & 1 \\ -2 & -1 \end{bmatrix}}$&$\mathsm{ \begin{bmatrix} -1 & 0 \\ -2 & 1 \end{bmatrix}}$&$\mathsm{ \begin{bmatrix} -1 & -1 \\ 0 & 1 \end{bmatrix}}$ &$\frac{(1+\check{y}^\vee_2)^2}{\check{y}^\vee_1 (\check{y}^\vee_2)^2}$ & $\frac{\check{y}^\vee_1 \check{y}^\vee_2}{1+ 2 \check{y}^\vee_2 + (\check{y}^\vee_2)^2 + \check{y}^\vee_1 (\check{y}^\vee_2)^2}$ \\
    \midrule
    $t_5$&$\mathsm{ \begin{bmatrix} 0 & 1 \\ -2 & 0 \end{bmatrix}}$&$\mathsm{ \begin{bmatrix} 1 & 0 \\ 2 & -1 \end{bmatrix}}$&$\mathsm{ \begin{bmatrix} 1 & 1 \\ 0 & -1 \end{bmatrix}}$&$\mathsm{ \begin{bmatrix} 1 & -1 \\ 2 & -1 \end{bmatrix}}$&$\mathsm{ \begin{bmatrix} -1 & -1 \\ 2 & 1 \end{bmatrix}}$ & $\frac{\check{y}^\vee_1 (\check{y}^\vee_2)^2}{(1+\check{y}^\vee_2)^2}$ & $\frac{1}{\hat{y}^\vee_2}$ \\
    \bottomrule
  \end{tabular}
  \caption{The same data given in Table \ref{B2examplegvecs} but for a cluster algebra with the opposite exchange matrices.}
  \label{B2opexamplegvecs}
\end{table}
}

Note also that we can observe in this example the tropical dualities described in Theorem 1.2 in \cite{NZ12},
\begin{align}
(G_{t,t'})^T&=(C^\vee_{t,t'})^{-1}\,, \qquad (G^\vee_{t,t'})^T = (C_{t,t'})^{-1}\,, \\
C_{t,t'} &= (C^{\rm op}_{t',t})^{-1}\,, \,\,\, \quad \qquad C^\vee_{t,t'} = (C^{{\rm op}, \vee}_{t',t})^{-1}\,.
\notag 
\end{align}

\end{example}

\subsection{Quasi-homomorphisms of cluster algebras}
\label{Quasihom}

We first recall Fraser's definition of quasi-homomorphisms of cluster algebras \cite{Fra16,Fra20}. 

Let $\mathcal{A}$ be a cluster algebra of geometric type. For elements $x,y\in \mathcal{A}$, $x$ is said to be proportional to $y$ (written $x \propto y$) if $x = m y$ for some Laurent monomial $m$ in frozen variables. 

Let $\mathcal{A}$ and $\mathcal{A}'$ be two cluster algebras of geometric type with $n$ active cluster variables in each seed, and with respective groups $\mathbb{P}$ and $\mathbb{P}'$ of Laurent monomials in frozen variables. 
%Let there be an isomorphism $t \mapsto t'$ from the labelled tree $\TT_n$ of $\mathcal{A}$ to the labelled tree $\TT_n'$ of $\mathcal{A}'$. 
An algebra homomorphism $f: \mathcal{A} \to \mathcal{A}'$ is called a quasi-homomorphism (Definition 2.1 of \cite{Fra20})
if $f(\mathbb{P}) \subset \mathbb{P}'$, and there is a labelled seed $\Sigma_t=({\bf x}_t, {\bf y}_t, B_t)$ for $\mathcal{A}$, a labelled seed $\Sigma_{t'}=({\bf x}_{t'}, {\bf y}_{t'}, B_{t'})$ for $\mathcal{A}'$, a permutation $\pi$ of the labels of $\Sigma_{t'}$ and a signature $\epsilon \in \{+1,-1\}$ such that
\begin{enumerate}[label=(\roman*)]
\item $f(x_{i;t}) \propto x_{\pi(i);t'}$, $1 \le i \le n$, 

\item $f(\hat{y}_{i;t}) = (\hat{y}_{\pi(i); t'})^\epsilon$, $1 \le i \le n$,

\item $B_t = \epsilon B_{t'}^\pi$, 
\end{enumerate} 
where $\hat{y}_{i;t}$ is defined in (\ref{eq:definition of y hat}).  
%We refer to the case of $\epsilon = 1$ as even parity and $\epsilon = -1$ as odd parity. 
Note\footnote{See the more detailed discussion in \cite{Fra16}.} that in (ii) above we consider the natural extension of $f$ to act on (the semifield of subtraction free) elements of $\mathcal{F}$.
Note also that properties (i) and (ii) straightforwardly imply property (iii) but we include (iii) for convenience as we make use of it in the following. Let $\mu = \mu_{k_r} \ldots \mu_{k_1}$ be a sequence of mutations and let $\mu^\pi = \mu_{\pi(k_r)} \ldots \mu_{\pi(k_1)}$ be the permuted sequence. If properties (i), (ii) (and (iii)) hold for labelled seeds $\Sigma_t$ and $\Sigma_{t'}$, then they also hold for labelled seeds $\Sigma_{\tilde{t}} = \mu \Sigma_t$ and $\Sigma_{\tilde{t}'} = \mu^\pi \Sigma_{t'}$. In other words, the seeds of $\mathcal{A}'$ mutate `in parallel' with those of $\mathcal{A}$ \cite{Fra16,Fra20}. 

Note that for any cluster algebra of geometric type  there is a canonical quasi-isomorphism ${\bf op}: \mathcal{A} \rightarrow \mathcal{A}^{\rm op}$ with $\pi = {\rm id}$ and $\epsilon=-1$ which maps seeds $\Sigma_t$ for $\mathcal{A}$ to seeds $\Sigma_t^{\rm op}$ for $\mathcal{A}^{\rm op}$ via ${\bf op}(x_{i;t}) = x_{i;t}^{\rm op}$.

Two quasi-homomorphisms $f$ and $g$ are said to be proportional if $f(x_{i;t}) \propto g(x_{i;t})$ for $1 \leq i \leq n$. A quasi-homomorphism $f : \mathcal{A} \rightarrow \mathcal{A}'$ is said to be a quasi-isomorphism if there exists a quasi-homomorphism $g : \mathcal{A}' \rightarrow \mathcal{A}$ such that the composition $g \circ f$ is proportional to the identity map on $\mathcal{A}$ and we call $g$ a quasi-inverse for $f$. As described in \cite{Fra16}, for cluster algebras of geometric type, the existence of a quasi-homomorphism relating the labelled seeds $\Sigma_t$ and $\Sigma_{t'}$ is a matter of linear algebra. One needs, in addition to property (iii) above, that the integer row span of the extended exchange matrix $\tilde{B}_t$ contains that of $\tilde{B}_{t'}$. A quasi-inverse will then exist if and only if the row spans agree. Here we will be interested in the relations between g-vectors for the corresponding cluster algebras. Since g-vectors are determined by the properties of cluster algebras with principal coefficients, the relevant integer row spans are complete. In the following we will therefore restrict our attention to quasi-isomorphisms.
A quasi-homomorphism from a cluster algebra to itself necessarily has a quasi-inverse \cite{Fra16} and we call it a quasi-automorphism. Note that this definition of quasi-automorphism generalises the notion of cluster automorphism as introduced in \cite{ASS12}.

\subsection{Tropicalisation of quasi-isomorphisms}

Let $\mathcal{A}$ be a cluster algebra of geometric type with skew-symmetrisable exchange matrix. Take any quasi-isomorphism $f: \mathcal{A} \to \mathcal{A}'$ of signature $\epsilon=1$. We have that $f$ relates a labelled seed $\Sigma_{t} = ({\bf x}_t, {\bf y}_t, B_t)$ for $\mathcal{A}$ to a labelled seed $\Sigma_{t'} = ({\bf x}_{t'}, {\bf y}_{t'}, B_{t'})$ for $\mathcal{A}'$ via a permutation $\pi$ such that $f(x_{i;t}) \propto x_{\pi(i);t'}$ and $B_t=B_{t'}^\pi$.

Consider an arbitrary cluster variable $x \in \mathcal{A}$. We can find some labelled seed $\Sigma_{\tilde{t}} = ({\bf x}_{\tilde{t}}, {\bf y}_{\tilde{t}}, B_{\tilde{t}})$ for $\mathcal{A}$ and some $i \in \{1,\ldots,n\}$ such that $x = x_{i;\tilde{t}} \in {\bf x}_{\tilde{t}}$. The labelled seed $\Sigma_{\tilde{t}}$ can be obtained by some sequence of mutations $\mu = \mu_{l_m} \ldots \mu_{l_1}$  from $\Sigma_{t}$, i.e. $\Sigma_{\tilde{t}} = \mu \Sigma_{t}$ (or $\tilde{t} = \mu t$). %For brevity we use the notation $x=x_{i,\tilde{t}'} = \tilde{\mu} (x_{i,t'})$. 
Now define the permuted sequence of mutations $\mu^\pi = \mu_{\pi(l_m)}\ldots \mu_{\pi(l_1)}$ and define a corresponding labelled seed for $\mathcal{A}'$ via $\Sigma_{\tilde{t}'} = \mu^\pi \Sigma_{t'}$ (or $\tilde{t}' = \mu^\pi t'$). 

Define also a map $h_f$ as the composition of $f$ with the map which sends all frozen variables in $\mathcal{A}'$ to one. The map $h_f$ sends active cluster variables to active cluster variables, e.g. $h_f(x_{i;t}) = x_{\pi(i);t'}$. Then we have
\begin{equation}
    h_f(x) = h_f(x_{i;\tilde{t}}) = x_{\pi(i);\tilde{t}'}.
\end{equation}

Now consider the two sets of matrices $G_{s,t}$ and $G_{s',t'}$ for arbitrary vertices $s$, $s'$.
%(associated to the cluster algebras with principal coefficients $\mathcal{A}_{\bullet}(B_{t})$ and $\mathcal{A}_{\bullet}(B_{t'})$ respectively)
The initial data for these sets are given by
\begin{equation}
   C_{t,t} = G_{t,t} = 1\!\!1_{n \times n} = C_{t',t'} = G_{t',t'} \,.
   %\qquad  {\bf g}(x_{i;t};t) = {\bf g}(x_{i;t'};t') = {\bf e}_i\,.
\end{equation}
The remaining g-vectors in each case are obtained by mutation following the rule (\ref{gvecmutation}). Note that since $B_t = B_{t'}^\pi$ we will find 
\begin{equation}
C_{\mu t, t} = C_{\mu^\pi t' ,t'}^\pi \,, \qquad G_{\mu t, t} = G_{\mu^\pi t' ,t'}^\pi\,,
\label{gvecsauto}
\end{equation}
for any sequences of mutations $\mu = \mu_{l_m} \ldots \mu_{l_1}$ and $\mu^\pi = \mu_{\pi(l_m)} \ldots \mu_{\pi(l_1)}$. It follows that we have
\begin{equation}
G^\pi_{\tilde{t}',t'} = G^\pi_{\mu^\pi t' ,t'} = G_{\mu t , t} = G_{\tilde{t},t}\,.
\end{equation}
Now, considering column $i$ on both sides, we have
\begin{equation}\label{gmaps}
{\bf g}(x_{\pi(i);\tilde{t}'};t')^\pi = {\bf g}(x_{i;\tilde{t}};t)\,,
\end{equation}
where ${\bf g}(x_{\pi(i);\tilde{t}'};t')^\pi$ is the vector whose $j$th element is element $\pi(j)$ of ${\bf g}(x_{\pi(i);\tilde{t}'};t')$. 

Let $\Sigma_t$ be obtained from some choice of initial labelled seed $\Sigma_{t_0}$ for $\mathcal{A}$ by a sequence of mutations $\mu_0$ ($\mu_0 t_0 = t$). The relevant seeds are described in Fig. \ref{quasiisofig}.
\begin{figure}[ht]
\scalebox{1}{
\begin{tikzpicture}[scale=1]
    \node at (0,0) (A) {$\mathcal{A}:$};
    \node at (2,0) (t0) {$\Sigma_{t_0}$};
    \node at (4,0) (t) {$\Sigma_{t}$};
    \node at (6,0) (tt) {$\Sigma_{\tilde{t}}$};
    \node at (0,-2) (A') {$\mathcal{A}':$};
    \node at (4,-2) (t') {$\Sigma_{t'}$};
    \node at (6,-2) (tt') {$\Sigma_{\tilde{t}'}$};
    \node at (3,0.3) (mu0) {$\mu_0$};
    \node at (5,0.3) (mu0) {$\mu$};
    \node at (3.7,-1) () {$h_f$};
    \node at (5.7,-1) () {$h_f$};
    \node at (5,-1.7) (mus) {$\mu^\pi$};
   
    \draw[->] (t0)--(t);
    \draw[->] (t0)--(t);
    \draw[->] (t)--(tt);
    \draw[->] (t)--(t');
    \draw[->] (tt)--(tt');
    \draw[->] (t')--(tt');
\end{tikzpicture} }
            \caption{The relevant seeds in $\mathcal{A}$ and $\mathcal{A}'$ and the action of the quasi-isomorphism $f$.}
           \label{quasiisofig}
\end{figure}

We now apply the inverse permutation to both sides of (\ref{gmaps}), and we recall from Lemma \ref{lem:for every mutation, change of g-vector is the same as change of tropicalisation of x-coordinates} that the action on g-vectors of a sequence of mutations of the initial seed may be computed via tropicalisation of the change of dual $\hat{y}$ variables,
\begin{equation}
    {\bf g}\bigl(h_f(x);t'\bigr) = {\bf g}(x ; t)^{\pi^{-1}} = \bigl(Q^+_{t,t_0} {\bf g}(x;t_0) \bigr)^{\pi^{-1}}\,.
    \label{hfactionong}
\end{equation}

We have proved the following theorem.
\begin{theorem} \label{thm:quasi automorphism sends g vectors of cluster variables to g vectors of cluster variables}
Let $\mathcal{A}$, $\mathcal{A}'$ be cluster algebras of geometric type with skew-symmetrisable exchange matrices and let $f: \mathcal{A} \to \mathcal{A}'$ be a quasi-isomorphism of signature $\epsilon=1$. Let $\Sigma_t$ and $\Sigma_{t'}$ be labelled seeds and $\pi$ a permutation obeying properties (i), (ii) and (iii) of Sec. \ref{Quasihom}. Let $x$ be a cluster variable for $\mathcal{A}$ and let $\Sigma_{t_0}$ be some choice of initial labelled seed for $\mathcal{A}$. Then 
\begin{equation}
{\bf g}\bigl( h_f(x);t'\bigr) = \bigl(Q^+_{t,t_0} {\bf g}(x;t_0) \bigr)^{\pi^{-1}}
\end{equation}
with the piecewise linear map $Q^+_{t,t_0}$ defined in (\ref{Qmaps}).
\end{theorem}

\begin{corollary}
Suppose $f$ is as above in Theorem \ref{thm:quasi automorphism sends g vectors of cluster variables to g vectors of cluster variables}. Then the map $f^{\rm op}_{\rm op} : \mathcal{A}^{\rm op} \rightarrow \mathcal{A}^{\prime {\rm op}}$ defined by $f^{\rm op}_{\rm op} = {\bf op} \circ f \circ {\bf op}$, is also a quasi-isomorphism with $\epsilon = 1$. The analogue of relation (\ref{hfactionong}) is
\begin{equation}
    {\bf g}^{\rm op}\bigl(h_f(x);t'\bigr)= {\bf g}\bigl(h_{f^{\rm op}_{\rm op}}(x^{\rm op});t'\bigr) = {\bf g}(x^{\rm op};t)^{\pi^{-1}}  = {\bf g}^{\rm op}(x ; t)^{\pi^{-1}} = \bigl(Q^-_{t,t_0} {\bf g}^{\rm op}(x;t_0)\bigr)^{\pi^{-1}}\,.
    \label{hfopopactionong}
\end{equation}

    Now suppose that $f$ is as in Theorem \ref{thm:quasi automorphism sends g vectors of cluster variables to g vectors of cluster variables} except with $\epsilon =-1$. Then the maps $f^{\rm op} : \mathcal{A}^{\rm op} \rightarrow \mathcal{A}'$, defined by $f^{\rm op} = f \circ {\bf op}$, and $f_{\rm op} : \mathcal{A} \rightarrow \mathcal{A}^{\prime {\rm op}}$, defined by $f_{\rm op} = {\bf op} \circ f$, are both quasi-isomorphisms with $\epsilon = 1$. We have
    \begin{align}
    &{\bf g}(h_f(x);t') = {\bf g}\bigl(h_{f^{\rm op}}(x^{\rm op});t'\bigr)
    = {\bf g}^{\rm op}(x ; t)^{\pi^{-1}} = \bigl( Q^{-}_{t,t_0} {\bf g}^{\rm op}(x;t_0)\bigr)^{\pi^{-1}}\,, \\
    &{\bf g}^{\rm op}(h_f(x);t') = {\bf g}\bigl(h_{f_{\rm op}}(x);t'\bigr) = {\bf g}(x ; t)^{\pi^{-1}} = \bigl(Q^+_{t,t_0} {\bf g}(x;t_0)\bigr)^{\pi^{-1}}\,.
    \notag 
    \end{align}
\end{corollary}

\begin{corollary}\label{cor-quasiauto}
    In the case of a quasi-automorphism we can identify $\mathcal{A}'$ with $\mathcal{A}$ and the vertices $t_0$ and $t'$ in the above. The map $\mu_0$ then becomes $\mu_f^{-1}$ where $\mu_f$ is the sequence of mutations which represents the quasi-automorphism when acting on $\Sigma_t$, i.e. $\mu_f \Sigma_t = \Sigma_{t'}$. We have
    \begin{align}\label{quasiautoposeps}
        &{\bf g}\bigl(h_f(x);t'\bigr) = \bigl(Q^+_{t,t'} {\bf g}(x;t')\bigr)^{\pi^{-1}}\,,  &&{\bf g}^{\rm op}\bigl(h_f(x);t'\bigr) = \bigl(Q^-_{t,t'} {\bf g}^{\rm op}(x;t')\bigr)^{\pi^{-1}}\,, &&(\epsilon = +1)\,, \\
        \label{quasiautonegeps}
        &{\bf g}\bigl(h_f(x);t'\bigr) = \bigl(Q^-_{t,t'} {\bf g}^{\rm op}(x;t')\bigr)^{\pi^{-1}}\,,  &&{\bf g}^{\rm op}\bigl(h_f(x);t'\bigr) = \bigl(Q^+_{t,t'} {\bf g}(x;t')\bigr)^{\pi^{-1}}\,, &&(\epsilon=-1)\,.
    \end{align}

\end{corollary}

\begin{example}
Consider a cluster algebra of type $C_2$ with all coefficients $y_{i;t}$ set to one. We observe that the map $f(x_{1;t_0}) = x_{1;t_2}$, $f(x_{2;t_0}) = x_{2;t_2}$ defines a quasi-automorphism. It relates the labelled seeds $\Sigma_{t_0}$ and $\Sigma_{t_2}$ with the permutation $\pi = {\rm id}$ and signature $\epsilon=1$, in the sense of properties (i), (ii) and (iii) of Sec. \ref{Quasihom}, since we have $B_{t_0} = B_{t_2}$ and $f(\hat{y}_{1;t_0}) = \hat{y}_{1;t_2}$ and $f(\hat{y}_{2;t_0}) = \hat{y}_{2;t_2}$. Note that the quasi-automorphism maps generic cluster variables according to $f(x_{i;t_j}) = x_{i;t_{j+2}}$ with $i=1,2$ and $j$ treated modulo six and that, in this case, we have $h_f(x) = f(x)$.

Now we compute the matrices $C_{t_i,t_2}$ and $G_{t_i,t_2}$. We find, as anticipated in (\ref{gvecsauto}), that $C_{t_i,t_2} = C_{t_{i-2},t_0}$ and $G_{t_i,t_2} = G_{t_{i-2},t_0}$ with all indices treated modulo six. 

The dual $\hat{y}$ variables for the relevant labelled seeds are given in the third row of Table \ref{B2examplegvecs}. We have
\begin{equation}
\hat{y}^\vee_{1;t_2} = \hat{y}^{\vee \prime}_{1} = \frac{\hat{y}^\vee_1 (\hat{y}^\vee_2)^2}{(1+\hat{y}^\vee_1 +\hat{y}^\vee_1 \hat{y}^\vee_2)^2}\,, \qquad  \hat{y}^\vee_{2;t_2} = \hat{y}^{\vee \prime}_{2} = \frac{1+\hat{y}^\vee_1}{\hat{y}^\vee_1 \hat{y}^\vee_2}\,.
\end{equation}
Inverting we find
\begin{equation}
\hat{y}^\vee_1 = \frac{1}{\hat{y}^{\vee \prime}_{1}(1+\hat{y}^{\vee \prime}_{2})^2}\,,
\qquad \hat{y}^\vee_2 = \frac{1 + \hat{y}^{\vee \prime}_{1} + 2 \hat{y}^{\vee \prime}_{1} \hat{y}^{\vee \prime}_{2} + \hat{y}^{\vee \prime}_{1} (\hat{y}^{\vee \prime}_{2})^2}{\hat{y}^{\vee \prime}_{2}}\,.
\label{exautoyhatduals}
\end{equation}

Tropicalising the rational maps (\ref{exautoyhatduals}) with `max' convention, we find
\begin{equation}
v^\vee_1 = -v^{\vee \prime}_{1} - 2\,{\rm max}(0,v^{\vee \prime}_{2})\,,
\qquad v^\vee_2 = {\rm max}(0,v^{\vee \prime}_{1} , v^{\vee \prime}_{1} + v^{\vee \prime}_{2}, v^{\vee \prime}_{1} + 2v^{\vee \prime}_{2}) - v^{\vee \prime}_{2}\,.
\end{equation}
Now, identifying $t=t_0$ and $t' = t_2$ in eq. (\ref{quasiautoposeps}), we find that the map
\begin{align}
        Q^+_{t,t'} : \begin{pmatrix} v^{\vee \prime}_1 \\ v^{\vee \prime}_2 \end{pmatrix} \mapsto \begin{pmatrix} v^\vee_{1} \\ v^\vee_{2}\end{pmatrix}
        \label{C2examplefongvecsplus}
    \end{align}
indeed maps the columns of $G_{t_i,t_2}$ (as listed in Table \ref{C2examplegvecs}) to those of $G_{t_{i+2},t_2}$, with indices treated modulo six (and hence also maps the columns of $G_{t_i,t_0}$ to those of $G_{t_{i+2},t_0}$). We have thus verified the first equality in (\ref{quasiautoposeps}) for this example.

Tropicalising (\ref{exautoyhatduals}) with `min' convention, we find
\begin{equation}
w^\vee_1 = -w^{\vee \prime}_{1} - 2\,{\rm min}(0,w^{\vee \prime}_{2})\,,
\qquad w^\vee_2 = {\rm min}(0,w^{\vee \prime}_{1} , w^{\vee \prime}_{1} + w^{\vee \prime}_{2}, w^{\vee \prime}_{1} + 2w^{\vee \prime}_{2}) - w^{\vee \prime}_{2}\,.
\end{equation}
We then find that the map
\begin{align}
        Q^-_{t,t'} : \begin{pmatrix} w^{\vee \prime}_1 \\ w^{\vee \prime}_2 \end{pmatrix} \mapsto \begin{pmatrix} w^\vee_{1} \\ w^\vee_{2}\end{pmatrix}
        \label{C2examplefongvecsminus}
    \end{align}
maps the columns of $G^{\rm op}_{t_i,t_2}$ (as listed in Table \ref{C2opexamplegvecs}) to those of $G^{\rm op}_{t_{i+2},t_2}$, with indices treated modulo six (and hence also maps the columns of $G^{\rm op}_{t_i,t_0}$ to those of $G^{\rm op}_{t_{i+2},t_0}$). This verifies the second equality in (\ref{quasiautoposeps}) for this example.

Similarly, we may define the analogous quasi-homomorphism on a cluster algebra of type $B_2$, again with coefficients set to one. The analysis is the same as above for $C_2$, except that the relevant $\hat{y}$ variables are then given in the third row of Table \ref{C2examplegvecs},
\begin{equation}
    \hat{y}_{1;t_2} = \hat{y}_1' = \frac{\hat{y}_1 \hat{y}_2}{1+2\hat{y}_1 +\hat{y}_1^2 + \hat{y}_1^2 \hat{y}_2}\, \qquad \hat{y}_{2;t_2} = \hat{y}_2' = \frac{(1+\hat{y}_1)^2}{\hat{y}_1^2 \hat{y}_2}\,.
\end{equation}
Inverting these equations gives us
\begin{equation}
    \hat{y}_1 = \frac{1}{\hat{y}_1'(1+\hat{y}_2')}\, \qquad \hat{y}_2 = \frac{(1+\hat{y}_1' + \hat{y}_1' \hat{y}_2')^2}{\hat{y}_2'}\,.
\end{equation}
Tropicalising with `max' convention then yields
\begin{equation}
    v_1 = - v_1' - {\rm max}(0,v_2')\, \qquad v_2 = 2\, {\rm max}(0, v_1', v_1' + v_2') - v_2'\,.
\end{equation}
The map 
\begin{align}
      Q^{\vee +}_{t,t'}:  \begin{pmatrix} v^{\prime}_1 \\ v^{\prime}_2 \end{pmatrix} \mapsto \begin{pmatrix} v_{1} \\ v_{2}\end{pmatrix}
        \label{B2examplefongvecsplus}
    \end{align}
then maps the columns of $G^\vee_{t_i,t_2}$ to those of $G^\vee_{t_{i+2},t_2}$ (hence, again, those of $G^\vee_{t_i;t_0}$ to those of $G^\vee_{t_{i+2};t_0}$) with all indices understood modulo six.  

Tropicalising instead with `min' convention then yields
\begin{equation}
    w_1 = - w_1' - {\rm min}(0,w_2')\, \qquad w_2 = 2\, {\rm min}(0, w_1', w_1' + w_2') - w_2'\,.
\end{equation}
The map 
\begin{align}
       Q^{\vee -}_{t,t'}: \begin{pmatrix} w^{\prime}_1 \\ w^{\prime}_2 \end{pmatrix} \mapsto \begin{pmatrix} w_{1} \\ w_{2}\end{pmatrix}
        \label{B2examplefongvecsminus}
    \end{align}
then maps the columns of $G^{{\rm op},\vee}_{t_i,t_2}$ to those of $G^{{\rm op},\vee}_{t_{i+2},t_2}$ (and those of $G^{{\rm op},\vee}_{t_i;t_0}$ to those of $G^{{\rm op}, \vee}_{t_{i+2};t_0}$). Thus we have verified the relations (\ref{quasiautoposeps}) applied to the dual of the previous case.  
\end{example}

\begin{example}
Consider a cluster algebra of type $A_2$ with all coefficients $y_{i;t}$ set to one. We define a quasi-automorphism $f$ via $f(x_{1;t_0}) = x_{2;t_1}$ and $f(x_{2;t_0}) = x_{1;t_1}$\,. It relates $\Sigma_{t_0}$ and $\Sigma_{t_1}$ via the permutation $\pi(1) = 2$, $\pi(2)=1$ with signature $\epsilon=1$. Note that the quasi-automorphism relates the cluster variables according to $f(x_{i;t_j}) = x_{\pi(i),t_{j+1}}$ with $i=1,2$ and $j=0,\ldots,4$. Since the exchange matrices in type $A_2$ are skew-symmetric, they are equal to their duals $B = B^\vee$. Therefore we need only look at the $\hat{y}$ variables in the same cluster algebra.

The relevant $\hat{y}$ variables are given in the second row of Table \ref{A2examplegvecs},
\begin{equation}
    \hat{y}_{1;t_1} = \hat{y}_1'= \frac{1}{\hat{y}_1}\,, \qquad \hat{y}_{2;t_1} = \hat{y}_2' = \frac{\hat{y}_1 \hat{y}_2}{1+\hat{y}_1}\,.
\end{equation}
Inverting these relations we find
\begin{equation}\label{A2yhatinv}
    \hat{y}_1 = \frac{1}{\hat{y}_1'}\, \qquad \hat{y}_2 = \hat{y}_2'(1+\hat{y}_1')\,.
\end{equation}
Tropicalising with `max' convention then gives 
\begin{equation}
    v_1 = - v_1'\, \qquad v_2 = v_2' + {\rm max}(0,v_1')\,.
    \label{A2vmap}
\end{equation}
Then, identifying $t=t_0$ and $t'=t_1$ and recalling that the permutation $\pi^{-1}=\pi$ swaps the two labels, we verify from Table \ref{A2examplegvecs} that the map $\pi^{-1} \circ Q^+_{t,t'}$,
\begin{align}
        \begin{pmatrix} v^{\prime}_1 \\ v^{\prime}_2 \end{pmatrix} \overset{Q^+_{t,t'}}{\longmapsto} \begin{pmatrix} v_{1} \\ v_{2}\end{pmatrix} \overset{\pi^{-1}}{\longmapsto} \begin{pmatrix} v_{2} \\ v_{1}\end{pmatrix}
        \label{A2examplefongvecs}
    \end{align}
maps column $1$ (resp. $2$) of $G_{t_i,t_1}$ to column $2$ (resp. $1$) of $G_{t_{i+1},t_1}$ (where $i$ runs from 0 to 4).

Tropicalising (\ref{A2yhatinv}) with `min' convention gives instead 
\begin{equation}
    w_1 = - w_1'\, \qquad w_2 = w_2' + {\rm min}(0,w_1')\,.
\end{equation}
Then the map $\pi^{-1} \circ Q^-_{t,t'}$
\begin{align}
        \begin{pmatrix} w^{\prime}_1 \\ w^{\prime}_2 \end{pmatrix} \overset{Q^-_{t,t'}}{\longmapsto} \begin{pmatrix} w_{1} \\ w_{2}\end{pmatrix} \overset{\pi^{-1}}{\longmapsto} \begin{pmatrix} w_{2} \\ w_{1}\end{pmatrix}
        \label{A2opexamplefongvecs}
    \end{align}
maps column 1 (resp. 2) of $G^{\rm op}_{t_i,t_1}$ to column 2 (resp. 1) of $G^{\rm op}_{t_{i+1},t_1}$ in Table \ref{A2opexamplegvecs}.

{
  \renewcommand{\arraystretch}{1.2}
\begin{table}
  \centering
  \begin{tabular}{cccccccc}
    \toprule
    $t$&$B_t$&$C_{t,t_0}$&$G_{t,t_0}$&$C_{t,t_1}$&$G_{t,t_1}$ & $\hat{y}_{1;t}$ & $\hat{y}_{2;t}$\\
    \midrule
    $t_0$&$ \biggl[\mathsm{ \begin{matrix} 0 & 1 \\ -1 & 0 \end{matrix}}\biggr]$&$\mathsm{ \begin{bmatrix} 1 & 0 \\ 0 & 1 \end{bmatrix}}$&$\mathsm{ \begin{bmatrix} 1 & 0 \\ 0 & 1 \end{bmatrix}}$&$\mathsm{ \begin{bmatrix} -1 & 0 \\ 0 & 1 \end{bmatrix}}$&$\mathsm{ \begin{bmatrix} -1 & 0 \\ 0 & 1 \end{bmatrix}}$& $\hat{y}_1$ & $\hat{y}_2$\\
    \midrule
    $t_1$ &$\mathsm{ \begin{bmatrix} 0 & -1 \\ 1 & 0 \end{bmatrix}}$&$\mathsm{ \begin{bmatrix} -1 & 1 \\ 0 & 1 \end{bmatrix}}$&$\mathsm{ \begin{bmatrix} -1 & 0 \\ 1 & 1 \end{bmatrix}}$&$\mathsm{ \begin{bmatrix} 1 & 0 \\ 0 & 1 \end{bmatrix}}$&$\mathsm{ \begin{bmatrix} 1 & 0 \\ 0 & 1 \end{bmatrix}}$ & $\frac{1}{\hat{y}_1}$ & $\frac{\hat{y}_1 \hat{y}_2}{1+\hat{y}_1}$\\
    \midrule
    $t_2$ &$\mathsm{ \begin{bmatrix} 0 & 1 \\ -1 & 0 \end{bmatrix}}$&$\mathsm{ \begin{bmatrix} 0 & -1 \\ 1 & -1 \end{bmatrix}}$&$\mathsm{ \begin{bmatrix} -1 & -1 \\ 1 & 0 \end{bmatrix}}$&$\mathsm{ \begin{bmatrix} 1 & 0 \\ 1 & -1 \end{bmatrix}}$&$\mathsm{ \begin{bmatrix} 1 & 1 \\ 0 & -1 \end{bmatrix}}$ & $\frac{\hat{y}_2}{1+\hat{y}_1 + \hat{y}_1 \hat{y}_2}$ & $\frac{1+\hat{y}_1}{\hat{y}_1 \hat{y}_2}$\\
    \midrule
    $t_3$&$\mathsm{ \begin{bmatrix} 0 & -1 \\ 1 & 0 \end{bmatrix}}$&$\mathsm{ \begin{bmatrix} 0 & -1 \\ -1 & 0 \end{bmatrix}}$&$\mathsm{ \begin{bmatrix} 0 & -1 \\ -1 & 0 \end{bmatrix}}$&$\mathsm{ \begin{bmatrix} -1 & 1 \\ -1 & 0 \end{bmatrix}}$&$\mathsm{ \begin{bmatrix} 0 & 1 \\ -1 & -1 \end{bmatrix}}$ & $\frac{1+\hat{y}_1 + \hat{y}_1 \hat{y}_2}{\hat{y}_2}$ & $\frac{1}{\hat{y}_1(1+\hat{y}_2)}$\\
    \midrule
    $t_4$&$\mathsm{ \begin{bmatrix} 0 & 1 \\ -1 & 0 \end{bmatrix}}$&$\mathsm{ \begin{bmatrix} 0 & 1 \\ -1 & 0 \end{bmatrix}}$&$\mathsm{ \begin{bmatrix} 0 & 1 \\ -1 & 0 \end{bmatrix}}$&$\mathsm{ \begin{bmatrix} 0 & -1 \\ -1 & 0 \end{bmatrix}}$&$\mathsm{ \begin{bmatrix} 0 & -1 \\ -1 & 0 \end{bmatrix}}$ &$\frac{1}{\hat{y}_2}$ & $\hat{y}_1(1+\hat{y}_2)$ \\
    \midrule
    $t_5$&$\mathsm{ \begin{bmatrix} 0 & -1 \\ 1 & 0 \end{bmatrix}}$&$\mathsm{ \begin{bmatrix} 0 & 1 \\ 1 & 0 \end{bmatrix}}$&$\mathsm{ \begin{bmatrix} 0 & 1 \\ 1 & 0 \end{bmatrix}}$&$\mathsm{ \begin{bmatrix} 0 & -1 \\ 1 & 0 \end{bmatrix}}$&$\mathsm{ \begin{bmatrix} 0 & -1 \\ 1 & 0 \end{bmatrix}}$ & $\hat{y}_2$ & $\hat{y}_1$ \\
    \bottomrule
  \end{tabular}
  \caption{The five inequivalent seeds (at vertices $t_0$ to $t_4$) for a cluster algebra of type $A_2$, showing $c$-vectors and g-vectors with respect to two choices of initial seed and the $\hat{y}$-variables in terms of those of the seed at $t_0$. Note the labelled seed at vertex $t_5$ is a permutation of that at $t_0$.}
  \label{A2examplegvecs}
\end{table}
}

{
  \renewcommand{\arraystretch}{1.2}
\begin{table}
  \centering
  \begin{tabular}{cccccccc}
    \toprule
    $t$&$B^{\rm op}_t$&$C^{\rm op}_{t,t_0}$&$G^{\rm op}_{t,t_0}$&$C^{\rm op}_{t,t_1}$&$G^{\rm op}_{t,t_1}$ & $\hat{y}^{\rm op}_{1;t}$ & $\hat{y}^{\rm op}_{2;t}$\\
    \midrule
    $t_0$&$ \biggl[\mathsm{ \begin{matrix} 0 & -1 \\ 1 & 0 \end{matrix}}\biggr]$&$\mathsm{ \begin{bmatrix} 1 & 0 \\ 0 & 1 \end{bmatrix}}$&$\mathsm{ \begin{bmatrix} 1 & 0 \\ 0 & 1 \end{bmatrix}}$&$\mathsm{ \begin{bmatrix} -1 & 1 \\ 0 & 1 \end{bmatrix}}$&$\mathsm{ \begin{bmatrix} -1 & 0 \\ 1 & 1 \end{bmatrix}}$& $\check{y}_1$ & $\check{y}_2$\\
    \midrule
    $t_1$ &$\mathsm{ \begin{bmatrix} 0 & 1 \\ -1 & 0 \end{bmatrix}}$&$\mathsm{ \begin{bmatrix} -1 & 0 \\ 0 & 1 \end{bmatrix}}$&$\mathsm{ \begin{bmatrix} -1 & 0 \\ 0 & 1 \end{bmatrix}}$&$\mathsm{ \begin{bmatrix} 1 & 0 \\ 0 & 1 \end{bmatrix}}$&$\mathsm{ \begin{bmatrix} 1 & 0 \\ 0 & 1 \end{bmatrix}}$ & $\frac{1}{\check{y}_1}$ & $(1+\check{y}_1)\check{y}_2$ \\
    \midrule
    $t_2$ &$\mathsm{ \begin{bmatrix} 0 & -1 \\ 1 & 0 \end{bmatrix}}$&$\mathsm{ \begin{bmatrix} -1 & 0 \\ 0 & -1 \end{bmatrix}}$&$\mathsm{ \begin{bmatrix} -1 & 0 \\ 0 & -1 \end{bmatrix}}$&$\mathsm{ \begin{bmatrix} 1 & 0 \\ 0 & -1 \end{bmatrix}}$&$\mathsm{ \begin{bmatrix} 1 & 0 \\ 0 & -1 \end{bmatrix}}$ & $\frac{1+\check{y}_2 + \check{y}_1 \check{y}_2}{\check{y}_1}$ & $\frac{1}{(1+\check{y}_1)\check{y}_2}$\\
    \midrule
    $t_3$&$\mathsm{ \begin{bmatrix} 0 & 1 \\ -1 & 0 \end{bmatrix}}$&$\mathsm{ \begin{bmatrix} 1 & -1 \\ 0 & -1 \end{bmatrix}}$&$\mathsm{ \begin{bmatrix} 1 & 0 \\ -1 & -1 \end{bmatrix}}$&$\mathsm{ \begin{bmatrix} -1 & 0 \\ 0 & -1 \end{bmatrix}}$&$\mathsm{ \begin{bmatrix} -1 & 0 \\ 0 & -1 \end{bmatrix}}$ & $\frac{\check{y}_1}{1+\check{y}_2 + \check{y}_1 \check{y}_2}$ & $\frac{1 + \check{y}_2}{\check{y}_1\check{y}_2}$\\
    \midrule
    $t_4$&$\mathsm{ \begin{bmatrix} 0 & -1 \\ 1 & 0 \end{bmatrix}}$&$\mathsm{ \begin{bmatrix} 0 & 1 \\ -1 & 1 \end{bmatrix}}$&$\mathsm{ \begin{bmatrix} 1 & 1 \\ -1 & 0 \end{bmatrix}}$&$\mathsm{ \begin{bmatrix} -1 & 0 \\ -1 & 1 \end{bmatrix}}$&$\mathsm{ \begin{bmatrix} -1 & -1 \\ 0 & 1 \end{bmatrix}}$ &$\frac{1}{\check{y}_2}$ & $\frac{\check{y}_1\check{y}_2}{1 + \check{y}_2}$ \\
    \midrule
    $t_5$&$\mathsm{ \begin{bmatrix} 0 & 1 \\ -1 & 0 \end{bmatrix}}$&$\mathsm{ \begin{bmatrix} 0 & 1 \\ 1 & 0 \end{bmatrix}}$&$\mathsm{ \begin{bmatrix} 0 & 1 \\ 1 & 0 \end{bmatrix}}$&$\mathsm{ \begin{bmatrix} 1 & -1 \\ 1 & 0 \end{bmatrix}}$&$\mathsm{ \begin{bmatrix} 0 & -1 \\ 1 & 1 \end{bmatrix}}$ & $\check{y}_2$ & $\check{y}_1$ \\
    \bottomrule
  \end{tabular}
  \caption{The same data as given in Table \ref{A2examplegvecs} but for a cluster algebra with the opposite exchange matrices.}
  \label{A2opexamplegvecs}
\end{table}
}
\end{example}

\begin{example}
Consider again a cluster algebra of type $A_2$ with all coefficients $y_{i;t}$ set to one. 
We now define a quasi-automorphism $f$ via $f(x_{1;t_0}) = x_{1;t_1}$ and $f(x_{2;t_0}) = x_{2;t_1}$\,. It relates $\Sigma_{t_0}$ and $\Sigma_{t_1}$ via the permutation $\pi = {\rm id}$ and with signature $\epsilon=-1$. If we define the five distinct cluster variables as follows,
\begin{equation}
 x_1 = x_{1;t_0}\,, \quad x_2 = x_{2;t_0} \,, \quad x_3 = x_{1;t_1}\,, \quad x_4 = x_{2;t_2}\,, \quad x_5 = x_{1;t_3}\,
\end{equation}
then we have
\begin{equation}
 f(x_1) = x_3\,, \quad f(x_2) = x_2 \,, \quad f(x_3) = x_1\,, \quad 
 f(x_4) = x_5\,, \quad f(x_5) = x_4\,.
\end{equation}
Reading from Tables \ref{A2examplegvecs} and \ref{A2opexamplegvecs}, we note here the g-vectors (with respect to vertex $t_1$) associated to the five cluster variables and their opposites,
\begin{align}
    \{{\bf g}(x_1;t_1), \ldots {\bf g}(x_5;t_1)\} &= \left\{ \begin{bmatrix} -1 \\ 0 \end{bmatrix}, \begin{bmatrix} 0 \\ 1 \end{bmatrix}, \begin{bmatrix} 1 \\ 0 \end{bmatrix}, \begin{bmatrix} 1 \\ -1 \end{bmatrix}, \begin{bmatrix} 0 \\ -1 \end{bmatrix}  \right\}\,,\\
\{{\bf g}^{\rm op}(x_1;t_1), \ldots {\bf g}^{\rm op}(x_5;t_1)\} &= \left\{ \begin{bmatrix} -1 \\ 1 \end{bmatrix}, \begin{bmatrix} 0 \\ 1 \end{bmatrix}, \begin{bmatrix} 1 \\ 0 \end{bmatrix}, \begin{bmatrix} 0 \\ -1 \end{bmatrix}, \begin{bmatrix} -1 \\ 0 \end{bmatrix}  \right\}\,.\notag 
\end{align}

We again identify $t=t_0$ and $t'=t_1$. The relevant $\hat{y}$ variables and the corresponding maps $Q^\pm_{t,t'}$ are the same as in the previous example. We note the following relations which verify (\ref{quasiautonegeps}) in this example,
\begin{align}
    Q^+_{t,t'} {\bf g}(x_1;t_1) &= {\bf g}^{\rm op}(x_3;t_1)\,, &&Q^-_{t,t'} {\bf g}^{\rm op} (x_1;t_1) = {\bf g}(x_3;t_1)\,,\\
    Q^+_{t,t'} {\bf g}(x_2;t_1) &= {\bf g}^{\rm op}(x_2;t_1)\,, &&Q^-_{t,t'} {\bf g}^{\rm op} (x_2;t_1) = {\bf g}(x_2;t_1)\,, \notag \\
    Q^+_{t,t'} {\bf g}(x_3;t_1) &= {\bf g}^{\rm op}(x_1;t_1)\,, &&Q^-_{t,t'} {\bf g}^{\rm op} (x_3;t_1) = {\bf g}(x_1;t_1)\,,\notag \\
    Q^+_{t,t'} {\bf g}(x_4;t_1) &= {\bf g}^{\rm op}(x_5;t_1)\,, &&Q^-_{t,t'} {\bf g}^{\rm op} (x_4;t_1) = {\bf g}(x_5;t_1)\,,\notag \\
    Q^+_{t,t'} {\bf g}(x_5;t_1) &= {\bf g}^{\rm op}(x_4;t_1)\,, &&Q^-_{t,t'} {\bf g}^{\rm op} (x_5;t_1) = {\bf g}(x_4;t_1)\,. \notag
\end{align}
\end{example}

\begin{remark} 
We expect that Theorem \ref{thm:quasi automorphism sends g vectors of cluster variables to g vectors of cluster variables} is also true for g-vectors of elements in the dual canonical bases of Grassmannian cluster algebras, and more generally is true for g-vectors of elements in triangular bases \cite{BZ14, Qin17} (or other bases with certain good properties) of cluster algebras.  
\end{remark}

\begin{remark}
Recall that for a quasi-automorphism $f: \mathcal{A} \to \mathcal{A}$ of a cluster algebra, we use $h_f$ to denote the map which maps a cluster variable $x$ to another cluster variable $h_f(x)$ by ignoring the frozen factors in $f(x)$. In all the above simple examples we actually have $h_f(x) = f(x)$. However we will later discuss the twist map and braid transformations in the context of Grassmannian cluster algebras where generically $f(x) \neq h_f(x)$.
\end{remark}

\section{Grassmannian cluster algebras} \label{sec:Grassmannian cluster algebras}

\subsection{Grassmannian cluster algebras} \label{subsec:Grassmannian cluster algebras}

For $k \le n$, the Grassmannian $\Gr(k,n)$ is the set of $k$-dimensional subspaces in an $n$-dimensional vector space. In this paper, with slight abuse of notation, we denote by $\Gr(k,n)$ the Zariski-open subset of the affine cone over the Grassmannian cut out by the non-vanishing of the frozen variables and denote by $\CC[\Gr(k,n)]$ its coordinate ring. This algebra is generated by Pl\"{u}cker coordinates 
\begin{align*}
P_{i_1, \ldots, i_{k}}, \quad 1 \leq i_1 < \cdots < i_{k} \leq n.
\end{align*}
Denote by $\CC[\Gr(k,n,\sim)]$ the quotient of $\CC[\Gr(k,n)]$ by the ideal generated by $P_{i,i+1, \ldots, k+i-1}-1$, $i \in [n-k+1]$.

It was shown by Scott \cite{Sco} that the ring $\CC[\Gr(k,n)]$ has a cluster algebra structure. Furthermore, in \cite{CDFL}, it is shown that the elements in the dual canonical basis of $\CC[\Gr(k,n,\sim)]$ (in particular, cluster monomials) are in bijection with semistandard Young tableaux in ${\rm SSYT}(k, [n],\sim)$. The elements in the dual canonical basis of $\CC[\Gr(k,n,\sim)]$ are in bijection with simple modules in a certain category of finite dimensional $U_q(\widehat{\mathfrak{sl}_k})$-modules \cite{HL10}. We now recall these results. 

The cluster algebra $\CC[\Gr(k,n)]$ has an initial seed $\Sigma_{t_0}$ with an initial quiver $Q$, where $Q$ has vertices $(0, 0)$, and $(a,b)$, $a \in [n-k]$, $b \in [k]$, and arrows $(0,0) \to (1,1)$, $(a-1,b) \to (a,b)$, $a \in [2, n-k]$, $b \in [k]$, $(a, b-1) \to (a,b)$, $a \in [n-k]$, $b \in [2, k]$, $(a+1,b+1) \to (a,b)$, $a \in [n-k-1]$, $b \in [k-1]$, see Figure \ref{fig:initial cluster for a quotient of Gr(5,10)}. The initial cluster variables are certain Pl\"{u}cker coordinates. The frozen variable at $(0,0)$ is $P_{1,\ldots,k}$. Then for $b=1,\ldots ,k$ the cluster variables in column $b$ are $P_{1,2,\ldots, k-b, k-b+2,\ldots,k+1}$, $\ldots$, $P_{1,2,\ldots, k-b, n-b+1,\ldots, n}$, with the last variable in each column being frozen, except for column $k$ where they are all frozen.

\begin{figure}
\scalebox{0.6}{
\begin{tikzpicture}[scale=0.8]
     \node at (-4,0) (v00) {\fbox{$P_{1,2,3,4,5}, (0,0)$}};
     \node at (0,-4) (v10) {$P_{1,2,3,4,6}, (1,1)$};
     \node at (0,-8) (v20) {$P_{1,2,3,4,7}, (2,1)$};
     \node at (0,-12) (v30) {$P_{1,2,3,4,8}, (3,1)$};
     \node at (0,-16) (v40) {$P_{1,2,3,4,9}, (4,1)$};
     \node at (0,-20) (v50) {\fbox{$P_{1,2,3,4,10}, (5,1)$}};
    
     \node at (4,-4) (v11) {$P_{1,2,3,5,6}, (1,2)$};
     \node at (4,-8) (v21) {$P_{1,2,3,6,7}, (2,2)$};
     \node at (4,-12) (v31) {$P_{1,2,3,7,8}, (3,2)$};
     \node at (4,-16) (v41) {$P_{1,2,3,8,9}, (4,2)$};
     \node at (4,-20) (v51) {\fbox{$P_{1,2,3,9,10}, (5,2)$}};
     
     \node at (8,-4) (v12) {$P_{1,2,4,5,6}, (1,3)$};
     \node at (8,-8) (v22) {$P_{1,2,5,6,7}, (2,3)$};
     \node at (8,-12) (v32) {$P_{1,2,6,7,8}, (3,3)$};
     \node at (8,-16) (v42) {$P_{1,2,7,8,9}, (4,3)$};
     \node at (8,-20) (v52) {\fbox{$P_{1,2,8,9,10}, (5,3)$}};
     
     \node at (12,-4) (v13) {$P_{1,3,4,5,6}, (1,4)$};
     \node at (12,-8) (v23) {$P_{1,4,5,6,7}, (2,4)$};
     \node at (12,-12) (v33) {$P_{1,5,6,7,8}, (3,4)$};
     \node at (12,-16) (v43) {$P_{1,6,7,8,9}, (4,4)$};
     \node at (12,-20) (v53) {\fbox{$P_{1,7,8,9,10}, (5,4)$}};
     
     \node at (16,-4) (v14) {\fbox{$P_{2,3,4,5,6}, (1,5)$}};
     \node at (16,-8) (v24) {\fbox{$P_{3,4,5,6,7}, (2,5)$}};
     \node at (16,-12) (v34) {\fbox{$P_{4,5,6,7,8}, (3,5)$}};
     \node at (16,-16) (v44) {\fbox{$P_{5,6,7,8,9}, (4,5)$}};
     \node at (16,-20) (v54) {\fbox{$P_{6,7,8,9,10}, (5,5)$}};
     
     \draw[thick,<-] (v10)--(v00);
     \draw[thick,<-] (v20)--(v10);
     \draw[thick,<-] (v30)--(v20);
     \draw[thick,<-] (v40)--(v30);
     \draw[thick,<-] (v50)--(v40);

     \draw[thick,<-] (v21)--(v11);
     \draw[thick,<-] (v31)--(v21);
     \draw[thick,<-] (v41)--(v31);
     \draw[thick,<-] (v51)--(v41);
     
     \draw[thick,<-] (v22)--(v12);
     \draw[thick,<-] (v32)--(v22);
     \draw[thick,<-] (v42)--(v32);
     \draw[thick,<-] (v52)--(v42);
     
     \draw[thick,<-] (v23)--(v13);
     \draw[thick,<-] (v33)--(v23);
     \draw[thick,<-] (v43)--(v33);
     \draw[thick,<-] (v53)--(v43);
     
     \draw[thick,<-] (v11)--(v10);
     \draw[thick,<-] (v12)--(v11);
     \draw[thick,<-] (v13)--(v12);
     \draw[thick,<-] (v14)--(v13);
     
     \draw[thick,<-] (v21)--(v20);
     \draw[thick,<-] (v22)--(v21);
     \draw[thick,<-] (v23)--(v22);
     \draw[thick,<-] (v24)--(v23);
     
     \draw[thick,<-] (v31)--(v30);
     \draw[thick,<-] (v32)--(v31);
     \draw[thick,<-] (v33)--(v32);
     \draw[thick,<-] (v34)--(v33);
     
     \draw[thick,<-] (v41)--(v40);
     \draw[thick,<-] (v42)--(v41);
     \draw[thick,<-] (v43)--(v42);
     \draw[thick,<-] (v44)--(v43);
     
     \draw[thick,<-] (v10)--(v21);
     \draw[thick,<-] (v21)--(v32);
     \draw[thick,<-] (v32)--(v43);
     \draw[thick,<-] (v43)--(v54);
     
     \draw[thick,<-] (v20)--(v31);
     \draw[thick,<-] (v31)--(v42);
     \draw[thick,<-] (v42)--(v53);
     
     \draw[thick,<-] (v30)--(v41);
     \draw[thick,<-] (v41)--(v52);
     
     \draw[thick,<-] (v40)--(v51);
     
     \draw[thick,<-] (v11)--(v22);
     \draw[thick,<-] (v22)--(v33);
     \draw[thick,<-] (v33)--(v44);
     
     \draw[thick,<-] (v12)--(v23);
     \draw[thick,<-] (v23)--(v34);
     
     \draw[thick,<-] (v13)--(v24);
     
\end{tikzpicture} }
            \caption{The quiver diagram for the initial seed for $\CC[\Gr(5,10)]$, where $(a,b)$'s in the vertices are used to denote the positions of cluster variables and frozen variables.}
            \label{fig:initial cluster for a quotient of Gr(5,10)}
\end{figure}
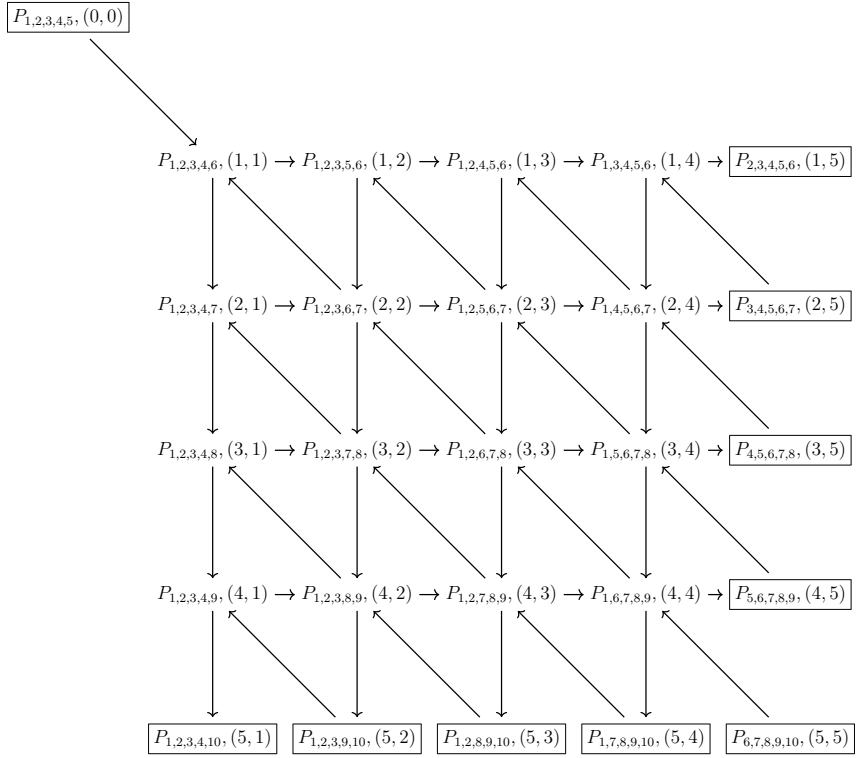

\subsection{Semistandard Young tableaux} A semistandard Young tableau is a Young tableau with weakly increasing rows and strictly increasing columns. For $k,n \in \ZZ_{\ge 1}$, we denote by ${\rm SSYT}(k, [n])$ the set of rectangular semistandard Young tableaux with $k$ rows and with entries in $[n] = \{1,\ldots, n\}$ (with arbitrarly many columns). The empty tableau is denoted by $\mathds{1}$. 

In the following, we recall a few definitions related to semistandard Young tableaux. For $S,T \in {\rm SSYT}(k, [n])$, let $S \cup T$ be the row-increasing tableau whose $i$th row is the union of the $i$th rows of $S$ and $T$ (as multisets), \cite{CDFL}. We call $S$ a factor of $T$, and write $S \subset T$, if the $i$th row of $S$ is contained in that of $T$ (as multisets), for $i \in [k]$. In this case, we define $\frac{T}{S}=S^{-1}T=TS^{-1}$ to be the row-increasing tableau whose $i$th row is obtained by removing that of $S$ from that of $T$ (as multisets), for $i \in [k]$. A tableau $T \in {\rm SSYT}(k, [n])$ is trivial if each entry of $T$ is one less than the entry below it. For any $T \in {\rm SSYT}(k, [n])$, we  denote by $T_{\text{red}} \subset T$ the semistandard tableau obtained by removing a maximal trivial factor from $T$. For trivial $T$, one has $T_{\text{red}} = \mathds{1}$. Let ``$\sim$'' be the equivalence relation on $S, T \in {\rm SSYT}(k, [n])$ defined by: $S \sim T$ if and only if $S_{\text{red}} = T_{\text{red}}$. We denote by ${\rm SSYT}(k, [n],\sim)$ the set of $\sim$-equivalence classes. A one-column tableau is called a fundamental tableau if its content is $[i,i+k] \setminus \{r\}$ for $r \in \{i+1, \ldots, i+k-1\}$. A tableau $T$ is said to have small gaps if each of its columns is a fundamental tableau. Then any tableau in $\SSYT(k,[n])$ is $\sim$-equivalent to a unique small gap semistandard tableau.

\subsection{\texorpdfstring{Dual canonical basis of $\CC[\Gr(k,n,\sim)]$}{Dual canonical basis of C[Gr(k,n,~)]}}

We now recall the explicit formula of $\ch(T)$ in the dual canonical basis of $\CC[\Gr(k,n)]$ in \cite[Theorem 5.8]{CDFL}. For that we need to first define $w_T \in S_m$, $P_{u; T'}$, $u \in S_m$, for every $T \in \SSYT(k, [n])$, where $T'$ is the unique small gap tableau which is $\sim$-equivalent to $T$, and $m$ is the number of columns of $T'$. 

Let ${\bf i} = i_1 \leq i_2 \dots \leq i_m$ be the entries in the first row of $T'$, and let $r_1,\dots,r_m$ be the elements such that the $a$th column of $T'$ has content $[i_a,i_a+k] \setminus \{r_a\}$. Let ${\bf j} = j_1 \leq j_2 \leq \dots \leq j_m$ be the elements $r_1,\dots,r_m$ written in weakly increasing order. 

For $u \in S_m$, define $P_{u;T'} \in \CC[\Gr(k,n)]$ as follows. Provided $j_a \in [i_{u(a)}, i_{u(a)}+k]$ for all $a \in [m]$, define the tableau $\alpha(u;T')$ to be the semistandard tableau whose columns have content $[i_{u(a)}, i_{u(a)}+k] \setminus \{j_a\}$ for $a \in [m]$, and define $P_{u; T'} = P_{\alpha(u;T')} \in \CC[\Gr(k,n)]$ to be the corresponding standard monomial. On the other hand, if $j_a \notin [i_{u_a}, i_{u(a)}+k]$ for some $a$, then the tableau $\alpha(u;T')$ is {\sl undefined} and $P_{u ;T'} = 0$. 

There is a unique $u \in S_m$ is of maximal length with the property that the sets $\{[i_{u(a)},i_{u(a)}+k] \setminus \{j_a \} \}_{a \in [m]}$ describe the columns of $T'$. This $u$ is denoted by $u = w_{T}$.

By \cite[Theorem 5.8]{CDFL}, the element $\widetilde{\ch}(T)$ in the dual canonical basis of $\CC[\Gr(k,n,\sim)]$ is given by 
\begin{align}\label{eq:formula of ch(T) in the quotient}
\widetilde{\ch}(T) = \widetilde{\ch}(T') = \sum_{u \in S_m} (-1)^{\ell(uw_T)} p_{uw_0, w_Tw_0}(1) P_{u; T'} \in \CC[\Gr(k,n,\sim)]
\end{align}
where $w_0$ is the longest word in $S_m$, $p_{u,v}(q)$ is a Kazhdan-Lusztig polynomial \cite{KL79}.  

The set $\{\widetilde{\ch}(T): T \in \SSYT(k,[n],\sim)\}$ is the dual canonical basis of $\CC[\Gr(k,n,\sim)]$, see Section 13 of \cite{HL10} and Theorems 3.17 and 5.8 in \cite{CDFL}. 

\subsection{\texorpdfstring{Dual canonical basis of $\CC[\Gr(k,n)]$}{Dual canonical basis of C[Gr(k,n)]}}

An explicit formula of the dual canonical basis of $\CC[\Gr(k,n)]$ is conjectured in the paragraph below Definition 5.9 in \cite{CDFL}. We will prove the conjecture in the rest of this section. 

For a one-column tableau $T$, denote by $P_T$ the Pl\"{u}cker coordinate with indices which are the entries of $T$. For $T \in {\rm SSYT}(k, [n])$ with columns $T_1,\dots,T_m$, let $P_T = P_{T_1} \cdots P_{T_m} \in \mathbb{C}[{\rm Gr}(k,n)]$.

We recall the definition of the weight of a tableau $T$, see Section 3 of \cite{CDFL}. For a Pl\"{u}cker coordinate $P_{i_1, \ldots, i_k}$, its weight is defined as $\wt(P_{i_1,\ldots, i_k}) = \sum_{j=2}^{k} (i_j-i_{j-1}-1)\omega_{k-j+1}$, where $\omega_j$'s are fundamental weights of $\mathfrak{sl}_k$. The weight of a tableau $T$ is defined as $\wt(T) = \sum_{i=1}^m \wt(P_{T_i})$, where $T_i$'s are columns of $T$. The weights of tableaux give a partial order on $\SSYT(k,[n])$: $T \le S$ if and only if $\wt(T) \le \wt(S)$. This partial order agrees with the partial order on $\SSYT(k,[n])$ defined as follows. For a tableau $T$, let ${\rm sh}(T)$ denote the shape of $T$. For $i \in [m]$, let $T[i]$ denote denote the restriction
of $T \in \SSYT(k,[n])$ to the entries in $[i]$. Then for every $i$, ${\rm sh}(T[i])$ is a Young diagram and it corresponds to a partition of an integer. For $T, S \in \SSYT(k,[n])$, $T \le S$ if and only if ${\rm sh}(T[i]) \le {\rm sh}(S[i])$ in the dominance order on partitions for every $i \in [n]$.

By standard monomial theory \cite{Ses}, elements in $\CC[\Gr(k,n)]$ are homogeneous polynomials in Pl\"{u}cker coordinates and $\{P_T: T \in \SSYT(k, [n])\}$ is a $\CC$-linear basis of $\CC[\Gr(k,n)]$. Using the isomorphism between the Grothendieck ring $K_0(\mathcal{C}_{\ell})$ and the quotient $\CC[\Gr(k,n,\sim)]$ of $\CC[\Gr(k,n)]$ (see Section \ref{subsec:correspondence between modules and tableaux}), it is shown in Sections 3 and 5 of \cite{CDFL} that the dual canonical basis of $\CC[\Gr(k,n,\sim)]$ is $\{ \widetilde{\ch}(T): T\in \SSYT(k,[n],\sim) \}$. 

The existence of the dual canonical basis for $\CC[\Gr(k,n)]$ follows from general results of Lusztig \cite{Lus90} and Kashiwara \cite{Kas91}. Let $L_n(\lambda)$ be note the irreducible polynomial representation of $GL_n(\CC)$ indexed by the weakly decreasing sequence $\lambda \in \ZZ_{\ge 1}^n$. Every representation $L_n(\lambda)$ has a (dual) canonical basis, see \cite{Lus90} and \cite{Kas91}. Let $\CC[\Gr(k,n)]_{(d)}$ denote the subspace spanned by degree $d$ monomials in Pl\"{u}cker coordinates. The action of $GL_n(C)$ on $\CC^n$ induces an action on $\CC[\Gr(k,n)]$ and $\CC[\Gr(k,n)]_{(d)} \cong L_n(d\omega_k)$ as $GL_n(\CC)$-representations, where $\omega_k$ is the $k$th fundamental weight for $GL_n(\CC)$. Therefore $\CC[\Gr(k,n)]$ inherits a (dual) canonical basis. 

The dual canonical basis of $\CC[\Gr(k,n)]$ can be also understood as a $\CC$-linear basis $B$ of $\CC[\Gr(k,n)]$ such that every element of $B$ is sent to an element in the dual canonical basis of $\CC[\Gr(k,n,\sim)]$ under the natural projection map. We explain this fact as follows. 

It is known that the transition matrix between the PBW basis and the canonical basis is upper triangular unipotent, see \cite{Lus90}. The $\CC$-linear basis $\{P_T: T \in \SSYT(k, [n])\}$ of $\CC[\Gr(k,n)]$ can be thought of as the PBW basis. Brundan \cite{Bru06} also gave a formula for the entries of the unitriangular transition matrix between the standard monomial and dual canonical bases of the irreducible polynomial representations of $U_q(\mathfrak{gl}_n)$ in terms of Kazhdan-Lusztig polynomials. Therefore every element $f$ in the dual canonical basis for $\CC[\Gr(k,n)]$ can be expressed as $f = \sum_{S} c_S P_S$, where $c_S \ne 0$ and $c_T = 1$ for the highest weight tableau $T$ appearing in $f = \sum_{S} c_S P_S$.  It follows that the dual canonical basis of $\CC[\Gr(k,n)]$ is also parametrized by $\SSYT(k,[n])$. Moreover, since the degree of every $P_S$ in the homogeneous polynomial $f = \sum_{S} c_S P_S$ is equal to the number of columns of $T$ ($T$ is the highest weight tableau among these $S$'s), we have that the polynomial in terms Pl\"{u}cker coordinates in the dual canonical basis of $\CC[\Gr(k,n)]$ corresponding to $T$ has degree equal to the number of columns of $T$. We have the following lemma.

\begin{lemma}
There is a one to one correspondence between semistandard Young tableaux in $\SSYT(k,[n])$ and elements in the dual canonical basis of $\CC[\Gr(k,n)]$. For every $T \in \SSYT(k,[n])$, the element in the dual canonical basis of $\CC[\Gr(k,n)]$ corresponding to $T$ is of the form $f_T = \sum_{S \in \SSYT(k,[n])} c_S P_S$, where $c_S \ne 0$, $c_T=1$, and $T$ is the highest weight tableau among the $S$'s. The degree of $f_T$ is equal to the number of columns of $T$.
\end{lemma}

Since $\CC[\Gr(k,n,\sim)]$ is a quotient of $\CC[\Gr(k,n)]$ by sending $P_{i,i+1, \ldots, k+i-1}$ to $1$, $i \in [n-k+1]$, we have that the dual canonical basis of $\CC[\Gr(k,n)]$ is a $\CC$-linear basis $B$ of $\CC[\Gr(k,n)]$ such that $B$ is the dual canonical basis of $\CC[\Gr(k,n,\sim)]$ up to sending $P_{i,i+1, \ldots, k+i-1}$ to $1$, $i \in [n-k+1]$. Moreover, the product of any dual canonical basis element and any $P_{i,i+1,\ldots,i+k-1}$, $i \in [n-k+1]$, is still a dual canonical basis element.
%Therefore the dual canonical basis of $\CC[\Gr(k,n)]$ is parametrised by the set $\SSYT(k, [n])$ of rectangular semistandard tableaux. 

Using the isomorphism between $K_0(\mathcal{C}_{\ell})$ and $\CC[\Gr(k,n,\sim)]$, and the highest weights of representations, we see that the tableau corresponding in $\SSYT(k,[n],\sim)$ to an element $f$ in the dual canonical basis of $\CC[\Gr(k,n,\sim)]$ can be obtained as follows. First write $f$ as $f = \sum_{S \in \SSYT(k,[n],\sim)} c_S P_{S}$, where $c_S \ne 0$, using Formula (\ref{eq:formula of ch(T) in the quotient}). Then the tableau $T$ corresponding to $f$ is the highest weight tableau $T$ appearing in $f = \sum_{S \in \SSYT(k,[n],\sim)} c_S P_{S}$. 

For any trivial tableau $S$, we define 
\begin{align} \label{eq:chT for trivial tableau}
\ch(S) = P_S = P_{S_1} \cdots P_{S_r},
\end{align}
where $S_i$'s are columns of $S$. Under the projection map $\CC[\Gr(k,n)] \to \CC[\Gr(k,n,\sim)]$, $\ch(S) = P_{S_1} \cdots P_{S_r}$ is sent to $1$. The dual canonical basis element in $\CC[\Gr(k,n)]$ corresponding to $S$ is $\ch(S)$. 

For $T \in \SSYT(k, [n])$, let $T'$ be the small gap tableau such that $T \sim T'$ and let $T'' = T' T^{-1}$. The tableau $T''$ is a trivial tableau. We have that $T' = T \cup T''$. Note that the number of columns of $T'$ is the same as the degree of $\widetilde{\ch}(T')$ in (\ref{eq:formula of ch(T) in the quotient}). We have that $\widetilde{\ch}(T')$ is the element in the dual canonical basis of $\CC[\Gr(k,n)]$ corresponding to $T'$. On the other hand, as discussed in the previous paragraphs, there is an element $f_{T}$ in the dual canonical basis of $\CC[\Gr(k,n)]$ corresponding to $T$.  
We need the following lemma.

\begin{lemma} \label{lem:product of dual canonical basis element and element corresponding to trivial tableau}
Let $T$ be a semistandard tableau corresponding to a dual canonical basis element $f_T$ in $\CC[\Gr(k,n)]$ and let $S$ be a trivial tableau. Then $f_T \ch(S)$ is a dual canonical basis element of $\CC[\Gr(k,n)]$ corresponding to $T \cup S$, where $\ch(S)$ is given in (\ref{eq:chT for trivial tableau}).
\end{lemma}

\begin{proof}
By Formula (\ref{eq:formula of ch(T) in the quotient}) and the fact that $S$ is a trivial tableau, the projection of $f_T\ch(S)$ to $\CC[\Gr(k,n,\sim)]$ is equal to $\widetilde{\ch}(T')$, where $T' \sim T$ is the unique small gap tableau corresponding to $T$. Therefore the homogeneous polynomial $f_T\ch(S)$ is in the dual canonical basis of $\CC[\Gr(k,n)]$. Since $f_T$ is in the dual canonical basis of $\CC[\Gr(k,n)]$, $f_T = \sum_{S'} c_{S'}P_{S'}$ for some tableaux $S'$'s, $c_{S'} \ne 0$, and the highest weight tableau on the right hand side is $T$. We have that $f_T\ch(S) = f_T P_S = \sum_{S'} c_{S'}P_{S'}P_S$. Since $\{P_{S''}: S'' \in \SSYT(k,[n])\}$ is a $\CC$-linear basis of $\CC[\Gr(k,n)]$, we have that $f_T\ch(S) = \sum_{S'} c_{S'}P_{S'}P_S = \sum_{S''} c_{S''}P_{S''}$. The highest weight term in $\sum_{S'} c_{S'}P_{S'}P_S$ is $P_TP_S$. Since $P_{T}P_S = P_{T \cup S} + \sum_{S'''} c_{S'''}P_{S'''}$, where $S'''$'s are some semistandard tableaux which have lower weight than $T \cup S$, we have that the highest weight term in $\sum_{S''} c_{S''}P_{S''}$ is $P_{T \cup S}$. Therefore the tableau corresponding to the dual canonical basis element $f_T\ch(S)$ is $T \cup S$. 
\end{proof}

By Lemma \ref{lem:product of dual canonical basis element and element corresponding to trivial tableau}, $f_{T}\ch(T'')$ is an element in the dual canonical basis of $\CC[\Gr(k,n)]$. Moreover, the tableau in $\SSYT(k,[n])$ corresponding to $f_{T}\ch(T'')$ is $T \cup T'' = T'$. Therefore $f_{T}\ch(T'') = \widetilde{\ch}(T')$. We must have that $\ch(T'')$ is divisible by $\widetilde{\ch}(T')$ and $f_{T} = \frac{1}{\ch(T'')} \widetilde{\ch}(T')$. We define 
\begin{align} \label{eq:chT for any tableau and for Grkn}
\ch(T) = \frac{1}{\ch(T'')} \widetilde{\ch}(T').
\end{align}
The element in the dual canonical basis of $\CC[\Gr(k,n)]$ corresponding to $T$ is $\ch(T)$. We have proved the following theorem which is conjectured in \cite{CDFL} (see the paragraph below Definition 5.9 in \cite{CDFL}).
\begin{theorem} \label{thm:dual canonical basis of CGrkn}
The dual canonical basis of $\CC[\Gr(k,n)]$ is $\{\ch(T): T \in \SSYT(k,[n])\}$.     
\end{theorem}

We have that the number of columns of $T$ is equal to the degree of $\ch(T)$ in Pl\"{u}cker coordinates. 

\section{Quantum affine algebras}
\label{sec:quantum affine algebras}

In this section, we recall results of quantum affine algebras \cite{CP94, FR98}, Hernandez-Leclerc's category $\mathcal{C}_{\ell}$, cluster algebra structure on the Grothendieck ring of $\mathcal{C}_{\ell}$ \cite{HL10}, and the parametrisation of dual canonical basis of $\CC[\Gr(k,n)]$ using semistandard Young tableaux \cite{CDFL}. Results of quantum affine algebras will be used in Section \ref{subsec:g vectors and tableaux} to recover a full g-vector (including entries corresponding to frozen variables) from a usual g-vector (without entries corresponding to frozen variables), and in Section \ref{subsec:degree of g vector and decomposition of dominant monomial} to study the degree of the polynomial in $\CC[\Gr(k,n)]$ corresponding to a g-vector.  

\subsection{Quantum affine algebras} \label{subsec:quantum affine algebras, qcharacters}

Let $\mathfrak{g}$ be a simple finite-dimensional Lie algebra and $I$ the set of vertices of the Dynkin diagram of $\mathfrak{g}$. Denote by $\{ \omega_i: i \in I \}$, $\{\alpha_i : i \in I\}$, $\{\alpha_i^{\vee} : i \in I\}$ the set of fundamental weights, the set of simple roots, the set of simple coroots, respectively. Denote by $P$ the integral weight lattice and
$P^+$ the set of dominant weights. The Cartan matrix is $C = (\alpha_j(\alpha_i^{\vee}))_{i,j \in I}$. Let $D = \diag(d_i: i \in I)$, where $d_i$'s are minimal positive integers such that $DC$ is symmetric.  

The {\sl quantum affine algebra} $U_q(\widehat{\mathfrak{g}})$ is a Hopf algebra that is a $q$-deformation of the universal enveloping algebra of $\widehat{\mathfrak{g}}$ \cite{Dri85, Jim85}. In this paper, we take $q$ to be a non-zero complex number which is not a root of unity. 

Denote by $\mathcal{P}$ the free abelian group generated by formal variables $Y_{i, a}^{\pm 1}$, $i \in I$, $a \in \CC^*$, denote by $\mathcal{P}^+$ the submonoid of $\mathcal{P}$ generated by $Y_{i, a}$, $i \in I$, $a \in \CC^*$. Let $\mathcal{C}$ denote the monoidal category of finite-dimensional representations of the quantum affine algebra $U_q(\widehat{\mathfrak{g}})$. Chari and Pressley \cite{CP95a} proved that any finite dimensional simple object in $\mathcal{C}$ is a highest $l$-weight module (denoted by $L(M)$) with a highest $l$-weight $M \in \mathcal{P}^+$. The elements in $\mathcal{P}^+$ are called dominant monomials.  

Frenkel and Reshetikhin \cite{FR98} introduced the $q$-character map which is an injective ring morphism $\chi_q$ from the Grothendieck ring of $\mathcal{C}$ to $\mathbb{Z}\mathcal{P} = \mathbb{Z}[Y_{i, a}^{\pm 1}]_{i\in I, a\in \mathbb{C}^*}$. For a $U_q(\widehat{\mathfrak{g}})$-module $V$, $\chi_q(V)$ encodes the decomposition of $V$ into common generalized eigenspaces for the action of a large commutative subalgebra of $U_q(\widehat{\mathfrak{g}})$ (the loop-Cartan subalgebra). These generalized eigenspaces are called $l$-weight spaces and generalized eigenvalues are called $l$-weights. One can identify $l$-weights with monomials in $\mathcal{P}$ \cite{FR98}. Then
the $q$-character of a $U_q(\widehat{\mathfrak{g}})$-module $V$ is given by (see \cite{FR98})
\begin{align*}
\chi_q(V) = \sum_{  m \in \mathcal{P}} \dim(V_{m}) m \in \mathbb{Z}\mathcal{P},
\end{align*}
where $V_{m}$ is the $l$-weight space with $l$-weight $m$. 

For $i \in I$, $a \in \mathbb{C}^*$, $k \in \ZZ_{\ge 1}$, the modules
\begin{align*}
X_{i,a}^{(k)} := L(Y_{i,a} Y_{i,aq^2} \cdots Y_{i,aq^{2k-2}})
\end{align*}
are called {\sl Kirillov-Reshetikhin modules}. The modules $X_{i,a}^{(1)} = L(Y_{i,a})$ are called {\sl fundamental modules}. 

A finite-dimensional $U_q(\widehat{\mathfrak{g}})$-module is called \textit{prime} if it is not isomorphic to a tensor product of two nontrivial $U_q(\widehat{\mathfrak{g}})$-modules (cf. \cite{CP97}). A simple $U_q(\widehat{\mathfrak{g}})$-module $L(M)$ is {\sl real} if $L(M) \otimes L(M)$ is simple (cf. \cite{Lec}).

\subsection{\texorpdfstring{Isomorphism of $K_0(\mathcal{C}_{\ell})$ and $\CC[\Gr(k,n,\sim)]$}{Isomorphism of K0(Cell) and C[Gr(k,n,~)]}}
\label{subsec:correspondence between modules and tableaux}

In \cite{HL10}, \cite{HL16}, Hernandez and Leclerc introduced a full subcategory $\mathcal{C}_{\ell}$ ($\ell \in \mathbb{Z}_{\geq 0}$) of $\mathcal{C}$. We recall the definition in the case of $\mathfrak{g}=\mathfrak{sl}_k$. 

Let $\mathfrak{g}=\mathfrak{sl}_k$ and $I=[1,k-1]$ be the set of vertices of the Dynkin diagram of $\mathfrak{g}$. We fix $a \in \CC^*$ and denote $Y_{i,s} = Y_{i,aq^s}$, $i \in I$, $s \in \ZZ$. Denote by $\mathcal{P}^+_\ell$ the submonoid of $\mathcal{P}^+$ generated by $Y_{i,i-2j-2}$, $i \in I$, $j \in [0, \ell]$. An object $V$ in $\mathcal{C}_{\ell}$ is a finite-dimensional $U_q(\widehat{\mathfrak{g}})$-module which satisfies the condition: for every composition factor $S$ of $V$, the highest $l$-weight of $S$ is a monomial in $\mathcal{P}^+_\ell$, \cite{HL10}. Simple modules in $\mathcal{C}_{\ell}$ are of the form $L(M)$ (cf. \cite{CP94}, \cite{HL10}), where $M \in \mathcal{P}_{\ell}^+$. Denote by $K_0(\mathcal{C}_{\ell})$ the Grothendieck ring of $\mathcal{C}_{\ell}$. By a slight abuse of notation, sometimes we write $[L(M)]$ ($M \in \mathcal{P}^+$) in $K_0(\mathcal{C}_{\ell})$ as $L(M)$ or as $[M]$. 

Hernandez and Leclerc proved that there is a cluster algebra structure on $K_0(\mathcal{C}_{\ell})$ \cite{HL10}. They also proved in Section 13 in \cite{HL10} that the cluster algebra $K_0(\mathcal{C}_{\ell})$ is isomorphic to the cluster algebra $\CC[\Gr(k,n,\sim)]$, where $n=k+\ell+1$. In \cite[Section 3]{CDFL}, it is shown that this isomorphism induces an isomorphism from the monoid $\SSYT(k,[n],\sim)$ to $\mathcal{P}_{\ell}^+$. Moreover, the dual canonical basis elements in $\CC[\Gr(k,n,\sim)]$ is in one to one correspondence with simple modules in $\mathcal{C}_{\ell}$. The dual canonical basis elements in $\CC[\Gr(k,n,\sim)]$ are in one to one correspondence with semistandard Young tableaux in $\SSYT(k,[n],\sim)$. In particular, cluster monomials in $\CC[\Gr(k,n,\sim)]$ correspond to semistandard Young tableaux in $\SSYT(k,[n],\sim)$. Explicitly, the correspondence between simple modules in $\mathcal{C}_{\ell}$ and semistandard Young tableaux in $\SSYT(k,[n],\sim)$ is given by the following map:
\begin{align*}
Y_{i,s} \mapsto T_{i,s},
\end{align*}
where $T_{i,s}$ is a fundamental tableau with entries $\{ \frac{i-s}{2}, \ldots, k+\frac{i-s}{2} \} \setminus \{ k - \frac{i+s}{2} \}$. We denote by $T_M$ the tableau corresponding to a simple module $L(M)$ and denote by $L(M_T)$ the simple module corresponding to a tableau $T$. For example, the simple module $L(Y_{1,-5}Y_{1,-3}Y_{2,-2}Y_{2,0})$ corresponds to the semistandard Young tableau corresponds to the tableau $\scalemath{0.6}{\begin{ytableau}
1 & 2 \\ 3 & 4 \\ 5 & 6
\end{ytableau}}$. 

We say that a g-vector of a simple module $L(M)$ in $\mathcal{C}_{\ell}$ is real (resp. prime) if the module $L(M)$ is real (resp. prime). We say that a tableau $T \in \SSYT(k,[n])$ is real (resp. prime) if $\ch(T)^2 = \ch(T \cup T)$ (resp. $\ch(T) \ne \ch(T')\ch(T'')$ for any non-empty tableaux $T', T''$).
We say that $T \in \SSYT(k,[n])$ is a cluster monomial (resp. cluster variable) if $\ch(T)$ is a cluster monomial (resp. variable) in $\CC[\Gr(k,n)]$. 
 
\subsection{Minimal affinizations} \label{subsec:minimal affinizations}

We recall the definition of certain modules of $U_q(\widehat{\mathfrak{g}})$ called minimal affinizations introduced by Chari in \cite{Cha95} and studied by Chari and Pressley in \cite{Cha95, CP95b, CP96a, CP96b}. Minimal affinization will be used in Section \ref{subsec:degree of g vector and decomposition of dominant monomial} to compute the degree of a g-vector (also the degree of the polynomial in $\CC[\Gr(k,n)]$ corresponding to $g$).

For a Laurent monomial $m = \prod_{i \in I, a \in \CC^*} Y_{i,a}^{u_{i,a}}$ in $\mathcal{P}$, denote $\omega(m) = \sum_{i \in I, a\in \CC^*} u_{i,a} \omega_i$. For a simple finite dimensional $U_q(\mathfrak{g})$-module $V$, a simple finite dimensional $U_q(\widehat{\mathfrak{g}})$-module $L(m)$ is said to be an affinization of $V$ if $\omega(m)$ is the highest weight of $V$. For a $U_q(\mathfrak{g})$-module $V$ and a domiant weight $\lambda \in P^+$, denote by $m_{\lambda}(V)$ the multiplicity of the simple $U_q(\mathfrak{g})$-module of highest weight $\lambda$ in $V$. Two affinizations are said to be equivalent if they are isomorphic as $U_q(\mathfrak{g})$-modules. Let $\mathcal{Q}_V$ be the set of equivalence classes of affinizations of $V$ and denote by $[L]$ the equivalence class of an affinization $L$ of $V$. There is a partial order on $\mathcal{Q}_V$ defined as follows: for $[L], [L'] \in \mathcal{Q}_V$, $[L] \le [L']$ if and only if for all $\mu \in P^+$, either $m_{\mu}(L) \le m_{\mu}(L')$ or there exists $\nu > \mu$ such that $m_{\nu}(L) < m_{\mu}(L')$. A minimal affinization of $V$ is a minimal element of $\mathcal{Q}_V$ with respect to this partial order. 

Minimal affinizations are classified in terms of highest $l$-weight monomials in \cite{Cha95, CP95b, CP96a, CP96b}. Minimal affinizations of $U_q(\widehat{\mathfrak{sl}_k})$ are evaluation modules of $U_q(\widehat{\mathfrak{sl}_k})$ \cite{Cha95}. Minimal affinizations of $U_q(\widehat{\mathfrak{sl}_k})$ are of the form $L(m)$, where
\begin{align} \label{eq:one type of minimal affinizations which corresponds to Plucker coordinates}
m = \prod_{i=1}^{k-1} \prod_{a_i=0}^{m_i-1} Y_{i, s+2a_i+i-1+2\sum_{j=1}^{i-1} m_j},
\end{align}
or 
\begin{align}
m = \prod_{i=1}^{k-1} \prod_{a_i=0}^{m_i-1} Y_{i, -s-2a_i-i+1-2\sum_{j=1}^{i-1} m_j},
\end{align}
for some $s \in \ZZ$, $m_j \in \ZZ_{\ge 0}$. In particular, according to the correspondence in Section \ref{subsec:correspondence between modules and tableaux}, minimal affinizations (\ref{eq:one type of minimal affinizations which corresponds to Plucker coordinates}) correspond to one-column tableaux in $\SSYT(k,[n])$ and Pl\"{u}cker coordinates in $\CC[\Gr(k,n)]$ \cite{CDFL}.

\section{Quasi-automorphism actions on g-vectors of Grassmannian cluster algebras} \label{sec:braid group actions on g vectors}

In this section, we study tropicalisation of quasi-automorphisms of Grassmannian cluster algebras. These include the cyclic rotation $\rho$, reflection $\theta$, twist map $\tau$ and the generators of Fraser's braid group action $\sigma_i$. As described in Section \ref{Quasihom}, these send cluster variables to cluster variables (up to factors monomial in frozen variables in the case of $\tau$ and $\sigma_i$).

For a map $f: \Gr(k,n) \to \Gr(k,n)$, we define the pullback $f^*: \CC[\Gr(k,n)] \to \CC[\Gr(k,n)]$ by $f^*(\varphi)(M) = \varphi(f(M))$, $M \in \Gr(k,n)$, $\varphi \in \CC[\Gr(k,n)]$. 

\begin{remark} 
We remind the reader that, as in Sec. \ref{subsec:Grassmannian cluster algebras},
our notation $\Gr(k,n)$ in this paper has the same meaning as $\widetilde{\Gr}^{\circ}(k,n)$ in \cite{Fra20}. It is the Zariski-open subset of the affine cone over the Grassmannian of $k$-dimensional subspaces in $\CC^n$ cut out by the non-vanishing of the frozen variables. 
\end{remark}

\subsection{Cyclic rotation}
The cyclic rotation map $\rho: \Gr(k,n) \to \Gr(k,n)$ is given by 
\begin{align*}
\bigl(z_1, \ldots, z_n\bigr) \mapsto \bigl(z_2, z_3, \ldots, z_{n}, (-1)^{k-1}z_1\bigr),
\end{align*}
where $z_i$'s are columns of the matrix $(z_1, \ldots, z_n)$. The inverse map $\rho^{-1}$ is given by
\begin{align*}
\bigl(z_1, \ldots, z_n\bigr) \mapsto \bigl((-1)^{k-1}z_n, z_1, z_2, \ldots, z_{n-1}\bigr).
\end{align*}
The cyclic rotation map $\rho: \Gr(k,n) \to \Gr(k,n)$ induces an algebra automorphism $\rho^*$ on $\CC[\Gr(k,n)]$. The automorphism $\rho^*$ is a cluster automorphism on $\CC[\Gr(k,n)]$ (we remind the reader that cluster automorphisms send cluster variables to cluster variables, send clusters to clusters, send frozen variables to frozen variables, without having to remove frozen factors), see Section 5 in \cite{Fra20}. The automorphism $\rho^*$ has signature $\epsilon = +1$ in the sense of Section \ref{Quasihom}.

\subsection{Reflection}
The reflection map $\theta: \Gr(k,n) \to \Gr(k,n)$ is given by 
\begin{align*}
(z_1, \ldots, z_n) \mapsto (-1)^{\binom{k}{2}}(z_n, z_{n-1}, \ldots, z_{1}).
\end{align*}
We have that $\theta^{-1} = \theta$. The reflection map $\theta: \Gr(k,n) \to \Gr(k,n)$ induces an algebra automorphism $\theta^*$ on $\CC[\Gr(k,n)]$. The automorphism $\theta^*$ is a cluster automorphism (up to reversing all the arrows in the quiver in a seed) on $\CC[\Gr(k,n)]$, see Section 5 in \cite{Fra20}. The automorphism $\theta^*$ sends cluster variables to cluster variables, sends frozen variables to frozen variables, and sends a seed $({\bf x}, Q)$ to the seed $(\theta^*({\bf x}), -Q)$, i.e. it has signature $\epsilon = -1$ in the sense of Section \ref{Quasihom}.

\subsection{Twist map}
We now recall the definition of the twist map $\tau$, see Section 2 of \cite{MS16}. For $z_1, \ldots, z_{k-1} \in \mathbb{C}^k$, the generalized cross product $z_1 \times \cdots \times z_{k-1}$ is the unique vector in $\mathbb{C}^k$ such that 
\begin{align*}
(z_1 \times \cdots \times z_{k-1})^T z = \det(z_1, \ldots, z_{k-1}, z).
\end{align*}

For a $k \times n$ matrix $p$ with columns $p_1, \ldots, p_n$, Marsh and Scott \cite{MS16} defined a twist $p'$ of $p$ to be the $k \times n$ matrix whose $i$th column vector is 
\begin{align*}
(p')_i = \beta_i \cdot (p_{\sigma^{k-1}(i)} \times p_{\sigma^{k-2}(i)} \times \cdots \times p_{\sigma(i)}),
\end{align*}
where $\sigma: \{1, \ldots, n\} \to \{1, \ldots, n\}$ is the map given by $i \mapsto i-1 \pmod n$ and
\[
\beta_i = \begin{cases} (-1)^{i(k-i)}, & i \le k-1, \\ 1, & i \ge k. \end{cases}
\]

The twist map on the set of $k \times n$ matrices induces a twist map $\tau^*$ on $\mathbb{C}[{\rm Gr}(k, n)]$. It is shown in \cite[Proposition 8.10]{MS16} that $\tau^*$ sends a cluster variable to a cluster variable (possibly multiplied by some frozen variables). It also follows from results in \cite{MS16} that we have
\begin{equation}
   (\tau^*)^2 (P_{i_1,\ldots,i_k})
   = (\rho^*)^{-k}(P_{i_1,\ldots,i_k}) \prod_{j=1}^k \Bigl(\prod_{r=1}^{k-1} P_{i_j-r,i_j-r+1,i_j-r+2,\ldots,i_j-r+k-1}\Bigr)\,,
   \label{tausquare}
\end{equation}
where indices on the Pl\"ucker coordinates should be treated modulo $n$ and we remind the reader that the set of indices on a Pl\"ucker coordinate should be written in increasing order (e.g. $P_{6789}$ means $P_{1678}$ when $n=8$).
The twist map $\tau^*$ is a quasi-automorphism of the cluster algebra $\CC[\Gr(k,n)]$ with signature $\epsilon = +1$, see \cite{Fra20, MS16}. From (\ref{tausquare}) it obeys $(\tau^*)^2 \propto (\rho^*)^{-k}$.

\subsection{\texorpdfstring{Braid group action on $\CC[\Gr(k,n)]$}{Fraser's braid group action on C[Gr(k,n)]}}
 \label{subsec:Fraser braid group action}

We now recall Fraser's braid group action on $\CC[\Gr(k,n)]$. 
Let $d={\rm gcd}(k,n)$. For $i \in [1,d-1]$, the map 
%$\sigma_i: (V^n)^{\circ} \to (V^n)^{\circ}$ 
$\sigma_i: \Gr(k,n) \to \Gr(k,n)$ 
is defined as follows (see Definition 5.2 and Equation (18) in Remark 5.6 of \cite{Fra20}):
\begin{equation}
\sigma_i :\, z_l
\mapsto
\begin{cases}
z_{l+1}, & l=i \: (\text{mod } d),\\
{\zeta}^{(j)}_i, & l=i+1 \: (\text{mod } d),\\
z_l, & \text{otherwise},
\end{cases}
\end{equation}
where $l \in [n]$, $j = \frac{l-i-1}{d}$, and
\begin{equation}
\label{braidzeta}
{\zeta}^{(j)}_i
=
\frac{\det(z_{i+jd}, z_{i+jd+2}, \ldots, z_{i+jd+k}) z_{i+jd+1}}{\det( z_{i+jd+1}, z_{i+jd+2}, \ldots, z_{i+jd+k} )} -  z_{i+jd}.
\end{equation}
The map $\sigma_i^{-1}$ is defined similarly:
\begin{equation} 
\sigma_i^{-1} :\, z_l
\mapsto
\begin{cases}
{\tilde \zeta}^{(j)}_i, &  l=i \: (\text{mod } d),\\
z_{l-1}, & l=i+1 \: (\text{mod } d),\\
z_l, & \text{otherwise}\, ,
\end{cases}
\end{equation}
where $l \in [n]$, $j = \frac{l-i}{d}$, and
\begin{equation}
\label{braidzetatilde}
{\tilde \zeta}^{(j)}_i
=
 \frac{\det(z_{i+jd-k+1}, \ldots, z_{i+jd-1},z_{i+jd+1})}{\det(z_{i+jd-k+1}, \ldots, z_{i+jd-1},z_{i+jd})} z_{i+jd} -  z_{i+jd+1}.
\end{equation}

In Theorem 5.3 of \cite{Fra20}, Fraser proved that the maps $\sigma_i$, $\sigma_i^{-1}$ are inverse regular automorphisms on %$(V^n)^{\circ}$
$\Gr(k,n)$
. They obey the braid relations
\begin{align}
\label{braidrels}
    &\sigma_i \sigma_{i+1} \sigma_i = \sigma_{i+1} \sigma_i \sigma_{i+1}\,\text{ for } i\in \{1\ldots,d-2\}, \\
    &\sigma_i \sigma_j = \sigma_j \sigma_i\,, \text{ for } |i-j| \geq 2 \text{ where } i,j \in \{1,\dots,d-1\}\,. \notag
\end{align}
Note also that they obey
\begin{align}
    \rho^{-1} \sigma_i \rho &= \sigma_{i+1}\,, \qquad i=1,\ldots,d-2\,,\\
    \rho^{-1} \sigma_{d-1} \rho &= \rho \sigma_1 \rho^{-1}\,. \notag
\end{align}
We may define an element $\sigma_d = \rho^{-1} \sigma_{d-1} \rho = \rho \sigma_1 \rho^{-1}$ such that conjugation by $\rho$ acts as a cyclic permutation of the $\sigma_i$.

When $d=k$ we have the additional property,
\begin{equation}
\label{cycfrombraids}
    \sigma_{d-1}\ldots\sigma_1: z_{l} \mapsto 
\begin{cases}
    \rho (z_{l})\frac{\det(z_{l-k+1},z_{l-k+2},\ldots,z_{l})}{\det(z_{l-k+2},z_{l-k+3},\ldots,z_{l+1})} \,, & l=0 \text{ (mod $k$)},\\
    \rho (z_{l})\,, & \text{otherwise.}
\end{cases}
\end{equation}

The pullbacks $\sigma_i^*$, $(\sigma_i^*)^{-1}$ are quasi-inverse quasi-automorphisms on $\CC[\Gr(k,n)]$ with signature $\epsilon = +1$. The relations (\ref{braidrels}) imply the braid relations for the $\sigma_i^*$, while (\ref{cycfrombraids}) implies $(\sigma_{d-1}\ldots \sigma_1)^* \propto \rho^*$ when $d=k$.

Note that we can consider a version $\tilde{\sigma}_i$ of the braid generators where the denominators are cleared in (\ref{braidzeta}) and (\ref{braidzetatilde}). The corresponding pullbacks $\tilde{\sigma}_i^*$ and $(\tilde{\sigma}_i^{-1})^*$ are proportional to $\sigma_i^*$ and $(\sigma_i^*)^{-1}$, quasi-inverse to each other and obey the braid relations up to proportionality.

\subsection{Quasi-automorphism actions on g-vectors} \label{subsec:braid group action on g vectors}

In this subsection, we study the tropical version of quasi-automorphism actions on Grassmannian cluster algebras, i.e. quasi-automorphism actions on g-vectors. In the following we will exclusively use the pullbacks of the maps introduced in the preceding sections and, by a slight abuse of notation, simply refer to them as $\rho,\theta,\tau,\sigma_i$ etc.

We choose the web matrix $W$ which is defined as $W = ( \mathds{1}_{k} \mid M )$, where $M=(M_{ij})_{k \times (n-k)}$, and 
\begin{align*}
M_{ij} = (-1)^{k+i} \sum_{0 \le \lambda_{k-i} \le \cdots \le \lambda_1 \le j-1} \prod_{b=1}^{k-i} \prod_{a=1}^{\lambda_b} \chi_{a,b},
\end{align*}
see \cite{Postnikov} and \cite{SW05}. This choice of web matrix is different from the choice of web matrix in \cite{DFGK20} (the $\chi_{a,b}$ here corresponds to $\chi_{a,b}$ in \cite{DFGK20}) because we use a different labeling of the vertices of the initial quiver in Section \ref{subsec:Grassmannian cluster algebras}.

Take the initial seed $\Sigma_{t_0}$ described in Section \ref{subsec:Grassmannian cluster algebras}. Evaluating the $\hat{y}$-variables of this labelled seed on the web matrix $W$, we have that
\begin{equation}
    (\hat{y}_{1;t_0},\ldots,\hat{y}_{(n-k-1)(k-1);t_0}) = ( \chi_{1,1}, \chi_{1,2}, \ldots, \chi_{1,k-1}, \chi_{2,1},\ldots, \chi_{n-k-1,k-1} )\,.
\end{equation}
Let $f$ be a quasi-automorphism on $\mathbb{C}[{\rm Gr}(k,n)]$ and $f^{-1}$ a choice of quasi-inverse. By evaluating the action of $f$ and $f^{-1}$ on the fundamental Pl\"ucker variables, we may extend their actions in the natural way to the field $\mathcal{F} = \mathbb{C}\bigl({\rm Gr}(k,n)\bigr)$ of rational functions on ${\rm Gr}(k,n)$. In particular we may evaluate the action of $f^{-1}$ on the $\hat{y}$ variables of the initial seed,
\begin{equation}
\label{xabvarsfinv}
    f^{-1} : ( \chi_{1,1}, \chi_{1,2}, \ldots, \chi_{n-k-1,k-1} ) \mapsto ( \xi_{1,1}, \xi_{1,2}, \ldots, \xi_{n-k-1,k-1} ).
\end{equation}
If $f$ (and hence $f^{-1}$) has signature $\epsilon = +1$ then (up to some permutation $\pi^{-1}$ of the node labels) the variables $\xi_{a,b}$ are the $\hat{y}$ variables of another labelled seed $\Sigma_{t}$ which is the image of the initial seed $\Sigma_{t_0}$ under $f^{-1}$,
\begin{equation}
( \xi_{1,1}, \xi_{1,2}, \ldots, \xi_{n-k-1,k-1} ) = (\hat{y}_{\pi^{-1}(1);t},\ldots,\hat{y}_{\pi^{-1}((n-k-1)(k-1));t})\,.
\end{equation}
If $f$ has signature $\epsilon=-1$ then the $\xi_{a,b}$ are the inverses of the $\hat{y}$ variables of $\Sigma_t$,
\begin{equation}
( \xi_{1,1}, \xi_{1,2}, \ldots, \xi_{n-k-1,k-1} ) = (\hat{y}_{\pi^{-1}(1);t}^{-1},\ldots,\hat{y}_{\pi^{-1}((n-k-1)(k-1));t}^{-1})\,.
\end{equation}
The variables $\xi_{a,b}$ are therefore birationally related to the original variables $\chi_{a,b}$. 

Now we wish to obtain the action of $f$ on g-vectors ${\bf g}(x) \equiv {\bf g}(x;t_0)$, i.e. we want to obtain ${\bf g}\bigl(h_f(x)\bigr)$, ${\bf g}^{\rm op}\bigl(h_f(x)\bigr)$ from ${\bf g}(x)$, ${\bf g}^{\rm op}(x)$.
Following Section \ref{sec:tropicalisation of quasi homomorphisms} we introduce the tropicalised variables,
\begin{align}
    v_i &= {\rm Trop}_+ \hat{y}_{i;t}, &&v'_i = {\rm Trop}_+ \hat{y}_{i;t_0}\,, \\
\notag    w_i &= {\rm Trop}_- \hat{y}_{i;t}, &&w'_i = {\rm Trop}_- \hat{y}_{i;t_0}\,.
\end{align}
We may also (somewhat redundantly) introduce the tropicalised versions of $\chi_{a,b}$ and $\xi_{a,b}$,
\begin{align}
    &v_{a,b} = {\rm Trop}_+ (\xi_{a,b}), &&v'_{a,b} = {\rm Trop}_+(\chi_{a,b}),\\
    \notag
    &w_{a,b} = {\rm Trop}_- (\xi_{a,b}), &&w'_{a,b} = {\rm Trop}_-(\chi_{a,b}).
\end{align}
Clearly we have
\begin{align}
    (v_{1,1},v_{1,2},\ldots,v_{n-k-1,k-1}) &= \epsilon (v_1,\ldots v_{(n-k-1)(k-1)})^{\pi^{-1}}, \\
    \notag
    (v'_{1,1},v'_{1,2},\ldots,v'_{n-k-1,k-1}) &= (v'_1,\ldots v'_{(n-k-1)(k-1)})
\end{align}
and the same relations obtained by swapping $w$ for $v$ everywhere.

The results of Section \ref{sec:tropicalisation of quasi homomorphisms}, specifically Corollary \ref{cor-quasiauto}, tell us that (setting $t'=t_0$ there and writing ${\bf g}(x;t_0) = {\bf g}(x)$ etc.)
\begin{align}
\label{gvecmapsGrass}
        &{\bf g}\bigl(h_f(x)\bigr) = \bigl(Q^+_{t,t_0} {\bf g}(x)\bigr)^{\pi^{-1}}, &&{\bf g}^{\rm op}\bigl(h_f(x)\bigr) = \bigl(Q^-_{t,t_0} {\bf g}^{\rm op}(x)\bigr)^{\pi^{-1}}, &&(\epsilon = +1), \\
        \notag
        &{\bf g}\bigl(h_f(x)\bigr) = \bigl(Q^-_{t,t_0} {\bf g}^{\rm op}(x)\bigr)^{\pi^{-1}},  &&{\bf g}^{\rm op}\bigl(h_f(x)\bigr) = \bigl(Q^+_{t,t_0} {\bf g}(x;)\bigr)^{\pi^{-1}}, &&(\epsilon=-1).
    \end{align}
We remind the reader that the piecewise linear maps $Q^\pm_{t,t_0}$ are obtained by tropicalising the relations between the $\hat{y}$ variables in the labelled seeds $\Sigma_t$ and $\Sigma_{t_0}$,
\begin{align}
    Q^+_{t,t_0}: (v'_{1},\ldots,v'_{(n-k-1)(k-1)}) &\mapsto (v_{1},\ldots,v_{(n-k-1)(k-1)}), \\
    \notag
Q^-_{t,t_0}: (w'_{1},\ldots,w'_{(n-k-1)(k-1)}) &\mapsto (w_{1},\ldots,w_{(n-k-1)(k-1)}).
\end{align}
The relations (\ref{gvecmapsGrass}) therefore become
\begin{align}
\label{gvecmapsGrass2}
        &{\bf g}\bigl(h_f(x)\bigr) = Q^+_f {\bf g}(x), &&{\bf g}^{\rm op}\bigl(h_f(x)\bigr) = Q^-_f {\bf g}^{\rm op}(x;t_0), &&(\epsilon = +1), \\
        \notag
        &{\bf g}\bigl(h_f(x)\bigr) = -Q^-_f {\bf g}^{\rm op}(x),  &&{\bf g}^{\rm op}\bigl(h_f(x)\bigr) = -Q^+_{f} {\bf g}(x), &&(\epsilon=-1),
    \end{align}
where $Q_f^\pm$ are the piecewise linear maps obtained by simply tropicalising the map (\ref{xabvarsfinv}),
\begin{align}
\label{Qfmaps}
Q^+_{f}: (v'_{1,1},\ldots,v'_{n-k-1,k-1}) &\mapsto (v_{1,1},\ldots,v_{n-k-1,k-1}), \\
\notag
Q^-_{f}: (w'_{1,1},\ldots,w'_{n-k-1,k-1}) &\mapsto (w_{1,1},\ldots,w_{n-k-1,k-1}).
\end{align}
We note that the permutation $\pi$ does not appear in (\ref{gvecmapsGrass2}) or (\ref{Qfmaps}).

 We conjecture the following extension to any element of the dual canonical basis. 
\begin{conjecture}
The relations (\ref{gvecmapsGrass2}) hold for any quasi-automorphism $f$ and any element $x$ in the dual canonical basis of $\CC[\Gr(k,n)]$, not only for cluster variables. 
\end{conjecture}

\begin{remark} \label{remark:theta_compose_with_tau_is_negative_of_reversing_g_vector}
    We note that the case $f=\theta \circ \tau = f^{-1}$ is particularly simple. We have
    \begin{align}
    f ( \chi_{1,1}, \chi_{1,2}, \ldots, \chi_{n-k-1,k-1} ) \mapsto (\chi_{n-k-1,k-1}^{-1}, \ldots , \chi_{1,1}^{-1})\,,
    \end{align}
    i.e. $f$ inverts and reverses the order of the $\hat{y}$ variables of the initial cluster. We therefore have a case of a quasi-automorphism of signature $\epsilon = -1$ where the seed $\Sigma_t$ above is $\Sigma_{t_0}$ and the permutation $\pi = \pi^{-1}$ is the map which reverses the order of the node labels. We also have $Q_f^+ = Q_f^- = - R_\pi$ where $R_\pi$ is the linear map (rather than piecewise linear map) which reverses the order of the components of vectors. The equations in the second line of (\ref{gvecmapsGrass2}) then become
    \begin{equation}
    \label{thetataugmaps}
        {\bf g}\bigl(\theta (h_{\tau}(x))\bigr) = {\bf g}^{\rm op}(x)^\pi, \qquad \qquad {\bf g}^{\rm op}\bigl(\theta(h_{\tau}(x))\bigr) = {\bf g}(x)^\pi,
    \end{equation}
    which are equivalent.
\end{remark}

\begin{example} \label{example:Gr36}
Consider the Grassmannian cluster algebra $\CC[\Gr(3,6)]$. There are six frozen cluster variables and a total of sixteen active cluster variables arranged across 50 seeds. The active variables $x$ and the corresponding g-vectors ${\bf g}(x) \equiv {\bf g}(x;t_0)$ for the algebra and its opposite are given in Table \ref{Gr36gvecs}.
{\small
  \renewcommand{\arraystretch}{1.2}
\begin{table}
  \centering
  \begin{tabular}{ccc}
    \toprule
    $x$&${\bf g}(x)$&${\bf g}^{\rm op}(x)$\\
    \midrule
    $P_{124}$ & (1,0,0,0) & (1,0,0,0) \\
    $P_{134}$ & (0,1,0,0) & (0,1,0,0) \\
    $P_{125}$ & (0,0,1,0) & (0,0,1,0) \\
    $P_{145}$ & (0,0,0,1) & (0,0,0,1) \\
    $P_{245}$ & (0,-1,0,1) & (1,-1,0,0) \\
    $P_{235}$ & (-1,0,1,0) & (0,-1,0,0) \\
    $P_{256}$ & (0,-1,0,0) & (0,0,1,-1) \\
    $P_{236}$ & (-1,0,0,0) & (0,0,0,-1) \\
    \bottomrule
  \end{tabular}
  \begin{tabular}{ccc}
    \toprule
    $x$&${\bf g}(x)$&${\bf g}^{\rm op}(x)$\\
    \midrule
    $P_{356}$ & (0,0,0,-1) & (-1,0,0,0) \\
    $P_{346}$ & (0,0,-1,0) & (0,1,0,-1) \\
    $P_{136}$ & (-1,1,0,0) & (0,0,-1,0) \\
    $P_{146}$ & (0,0,-1,1) & (1,0,-1,0) \\
    $P_{135}$ & (-1,1,1,0) & (-1,0,0,1) \\
    $P_{246}$ & (0,-1,-1,1) & (1,0,0,-1) \\
    $q_1$ & (1,0,0,-1) & (0,1,1,-1) \\
    $q_2$ & (-1,0,0,1) & (1,-1,-1,0) \\
    \bottomrule
  \end{tabular}
  \caption{The active cluster variables for the ${\rm Gr}(3,6)$ cluster algebra comprise the fourteen non-frozen Pl\"ucker variables and two quadratics $q_1 = P_{124}P_{356} - P_{123}P_{456}$ and $q_2 = P_{235}P_{146} - P_{234}P_{156}$. We note here the associated ${\bf g}(x)$ and ${\bf g}^{\rm op}(x)$.}
  \label{Gr36gvecs}
\end{table}
}

Writing the $\hat{y}$ variables of the initial seed we have
\begin{equation}
    (\hat{y}_1,\hat{y}_2,\hat{y}_3,\hat{y}_4) = (\chi_{1,1}, \chi_{1,2}, \chi_{2,1},\chi_{2,2}) = \Bigl(\tfrac{P_{123} P_{145}}{P_{125} P_{134}}, \tfrac{P_{124} P_{345}}{P_{145}P_{234}}, \tfrac{P_{124}P_{156}}{P_{126}P_{145}} , \tfrac{P_{125}P_{134}P_{456}}{P_{124}P_{156}P_{345}} \Bigr)\,.
\end{equation}
Applying the inverse cyclic rotation map $\rho^{-1}$ to the above expressions and evaluating the result on the web matrix $W$ we find
\begin{equation}
{\renewcommand{\arraystretch}{1.2}
    \rho^{-1} :
    \begin{bmatrix} 
    \chi_{1,1} \\\chi_{1,2} \\ \chi_{2,1} \\ \chi_{2,2} 
    \end{bmatrix} 
    \mapsto 
            \begingroup
    \renewcommand{\arraystretch}{1.8}
    \begin{bmatrix}
    {
    \frac{\chi_{1,2} (1 + \chi_{2,1} + \chi_{2,1} \chi_{2,2})}{(1 + \chi_{2,1}) (1 + \chi_{1,1} + \chi_{1,1} \chi_{1,2} + \chi_{1,1} \chi_{2,1} + \chi_{1,1} \chi_{1,2} \chi_{2,1} + 
   \chi_{1,1} \chi_{1,2} \chi_{2,1} \chi_{2,2})}
   }\\
   { 
   \frac{(1 + \chi_{1,1} + \chi_{1,1} \chi_{2,1}) \chi_{2,2}}{1 + \chi_{2,1} + \chi_{2,1} \chi_{2,2}}
   }
   \\
   \frac{1 + \chi_{1,1} + \chi_{1,1} \chi_{2,1}}{\chi_{1,1} \chi_{1,2} (1 + \chi_{2,1} + \chi_{2,1} \chi_{2,2})}
   \\
   \frac{(1 + \chi_{2,1}) (1 + \chi_{1,1} + \chi_{1,1} \chi_{1,2} + \chi_{1,1} \chi_{2,1} + \chi_{1,1} \chi_{1,2} \chi_{2,1} + 
   \chi_{1,1} \chi_{1,2} \chi_{2,1} \chi_{2,2})}{\chi_{2,1} (1 + \chi_{1,1} + \chi_{1,1} \chi_{2,1}) \chi_{2,2}}
   \end{bmatrix}
   \endgroup.
   }
\end{equation}
We verify (\ref{gvecmapsGrass2}) holds (with $\epsilon=+1$), i.e. that tropicalising the above rational map with `max' convention gives the piecewise linear map $Q_{\rho}^+$ which correctly maps each ${\bf g}(x)$ to ${\bf g}\bigl( \rho(x)\bigr)$ in Table \ref{Gr36gvecs}. Similarly, tropicalising with `min' convention gives $Q_{\rho}^-$ which maps each ${\bf g}^{\rm op}(x)$ to ${\bf g}^{\rm op}\bigl(\rho(x)\bigr)$. Note that in this case we have $h_{\rho} = \rho$. For example, we have ${\bf g}(P_{124}) = (1,0,0,0)$ and ${\bf g}^{\rm op}(P_{124}) = (1,0,0,0)$. Applying $Q^\pm_{\rho}$ we find
\begin{align}
    Q^+_{\rho} {\bf g}(P_{124}) &= (-1,0,1,0) = {\bf g}(P_{235}), \\
    \notag
    Q^-_{\rho} {\bf g}^{\rm op}(P_{124}) &= (0,-1,0,0) = {\bf g}^{\rm op}(P_{235}),
\end{align}
and we note that indeed $\rho (P_{124}) = P_{235}$.

Similarly, applying the reflection map $\theta$ (note $\theta^{-1} = \theta$) to the $\hat{y}$ variables of the initial seed, we have
\begin{equation}
{\renewcommand{\arraystretch}{1.2}
    \theta :
    \begin{bmatrix} 
    \chi_{1,1} \\ \chi_{1,2} \\ \chi_{2,1} \\ \chi_{2,2} 
    \end{bmatrix} 
    \mapsto 
            \begingroup
    \renewcommand{\arraystretch}{1.8}
    \begin{bmatrix}
    {
    \frac{\chi_{2,2} (1 + \chi_{1,1} + \chi_{1,1} \chi_{1,2} + \chi_{1,1} \chi_{2,1} + \chi_{1,1} \chi_{1,2} \chi_{2,1} + \chi_{1,1} \chi_{1,2} \chi_{2,1} \chi_{2,2})}{(1 + \chi_{1,2} + \chi_{1,2} \chi_{2,2}) (1 + \chi_{2,1} + \chi_{2,1} \chi_{2,2})}
   }\\
   { 
   \frac{\chi_{2,1} (1 + \chi_{2,2} + \chi_{1,1} \chi_{2,2})}{1 + \chi_{1,1} + \chi_{1,1} \chi_{1,2} + \chi_{1,1} \chi_{2,1} + \chi_{1,1} \chi_{1,2} \chi_{2,1} + \chi_{1,1} \chi_{1,2} \chi_{2,1} \chi_{2,2}}
   }\\
   \frac{\chi_{1,2} (1 + \chi_{2,2} + \chi_{1,1} \chi_{2,2})}{1 + \chi_{1,1} + \chi_{1,1} \chi_{1,2} + \chi_{1,1} \chi_{2,1} + \chi_{1,1} \chi_{1,2} \chi_{2,1} + \chi_{1,1} \chi_{1,2} \chi_{2,1} \chi_{2,2}}
   \\
   \frac{\chi_{1,1} (1 + \chi_{1,2} + \chi_{1,2} \chi_{2,2}) (1 + \chi_{2,1} + \chi_{2,1} \chi_{2,2})}{1 + \chi_{2,2} + \chi_{1,1} \chi_{2,2}}
   \end{bmatrix}
    \endgroup
    .
   }
\end{equation}
We verify (\ref{gvecmapsGrass2}) holds (with $\epsilon=-1$), i.e. that tropicalising the above rational map with `max' convention gives the piecewise linear map $Q_{\theta}^+$ which correctly maps each ${\bf g}(x)$ to the negative of ${\bf g}^{\rm op}\bigl( \theta(x)\bigr)$ in Table \ref{Gr36gvecs}. Similarly, tropicalising with `min' convention gives $Q_{\theta}^-$ which maps each ${\bf g}^{\rm op}(x)$ to the negative of ${\bf g}\bigl(\theta(x)\bigr)$. Note that in this case we have $h_{\theta} = \theta$. For example, we have
\begin{align}
    -Q^+_{\theta} {\bf g}(P_{124}) &= (-1,0,0,0) = {\bf g}^{\rm op}(P_{356}), \\
    \notag
    -Q^-_{\theta} {\bf g}^{\rm op}(P_{124}) &= (0,0,0,-1) = {\bf g}(P_{356}),
\end{align}
and we note that $\theta(P_{124}) = P_{356}$.

Applying the inverse twist map $\tau^{-1}$ to the $\hat{y}$ variables of the initial seed, we have
\begin{equation}
{\renewcommand{\arraystretch}{1.2}
    \tau^{-1} :
    \begin{bmatrix} 
    \chi_{1,1} \\\chi_{1,2} \\ \chi_{2,1} \\ \chi_{2,2} 
    \end{bmatrix} 
    \mapsto 
        \begingroup
    \renewcommand{\arraystretch}{1.8}
    \begin{bmatrix}
    {
    \frac{1 + \chi_{1,1} + \chi_{1,1} \chi_{1,2} + \chi_{1,1} \chi_{2,1} + \chi_{1,1} \chi_{1,2} \chi_{2,1} + \chi_{1,1} \chi_{1,2} \chi_{2,1} \chi_{2,2}}{(1 + \chi_{1,1} + \chi_{1,1} \chi_{1,2}) (1 + \chi_{1,1} + \chi_{1,1} \chi_{2,1}) \chi_{2,2}}
   }
   \\
   { 
   \frac{\chi_{2,1} (1 + \chi_{2,2} + \chi_{1,1} \chi_{2,2})}{1 + \chi_{1,1} + \chi_{1,1} \chi_{1,2} + \chi_{1,1} \chi_{2,1} + \chi_{1,1} \chi_{1,2} \chi_{2,1} + \chi_{1,1} \chi_{1,2} \chi_{2,1} \chi_{2,2}}
   }
   \\
   \frac{\chi_{1,2} (1 + \chi_{2,2} + \chi_{1,1} \chi_{2,2})}{1 + \chi_{1,1} + \chi_{1,1} \chi_{1,2} + \chi_{1,1} \chi_{2,1} + \chi_{1,1} \chi_{1,2} \chi_{2,1} + \chi_{1,1} \chi_{1,2} \chi_{2,1} \chi_{2,2}}
   \\
   \frac{(1 + \chi_{1,1} + \chi_{1,1} \chi_{1,2}) (1 + \chi_{1,1} + \chi_{1,1} \chi_{2,1})}{\chi_{1,1} \chi_{1,2} \chi_{2,1} (1 + \chi_{2,2} + \chi_{1,1} \chi_{2,2})}
   \end{bmatrix}
   \endgroup .
   }
\end{equation}
We verify (\ref{gvecmapsGrass2}) holds (with $\epsilon=+1$), i.e. that tropicalising the above rational map with `max' convention gives the piecewise linear map $Q_{\tau}^+$ which correctly maps each ${\bf g}(x)$ to ${\bf g}\bigl( h_{\tau}(x)\bigr)$ in Table \ref{Gr36gvecs}. Similarly, tropicalising with `min' convention gives $Q_{\tau}^-$ which maps each ${\bf g}^{\rm op}(x)$ to ${\bf g}^{\rm op}\bigl(h_{\tau}(x)\bigr)$. Note that in this case we have $h_{\tau} \neq \tau$, as $\tau$ generally introduces frozen factors. For example, we have
\begin{align}
    &Q_{\tau}^+ {\bf g}(P_{124}) = (-1,0,0,0) = {\bf g}(P_{236}),\\
    \notag
    &Q_{\tau}^- {\bf g}^{\rm op}(P_{124}) = (0,0,0,-1) = {\bf g}^{\rm op}(P_{236}),
\end{align}
and we note that $\tau(P_{124}) = P_{156}P_{236}$ and therefore ($P_{156}$ being a frozen variable) we have $h_{\tau}(P_{124}) = P_{236}$. We also verify that (\ref{thetataugmaps}) holds.

Since $d=\gcd(k,n)=\gcd(3,6)=3$, we have two generators of the braid group ${\rm Br}_d$ to consider in this example. Applying $\sigma_1^{-1}$ to the $\hat{y}$ variables of the initial seed we have
\begin{equation}
{\renewcommand{\arraystretch}{1.2}
    \sigma_1^{-1} :
    \begin{bmatrix} 
    \chi_{1,1} \\\chi_{1,2} \\ \chi_{2,1} \\ \chi_{2,2} 
    \end{bmatrix} 
    \mapsto 
    \begingroup
    \renewcommand{\arraystretch}{1.5}
    \begin{bmatrix}
    {
    \frac{\chi_{2,2}}{1 + \chi_{1,2} + \chi_{1,2} \chi_{2,2}}
   }\\
   { 
   \frac{1 + \chi_{1,2}}{\chi_{1,2} \chi_{2,2}}
   }\\
   {\scriptstyle \chi_{1,1} (1 + \chi_{1,2})}
   \\
   \frac{\chi_{2,1} (1 + \chi_{1,2} + \chi_{1,2} \chi_{2,2})}{1 + \chi_{1,2}}
   \end{bmatrix}
   \endgroup 
   .
   }
\end{equation}
Tropicalising the above map gives the maps $Q^\pm_{\sigma_1}$, which in this case are simple enough to write explicitly. We have
\begin{equation}
{\renewcommand{\arraystretch}{1.2}
    Q^+_{\sigma_1} :
    \begin{bmatrix} 
    v_{1,1} \\v_{1,2} \\ v_{2,1} \\ v_{2,2} 
    \end{bmatrix} 
    \mapsto 
    \begin{bmatrix}
    {v_{2,2} - {\rm max}(0,v_{1,2},v_{1,2}+v_{2,2})
   }\\
   { {\rm max}(0,v_{1,2})-v_{1,2}-v_{2,2}
   }\\
   v_{1,1} + {\rm max}(0,v_{1,2})
   \\
   v_{2,1} + {\rm max}(0,v_{1,2},v_{1,2}+v_{2,2})-{\rm max}(0,v_{1,2})
   \end{bmatrix}
   }
\end{equation}
and similarly for $Q^-_{\sigma_1}$ with `max' replaced by `min'. We verify that $Q^+_{\sigma_1}$ maps each ${\bf g}(x)$ to ${\bf g}\bigl(h_{\sigma_1}(x)\bigr)$ and that $Q^-_{\sigma_1}$ maps each ${\bf g}^{\rm op}(x)$ to ${\bf g}^{\rm op}\bigl(h_{\sigma_1}(x)\bigr)$ in Table \ref{Gr36gvecs}. For example, we have
\begin{align}
    Q^+_{\sigma_1} {\bf g}(P_{124}) = (0,0,1,0) = {\bf g}(P_{125}),\\
    \notag
    Q^-_{\sigma_1} {\bf g}^{\rm op}(P_{124}) = (0,0,1,0) = {\bf g}^{\rm op}(P_{125})
\end{align}
and we note that $\tilde{\sigma}_1(P_{124}) = P_{125}P_{234}$ so that $h_{\sigma_1}(P_{124}) = P_{125}$.
We omit the similar analysis for $\sigma_2$. It should be noted that, since $\mathbb{C}[{\rm Gr}(3,6)]$ is a cluster algebra of finite type, the generators $Q^\pm_{\sigma_i}$ necessarily have finite order. Indeed we find $(Q^\pm_{\sigma_1})^4 = {\rm id}$. It is not the case that the braid generators have finite order for general Grassmannians.
\end{example}

\section{g-vectors and tableaux} \label{sec:g vectors and tableaux}
In this section, we show that from any truncated g-vector, one can recover the corresponding semistandard Young tableau. Moreover, we obtain the degree of the polynomial in Pl\"{u}cker coordinates in the dual canonical basis of $\CC[\Gr(k,n)]$ corresponding to any truncated g-vector.  

\subsection{g-vectors} \label{subsec:g vectors of modules and tableaux}
By results in \cite[Section 5.2.2]{HL16}, see also \cite[Section 2.6]{DR20}, \cite[Section 7]{CDFL}, we have that for any simple $U_q(\widehat{\mathfrak{g}})$-module 
$L(M)$, its g-vector can be obtained as follows. The dominant monomial $M$ can be written as $M = \prod_{i,s} Y_{i,s}^{a_{i,s}}$ for some non-negative integers $a_{i,s}$, where the product runs over all fundamental modules $L(Y_{i,s})$ in $\mathcal{C}_{\ell}$. On the other hand, $M$ can also be written as $M = \prod_{M_i} M_i^{g_i}$ for some integers $g_i$, where the product runs over all initial cluster variables and frozen variables $L(M_i)$ in $\mathcal{C}_{\ell}$. The $g_i$'s are the unique solution of $\prod_{i,s} Y_{i,s}^{a_{i,s}} = \prod_{M_i} M_i^{g_i}$. With a chosen order, $g_i$'s form the g-vector of $L(M)$. A simple module $L(M)$ is determined by its g-vector uniquely. 

In the case of Grassmannian cluster algebras, in Section 7 of \cite{CDFL}, it is shown that given any tableau, one can recover its g-vector as follows. Any tableau $T \in \SSYT(k, [n])$ can be written uniquely as $S_1^{e_1} \cup \cdots \cup S_m^{e_m}$ for some integers $e_1, \ldots, e_m \in \ZZ$, where $S_1, \ldots, S_m$ are the tableaux in the initial cluster (we choose an order of the initial cluster variables). The vector $(e_1, \ldots, e_m)$ is the g-vector of $T$.

\subsection{g-vectors and tableaux} \label{subsec:g vectors and tableaux}

We say that a tableau $T$ does not have frozen factors if there are no frozen tableau $T'$ such that $T = T' \cup T''$ for some tableau $T''$.

When we want to emphasise that a g-vector ${\bf g}$ does not have entries corresponding to frozen variables, we call ${\bf g}$ a truncated g-vector. It is the usual g-vector in \cite{FZ07}.  

\begin{lemma} \label{lem: relation of g vectors and tableaux}
For any vector ${\bf g}$ in $\ZZ^{(k-1)(n-k-1)}$, there is a unique tableau $T_{\bf g}$ which does not have frozen factors such that the first $(k-1)(n-k-1)$ components of ${\bf g}(T_{\bf g})$ is ${\bf g}$.
\end{lemma}

\begin{proof}
Let $m=(k-1)(n-k-1)$ and let $M_1, \ldots, M_{m}$ be the monomials corresponding to the cluster variables $a_1, \ldots, a_m$ in the initial seed respectively. We write ${\bf g} = (g_1, \ldots, g_m)$. The monomial $M = M_1^{g_1} \cdots M_m^{g_m}$ can be written as $M = AB^{-1}$ for two dominant monomials $A, B$. The monomial $B$ is of the form $B = \prod_{i \in I} Y_{i, s_1}^{u_{i,1}} \cdots Y_{i, s_{r_i}}^{u_{i,r_i}}$ for some positive integers $r_i$, $u_{i,j}$. Let $B' = \prod_{i \in I} (Y_{i, \xi(i)} \cdots Y_{i, \xi(i)+2 \ell})^{ \max( u_{i, j}: j = 1, \ldots, r_i ) }$. Then $M B'$ is a dominant monomial and $L(MB')$ does not have frozen factors. By Theorem 3.17 in \cite{CDFL}, there is a unique tableau $T_{\bf g}$ which does not have frozen factors corresponding to $L( M B' )$. The procedure to obtain $T_{\bf g}$ from $L( M B' )$ is as follows: for each $Y_{i,s}$ in $MB'$, take the corresponding one-column tableau, and take the $\cup$-product of all these tableaux, we obtain a tableau $T'$. Removing all frozen factors from $T'$, we obtain $T_{\bf g}$. By the above procedure, and the relation of g-vectors and dominant monomials in Section \ref{subsec:correspondence between modules and tableaux}, the relation of g-vectors and tableaux in Section \ref{subsec:g vectors and tableaux}, we have that the first $m$ components of ${\bf g}(T_{\bf g})$ is ${\bf g}$.  
\end{proof}

\begin{example}
Consider the case of $k=4$, $n=8$ and take an initial cluster whose cluster variables and frozen variables are
\begin{align*}
& P_{1235},\ P_{1245},\ P_{1345},\ P_{1236},\ P_{1256},\ P_{1456},\ P_{1237},\ P_{1267},\ P_{1567},\\
& P_{1234},\ P_{2345},\ P_{3456},\ P_{4567},\ P_{5678},\ P_{1238},\ P_{1278},\ P_{1678}
\end{align*}   
The initial cluster variables correspond to the modules with the following dominant monomials respectively:
\begin{align*}
& Y_{1,-1}, \ Y_{2,0}, \  Y_{3,1}, \ Y_{1,-1}Y_{1,-3}, \  Y_{2,0}Y_{2,-2}, \  Y_{3,1}Y_{3,-1}, \ Y_{1,-1}Y_{1,-3}Y_{1,-5},  \ Y_{2,0}Y_{2,-2}Y_{2,-4},  \\ 
& Y_{3,-3}Y_{3,-1}Y_{3,1}, \ 1, \ 1, \ 1, \ 1, \ 1, \ Y_{1,-1}Y_{1,-3}Y_{1,-5}Y_{1,-7}, \ Y_{2,0}Y_{2,-2}Y_{2,-4}Y_{2,-6}, \ Y_{3,1}Y_{3,-1}Y_{3,-3}Y_{3,-5}.
\end{align*} 
Let ${\bf g} = (-1, 0, 0, -1, 0, 1, 1, 0, 0) \in \ZZ^{9}$. The vector ${\bf g}$ produces a monomial 
\begin{align*}
M =Y_{1,-1}Y_{1,-3}Y_{1,-5}Y_{3,-1}Y_{3,1} (Y_{1,-1}^2 Y_{1,-3})^{-1} = Y_{1,-5}Y_{3,-1}Y_{3,1}Y_{1,-1}^{-1} = A B^{-1},
\end{align*}
where $A = Y_{1,-5}Y_{3,-1}Y_{3,1}$, $B = Y_{1,-1}$. We multiply the dominant monomial $B' = Y_{1,-1}Y_{1,-3}Y_{1,-5}Y_{1,-7}$ of the frozen variable $L(Y_{1,-1}Y_{1,-3}Y_{1,-5}Y_{1,-7})$ to $M$ and obtain 
\begin{align*}
MB' = Y_{1,-3}Y_{1,-5}^2Y_{1,-7}Y_{3,-1}Y_{3,1}.
\end{align*}
This module corresponds to the tableau $T_{\bf g} = [[1, 3, 4, 7], [2, 4, 5, 8]]$ (each list is a column). The tableau $T_{\bf g}$ does not have frozen factors and its g-vector is
\begin{align*}
{\bf g}(T_{\bf g}) = (-1, 0, 0, -1, 0, 1, 1, 0, 0, 0, 1, 0, 0, 0, 1, 0, 0). 
\end{align*}
We see that the first $9$ entries of the above g-vector is the same as the given ${\bf g}$. 
\end{example}

\subsection{Two monoids of sets of one-column tableaux} \label{subsec:Two monoids of sets of one-column tableaux}

Denote by ${\rm S}(k, [n])$ the set of all sets of one-column tableaux with $k$ rows and with entries in $[n]$. We define a relation ``$\sim$'' on the set ${\rm S}(k, [n])$: for $S=\{S_1, \ldots, S_c\}$, $T=\{T_1, \ldots, T_{d}\} \in {\rm S}(k, [n])$, $S \sim T$ if $(\cup_i S)_{\rm red} = (\cup_i T_i)_{\rm red}$ (see Section \ref{subsec:Grassmannian cluster algebras} for the notation). Denote by ${\rm S}(k, [n], \sim)$ the set of $\sim$-equivalence classes in ${\rm S}(k, [n])$. It is clear that both of the sets ${\rm S}(k, [n])$ and ${\rm S}(k, [n], \sim)$ are monoids in which the product of two sets $S, T$ is given by $S \cup T$. 

It is easy to see that the monoid ${\rm S}(k, [n], \sim)$ is isomorphic to ${\rm SSYT}(k, [n], \sim)$ via $S=\{S_1, \ldots, S_c\} \mapsto \cup_i S_i$. 

There is a surjective map from the monoid ${\rm S}(k, [n])$ to ${\rm SSYT}(k, [n])$ via $S=\{S_1, \ldots, S_c\} \mapsto \cup_i S_i$. This map is not injective. 

For a given set of one-column tableaux $T=\{T_1, \ldots, T_r\} \in {\rm S}(k, [n])$, we define a dominant monomial $M_T$ corresponding to $T$ by $M_T = M_{ \cup_i T_i }$. It is clear that for two representatives $T, T'$ of an element in ${\rm S}(k, [n], \sim)$, we have $M_T = M_{T'}$. Therefore the map $T \mapsto M_T$ is also well-defined on  ${\rm S}(k, [n], \sim)$. 

The map ${\rm S}(k, [n]) \to \mathcal{P}_{\ell}^+$ ($n=k+\ell+1$) given by $T \mapsto M_T$ is surjective but not injective. The fibres of this map will be used in Section \ref{subsec:degree of g vector and decomposition of dominant monomial} to explain the degree of the polynomial in $\CC[\Gr(k,n)]$ in Pl\"{u}cker coordinates corresponding to a g-vector. 

\subsection{Degree of g-vectors and decompositions of a dominant monomial} \label{subsec:degree of g vector and decomposition of dominant monomial}
In Section \ref{subsec:g vectors and tableaux}, we construct a semistandard tableau $T_{\bf g}$ from a given truncated g-vector ${\bf g}$. In this subsection, we explain the relation between the number of columns of $T_{\bf g}$ and the minimum number of factors of decomposing $M_{T_{\bf g}}$ into increasing minimal affinizations (see the definition of increasing minimal affinizations below). 

We call the number of columns of $T_{\bf g}$ the degree of a given truncated g-vector ${\bf g}$. We say that the polynomial $\ch(T_{\bf g})$ in the dual canonical basis of $\CC[\Gr(k,n)]$ corresponds to the given truncated g-vector $g$. The degree of $g$ is equal to the degree of the polynomial $\ch(T_{\bf g})$ in Pl\"{u}cker coordinates. 
%We expect that it equals to the degree of the polynomial in Pl\"{u}cker coordinates in the dual canonical basis of $\CC[\Gr(k,n)]$ corresponding to $g$.

For a dominant monomial $M$, if $L(M)$ is a minimal affinization of $U_q(\widehat{\mathfrak{g}})$ (see Section \ref{subsec:minimal affinizations}), then we also call $M$ itself a minimal affinization. 
We call a minimal affinization $M$ increasing minimal affinization if we order the first indices of factors $Y_{i,s}$ in $M$ increasingly, then the second indices are also increasing. 

For example, $Y_{1,-3}Y_{2,0}$ is an increasing minimal affinization, and $Y_{2,-4}Y_{1,-1}$ is not an increasing minimal affinization. According to Theorem 3.17 in \cite{CDFL}, an increasing minimal affinization corresponds to a one-column tableau (after removing all trivial factors) but a minimal affinization which is not increasing corresponds to a tableau with more than one column (after removing all trivial factors). 

\begin{definition}
For every dominant monomial $M \in \mathcal{P}^+_{\ell}$, we define a set $P(M)$ of $M$, where each element of $P(M)$ is a decomposition of $M$ into a set of increasing minimal affinizations. 
\end{definition}

The set $P(M)$ is a poset with respect to the following relation: for two decompositions $\mathcal{M}=\{M_1, \ldots, M_r\}$, $\mathcal{M}'=\{M_1', \ldots, M_{r'}'\}$ of $M$, $\mathcal{M} \le \mathcal{M}'$ if there exists $M_t$ in $\mathcal{M}$ such that $M_t = \prod_{j=1}^a M_{i_j}'$ for some $M_{i_j}'$ in $\mathcal{M}'$. 

\begin{lemma} \label{lem:the numbers of terms in minimal elements in PM are the same}
For any dominant monomial $M \in \mathcal{P}^+_{\ell}$, the numbers of terms in minimal elements in $P(M)$ are the same.
\end{lemma}

\begin{proof}

Given two representatives $T=\{T_1, \ldots, T_r\}$, $T'=\{T_1', \ldots, T_{r'}'\}$ of an element in ${\rm S}(k, [n], \sim)$, according to Section \ref{subsec:Two monoids of sets of one-column tableaux}, we have that $(\cup_i T_i)_{\rm red} = (\cup_i T_i')_{\rm red}$ and $M_T = M_{T'}$. For any element $S \in {\rm S}(k, [n], \sim)$, there is a one to one correspondence between representatives of $S$ and different decompositions of $M_S$ into increasing minimal affinizations: for any representative $T= \{T_1, \ldots, T_r\}$ of $S$, the corresponding decomposition of $M_S$ into increasing minimal affinizations is $\{ M_{T_1}, \ldots, M_{T_r} \}$. 

For any two representatives $T, T'$ of an element in ${\rm S}(k, [n], \sim)$ such that $\cup_i T_i$ and $\cup_i T_i'$ do not have frozen factors, we have that $\cup_i T_i$ and $\cup_i T_i'$ have the same content and have the same number of columns. Therefore the numbers of terms in the decompositions of $M_{T}=M_{T'}$ corresponding to $T, T'$ are the same. 
\end{proof}

\begin{example}
Let $M$ be the dominant monomial $Y_{1,-5}Y_{1,-3}Y_{2,-2}Y_{2,0}$. The set $P(M)$ consists of the following decompositions of $M$ into increasing minimal affinizations:
\begin{align*}
& \{ Y_{1,-5}, \ Y_{1,-3}, \ Y_{2,-2}, \ Y_{2,0} \}, \quad \{ Y_{1,-5}Y_{1,-3}, \ Y_{2,-2}, \ Y_{2,0} \}, \quad \{ Y_{1,-5}, \ Y_{1,-3}, \ Y_{2,-2}Y_{2,0} \}, \\
& \{ Y_{1,-5}Y_{2,-2}, \ Y_{1,-3}, \ Y_{2,0} \}, \quad \{ Y_{1,-5}, \ Y_{2,-2}, \ Y_{1,-3}Y_{2,0} \}, \quad \{ Y_{1,-5}Y_{2,-2}, \ Y_{1,-3}Y_{2,0} \}, \\
& \{ Y_{1,-5}Y_{1,-3}, \ Y_{2,-2}Y_{2,0} \}, \quad \{ Y_{1,-5}Y_{2,-2}Y_{2,0}, \ Y_{1,-3} \},  \quad \{ Y_{1,-5}Y_{1,-3}Y_{2,0}, \ Y_{2,-2} \}. 
\end{align*}
The smallest size of the above sets is $2$. This agrees with that the tableau $\vcenter{\hbox{\scalebox{0.6}{ \begin{ytableau}
    1 & 2 \\ 3 & 4 \\ 5 & 6
\end{ytableau} }}}$ in $\SSYT(3, [6])$ corresponding to $M$ has $2$ columns. The tableau in $\SSYT(k, [n])$ ($k \ge 3$) corresponding to $M$ also has $2$ columns. 
\end{example}

By Lemma \ref{lem:the numbers of terms in minimal elements in PM are the same}, we have that
\begin{corollary} \label{corollary:tableau without trivial factors has the number of columns equals to the number of decomposition of M into increasing minimal affinizations}
Given a tableau $T$ without trivial factors, the number of columns of $T$ is equal to the number of terms in any minimal element in $P(M_{T})$.
\end{corollary}

\begin{remark}
A tableau without frozen factors does not have trivial factors. Therefore Corollary \ref{corollary:tableau without trivial factors has the number of columns equals to the number of decomposition of M into increasing minimal affinizations} works for tableaux without frozen factors. 
\end{remark}

\section{Quasi-automorphism actions on tableaux} \label{sec:braid group actions on tableaux}

In this section, we study braid group action (resp. cyclic rotation map, twist map, reflection map) on $\SSYT(k,[n])$. 

In Section \ref{sec:braid group actions on g vectors}, we defined the braid group action (resp. cyclic rotation map, twist map, reflection map) on (truncated) g-vectors for $\CC[\Gr(k,n)]$. In Section \ref{sec:g vectors and tableaux}, we explained the relation between g-vectors and tableaux.  

We denote $f({\bf g})=Q_f^+({\bf g})$ for any quasi-automorphism $f$, where $Q_f^+$ is in (\ref{Qfmaps}).

\begin{definition} \label{def:quasi-automorphisms on tableaux}
For any tableau $T \in \SSYT(k,[n])$ without frozen factors and any quasi-automorphism on $\CC[\Gr(k,n)]$, we define $f(T) = T_{f({\bf g}_T)}$, where ${\bf g}_T$ is the truncated g-vector of $T$, and $T_{f({\bf g}_T)}$ is the tableau without frozen factors whose truncated g-vector is $f({\bf g}_T)$. 
\end{definition}
For example, in the above definition, $f$ can be a composition of $\sigma, \rho, \theta, \tau$. 

\subsection{Braid group action on tableaux}

\begin{example} \label{example:Gr36 tableau}
We continue the example of $\Gr(3,6)$ in Example \ref{example:Gr36}. Recall that we use the order $P_{124}$, $P_{134}$, $P_{125}$, $P_{145}$ of initial cluster variables of $\CC[\Gr(3,6)]$ when we write a g-vector. The tableau corresponding to the Pl\"{u}cker coordinate $P_{135}$ is $T=\vcenter{\hbox{\scalebox{0.6}{\begin{ytableau}
    1 \\ 3 \\ 5
\end{ytableau}}}}$. Its truncated g-vector is ${\bf g}=(-1, 1, 1, 0)$, see Table \ref{Gr36gvecs}. We have $\sigma_1({\bf g}) = (-1, 0, 0, 1)$. The tableau corresponding to $\sigma_1({\bf g})$ is $\vcenter{\hbox{\scalebox{0.6}{\begin{ytableau}
    1 & 2 \\ 3 & 4 \\ 5 & 6
\end{ytableau}}}}$. Therefore $\sigma_1(T) = \vcenter{\hbox{\scalebox{0.6}{\begin{ytableau}
    1 & 2 \\ 3 & 4 \\ 5 & 6
\end{ytableau}}}}$. 
We have $\sigma_1(\ch(T)) = \sigma_1(P_{135}) = P_{145}P_{236}-P_{123}P_{456}$ (after removing frozen factors). On the other hand, $\ch(\sigma_1(T)) = \ch( \vcenter{\hbox{\scalebox{0.6}{\begin{ytableau}
    1 & 2 \\ 3 & 4 \\ 5 & 6
\end{ytableau}}}} ) = P_{145}P_{236}-P_{123}P_{456}$. Therefore $\sigma_1(\ch(T)) = \ch(\sigma_1(T))$. 
\end{example}

\begin{example}
We consider the case of $\CC[\Gr(3,6)]$ which is a finite type cluster algebra. We describe the braid group action on cluster variables of $\CC[\Gr(3,6)]$ in terms of tableaux in Table \ref{table:braid group action on Gr36 in terms of tableaux}. There are two orbits of cluster variables in $\CC[\Gr(3,6)]$:
\begin{align*}
\scalemath{0.7}{\begin{ytableau}
1 \\
2 \\
4
\end{ytableau}, \ \begin{ytableau}
1 \\
4 \\
6
\end{ytableau}, \ \begin{ytableau}
2 \\
3 \\
5
\end{ytableau}, \ \begin{ytableau}
1 \\
3 \\
4
\end{ytableau}, \ \begin{ytableau}
1 \\
4 \\
5
\end{ytableau}, \ \begin{ytableau}
1 \\
2 \\
5
\end{ytableau}, \ \begin{ytableau}
3 \\
5 \\
6
\end{ytableau}, \ \begin{ytableau}
2 \\
5 \\
6
\end{ytableau}, \ \begin{ytableau}
2 \\
3 \\
6
\end{ytableau}, \ \begin{ytableau}
3 \\
4 \\
6
\end{ytableau}, \ \begin{ytableau}
2 \\
4 \\
5
\end{ytableau}, \ \begin{ytableau}
1 \\
3 \\
6
\end{ytableau}, }
\end{align*}
and
\begin{align*}
\scalebox{0.7}{
\begin{ytableau}
1 \\
3 \\
5
\end{ytableau}, \ \begin{ytableau}
1 & 3 \\
2 & 5 \\
4 & 6
\end{ytableau}, \ \begin{ytableau}
2 \\
4 \\
6
\end{ytableau}, \ \begin{ytableau}
1 & 2 \\
3 & 4 \\
5 & 6
\end{ytableau}. }
\end{align*}

We describe the braid group action on cluster variables in Table \ref{table:braid group action on Gr36 in terms of tableaux}. This tableau agrees with Table 1 in \cite{CHL} which is the quantum version of $\sigma_i$.
\end{example}

\begin{table}
\begin{tabular}{|c|c|c|c|c|c|c|c|c|}
\hline
tableau $x$ & $\sigma_1(x)$ & $\sigma_2(x)$ \\
\hline
 $[1,2,4]$ & $[1,2,5]$  & $[1, 3, 4]$ \\
\hline 
 $[1,2,5]$ & $[1,4,5]$ & $[1, 3, 6]$ \\
\hline
 $[1,3,4]$ & $[2,3,5]$ & $[1, 4, 5]$  \\
\hline
 $[1,3,5]$ & $[[1,3,5],[2,4,6]]$ & $[[1, 3, 5], [2, 4, 6]]$   \\
\hline
 $[1,3,6]$ & $[2,3,6]$ & $[2, 4, 5]$  \\
\hline
 $[1,4,5]$ & $[2,4,5]$ &  $[1, 4, 6]$  \\
\hline
 $[1,4,6]$ & $[2,5,6]$ &  $[1, 2, 4]$  \\ 
\hline
 $[2,3,5]$ & $[1,4,6]$ & $[2, 3, 6]$  \\
\hline
 $[2,3,6]$ & $[3,4,6]$ & $[2, 5, 6]$  \\
\hline
 $[2,4,5]$ & $[1,2,4]$ &  $[3, 4, 6]$  \\
\hline
 $[2,4,6]$ & $[[1,2,4],[3,5,6]]$ & $[[1, 2, 4], [3, 5, 6]]$  \\
\hline
 $[2,5,6]$ & $[1,3,4]$ & $[3, 5, 6]$  \\
\hline
 $[3,4,6]$ & $[3,5,6]$ & $[1,2,5]$   \\
\hline
 $[3,5,6]$ & $[1,3,6]$ & $[2, 3, 5]$  \\
\hline
$[1,2,4], [3,5,6]$ & $[1,3,5]$ &  $[1,3,5]$  \\
\hline
$[1,3,5],[2,4,6]$ & $[2,4,6]$ & $[2, 4, 6]$ \\
\hline
\end{tabular}
\caption{Braid group action on $\CC[\Gr(3,6)]$ in terms of tableaux. In the table, $[[1,3,5],[2,4,6]]$ means the tableau with columns $[1,3,5]$ and $[2,4,6]$.}
\label{table:braid group action on Gr36 in terms of tableaux}
\end{table}

\subsection{Tropical cyclic rotations and promotions of tableaux} \label{sec:Tropical cyclic rotations and promotions of tableaux}
In this subsection, we study tropical cyclic rotations on $\CC[\Gr(k,n)]$ and their relation with promotions of tableaux. 

Promotion is an operator on the set of semistandard Young tableaux defined in terms of ``jeu de taquin'' sliding moves \cite{Sch63, Sch72, Sch77}. Gansner \cite{Ga80} proved that promotion can also be described using using Bender-Knuth involutions \cite{BK}. 
In this paper, we only need the promotion operator on $\SSYT(k, [n])$. 

The $i$th ($i \in [n]$) Bender-Knuth involution ${\rm BK}_i: \SSYT(k, [n]) \to \SSYT(k, [n])$, is defined by the following procedure: for $i$ and $i+1$ which are not in the same column, we replace $i$ by $i+1$ and replace $i+1$ by $i$, and then reorder $i$, $i+1$ in each row such that the resulting tableau is semistandard. The promotion ${\rm pr}(T)$ of $T$ is defined by 
\begin{align*}
{\rm pr}(T) = {\rm BK}_1 \circ \cdots \circ {\rm BK}_{n-1}(T),    
\end{align*}
see also Definition A.3 in \cite{Hop20}.  

We conjecture the following.
\begin{conjecture}
For any tableau $T \in \SSYT(k,[n])$ without frozen factors, $\rho(T)$ is equal to ${\rm pr}(T)$. 
\end{conjecture}

\begin{example}  
We consider the Grassmannian cluster algebra $\CC[\Gr(3,6)]$ in Example \ref{example:Gr36}. The tableau corresponding to the Pl\"{u}cker coordinate $P_{124}$ is $T=\vcenter{\hbox{\scalebox{0.6}{\begin{ytableau}
    1 \\ 2 \\ 4
\end{ytableau}}}}$. Its truncated g-vector is ${\bf g}=(1, 0, 0, 0)$. 
We have $\rho({\bf g}) = (-1, 0, 1, 0)$. The tableau corresponding to $\rho({\bf g})$ is $\vcenter{\hbox{\scalebox{0.6}{\begin{ytableau}
   2 \\ 3 \\ 5
\end{ytableau}}}}$. It is exactly ${\rm Pr}(T)$. 
\end{example}

\subsection{Tropical twist map}\label{sec:tropical Marsh Scott twist map}

In this subsection, we study a tropical version of the twist map on Grassmannian cluster algebras. 

The twist map $\tau$ on $\CC[\Gr(k,n)]$ is periodic up to removing frozen variables, see Section 4 in \cite{MS16}. Therefore $\tau$ on truncated g-vectors and on $\SSYT(k,[n])$ is periodic. Therefore there exists a positive integer $d$ such that $\tau^d(T)=T$. We have $\tau^{-1}(T) = \tau^{d-1}(T)$. 

\begin{example}
We consider the Grassmannian cluster algebra $\CC[\Gr(3,6)]$ in Example \ref{example:Gr36}. The tableau corresponding to the Pl\"{u}cker coordinate $P_{124}$ is $T=\vcenter{\hbox{\scalebox{0.6}{\begin{ytableau}
    1 \\ 2 \\ 4
\end{ytableau}}}}$. Its truncated g-vector is ${\bf g}=(1, 0, 0, 0)$. 
We have that $\tau({\bf g}) = (-1, 0, 0, 0)$. The g-vector $\tau({\bf g})$ is the g-vector of the cluster variable $P_{236}$ and $P_{236} = \tau(P_{124})$. The tableau corresponding to $\tau({\bf g})$ is $\vcenter{\hbox{\scalebox{0.6}{\begin{ytableau}
   2 \\ 3 \\ 6
\end{ytableau}}}}$. Therefore $\tau(T) = \vcenter{\hbox{\scalebox{0.6}{\begin{ytableau}
   2 \\ 3 \\ 6
\end{ytableau}}}}$. 
\end{example}

\subsection{Tropical reflection maps and evacuations of tableaux} \label{sec:tropical reflection map}

In this subsection, we study tropical reflection maps on $\CC[\Gr(k,n)]$ and their relation with evacuations of tableaux. 

Using tropicalisation in Section \ref{sec:tropicalisation of quasi homomorphisms}, we have a map $\theta$ sending g-vectors of cluster variables to some integer vectors. Recall that we denote $Q_f^+({\bf g})=f({\bf g})$ for any quasi-automorphism $f$ and any g-vector ${\bf g}$. We have the conjecture. 
\begin{conjecture} \label{conj:theta g is g of theta of tau inverse of b}
For any cluster variable $b$ in $\CC[\Gr(k,n)]$, we have that 
\begin{align*}
& \theta({\bf g}(b)) = {\bf g}( \theta( \tau^{-1}(b) ) ), \ {\bf g}(\theta(\tau^{-1}(b))) = - {\bf g}^{\rm op}(\theta(b)), \ {\bf g}^{\rm op}(a) = - {\bf g}(\theta(\tau^{-1}(\theta(a)))), \\
& {\bf g}(c) = - {\bf g}^{\rm op}(\theta \tau \theta (c)), \ {\bf g}^{\rm op}(\theta(c)) = - {\bf g}^{\rm op}(\tau^{-2} c)^\pi, \ {\bf g}(\theta(b)) = - {\bf g}(b)^\pi, \ \theta\tau({\bf g}(b)) = - {\bf g}(b)^\pi,
\end{align*} 
where $\pi$ is the permutation which reverses the order of the entries of a g-vector.
\end{conjecture}
We explain that the first identity in Conjecture \ref{conj:theta g is g of theta of tau inverse of b} implies all other identities in the conjecture. By (\ref{gvecmapsGrass2}), Conjecture \ref{conj:theta g is g of theta of tau inverse of b} is equivalent to: for any cluster variable $b$ in $\CC[\Gr(k,n)]$, we have ${\bf g}(\theta(\tau^{-1}(b))) = - {\bf g}^{\rm op}(\theta(b))$. Set $a=\theta(b)$ and use the fact that $\theta$ is an inovolution, we see that this conjecture is also equivalent to: for any cluster variable $a$ in $\CC[\Gr(k,n)]$, ${\bf g}^{\rm op}(a) = - {\bf g}(\theta(\tau^{-1}(\theta(a))))$. Set $a = \theta \tau \theta (c)$, from ${\bf g}^{\rm op}(a) = - {\bf g}(\theta(\tau^{-1}(\theta(a))))$, we obtain ${\bf g}(c) = - {\bf g}^{\rm op}(\theta \tau \theta (c))$. By the second equation in (\ref{thetataugmaps}), ${\bf g}^{\rm op}\bigl(\theta(\tau(x))\bigr) = {\bf g}(x)^\pi$, where we write $h_{\tau}$ as $\tau$ for simplicity. Let $x= \theta (c)$. Then ${\bf g}(c) = - {\bf g}(\theta (c))^\pi$.  That is, ${\bf g}(\theta(c)) = - {\bf g}(c)^\pi$. By ${\bf g}^{\rm op}(a) = - {\bf g}(\theta \tau^{-1} \theta (a)) = - {\bf g}(\theta \tau \tau^{-2} \theta (a))$ and the first equation in (\ref{thetataugmaps}): ${\bf g}(\theta \tau (x)) = {\bf g}^{\rm op}(x)^\pi$, we obtain ${\bf g}^{\rm op}(a) = -{\bf g}^{\rm op}(\tau^{-2} \theta(a))^\pi$. Set $a = \theta(c)$, we obtain ${\bf g}^{\rm op}(\theta(c)) = - {\bf g}^{\rm op}(\tau^{-2} c)^\pi$. Let $b = \tau(a)$. By $\theta({\bf g}(b)) = {\bf g}(\theta(\tau^{-1}(b)))$, we have that $\theta({\bf g}(\tau(a))) = {\bf g}(\theta(a))$. Therefore $\theta\tau({\bf g}(a)) = {\bf g}(\theta(a))$. By ${\bf g}(\theta(a)) = - {\bf g}(a)^\pi$, we have that $\theta\tau({\bf g}(a)) = - {\bf g}(a)^\pi$. 

Conjecture \ref{conj:theta g is g of theta of tau inverse of b} has an interesting corollary: there exists a seed whose set of 
g-vectors consists precisely of the negatives of the standard unit vectors (recall that the 
g-vectors of the initial cluster variables are the unit vectors). Indeed, starting from an initial seed and applying $\theta \circ \tau$ to the g-vectors of the initial cluster variables produces such a seed, as follows from the identity $\theta\tau({\bf g}(b)) = - {\bf g}(b)^\pi$.

\begin{example}
Let $b=P_{125}$ in $\CC[\Gr(3,6)]$. Then ${\bf g}(b)=(0,0,1,0)={\bf g}^{\rm op}(b)$, $\theta(b)=P_{256}$, $\tau(b)=P_{346}$, $\tau^{-1}(b)=P_{136}$, $\theta(\tau(b))=P_{134}$, $\theta(\tau^{-1}(b))=P_{146}$, ${\bf g}(P_{134})=(0,1,0,0)$, ${\bf g}(P_{146})=(0,0,-1,1)$,
$Q_{\theta}^+({\bf g}(b)) = (0,0,-1,1) = {\bf g}(\theta(\tau^{-1}(b)))$.  
\end{example}

Evacuation of tableaux is an involution which was introduced by Sch\"{u}tzenberger \cite{Sch63, Sch72}. It can be defined using Bender-Knuth moves (see Section \ref{sec:Tropical cyclic rotations and promotions of tableaux}) as follows: for $T \in \SSYT(k,n)$, the evacuation of $T$ is 
\begin{align*}
\scalemath{0.95}{{\rm eva}(T) = {\rm BK}_1 \circ ( {\rm BK}_2 \circ {\rm BK}_1 ) \circ \cdots \circ ( {\rm BK}_{n-2} \circ \cdots \circ {\rm BK}_{2} \circ {\rm BK}_1 ) \circ ( {\rm BK}_{n-1} \circ \cdots \circ {\rm BK}_{2} \circ {\rm BK}_1 )(T), }
\end{align*}
see also Definition A.3 in \cite{Hop20}. 

We conjecture the following. 
\begin{conjecture} \label{conjecture:thetaT_is_eva_tau_inverse_T}
For any $T \in \SSYT(k,[n])$, we have that $\theta( \ch(T)) = \ch({\rm eva}(T))$ and $\theta(T) = {\rm eva}(\tau^{-1}(T))$.
\end{conjecture}

\begin{example}
We consider the Grassmannian cluster algebra $\CC[\Gr(3,6)]$ in Example \ref{example:Gr36}. The tableau corresponding to the Pl\"{u}cker coordinate $P_{125}$ is $T=\vcenter{\hbox{\scalebox{0.6}{\begin{ytableau}
    1 \\ 2 \\ 5
\end{ytableau}}}}$. We have $\ch(T)=P_{125}$, $\theta(P_{125}) = P_{256}$, ${\rm eva}(T) = \vcenter{\hbox{\scalebox{0.6}{\begin{ytableau}
    2 \\ 5 \\ 6
\end{ytableau}}}}$. Therefore $\theta( \ch(T)) = \ch({\rm eva}(T))$.
We have ${\bf g} = {\bf g}(P_{125})=(0,0,1,0)$ and $\theta({\bf g}) = Q_{\theta}^+({\bf g}) = (0, 0, -1, 1)$. The vector $Q_{\theta}^+({\bf g})$ is the g-vector of $\vcenter{\hbox{\scalebox{0.6}{\begin{ytableau}
   1 \\ 4 \\ 6
\end{ytableau}}}}$. Therefore $\theta(T)=\vcenter{\hbox{\scalebox{0.6}{\begin{ytableau}
   1 \\ 4 \\ 6
\end{ytableau}}}}$. On the other hand, $\tau^{-1}(T) = \vcenter{\hbox{\scalebox{0.6}{\begin{ytableau}
   1 \\ 3 \\ 6
\end{ytableau}}}}$ and ${\rm eva}(\tau^{-1}(T)) = \vcenter{\hbox{\scalebox{0.6}{\begin{ytableau}
   1 \\ 4 \\ 6
\end{ytableau}}}}$. Therefore $\theta(T) = {\rm eva}(\tau^{-1}(T))$. 
\end{example}

We can define another version of $\theta$ on tableaux, denote it by $\widetilde{\theta}(T)$, as $\widetilde{\theta}(T) = T_{ - ( {\bf g}_T )^{\pi}} $, where the permutation $\pi$ is the reverse of the order of entries in a g-vector. We now explain that under the assumption that Conjecture \ref{conj:theta g is g of theta of tau inverse of b} is true, the second statement of Conjecture \ref{conjecture:thetaT_is_eva_tau_inverse_T} is equivalent to the following conjecture.
\begin{conjecture} \label{conj:conjecture for another version of theta on T}
For any $T \in \SSYT(k,[n])$, we have that $\widetilde{\theta}(T) = {\rm eva}(T)$. 
\end{conjecture}
According to Conjecture \ref{conj:theta g is g of theta of tau inverse of b}, $\theta \tau ({\bf g}) = - {\bf g}^{\pi}$. Therefore $\widetilde{\theta}(T)=T_{-({\bf g}_T)^{\pi}} = T_{\theta \tau ({\bf g}_T)} = \theta \tau(T)$.  Conjecture \ref{conjecture:thetaT_is_eva_tau_inverse_T} says that $\theta(T) = {\rm eva}( \tau^{-1}(T) )={\rm eva}(T_{\tau^{-1}({\bf g}_T)})$. Replace $T$ by $\tau(T)$, we obtain $\theta \tau (T) = {\rm eva}( T_{\tau^{-1}({\bf g}_{\tau(T)})}) = {\rm eva} (T_{{\bf g}_T})$ (note that the first $T$ in ${\rm eva}(T_{\tau^{-1}({\bf g}_T)})$ is just the notation for tableau and it is not a variable, so when we replace $T$ by $\tau(T)$, the first $T$ in ${\rm eva}(T_{\tau^{-1}({\bf g}_T)})$ does not change), we see that the second statement of Conjecture \ref{conjecture:thetaT_is_eva_tau_inverse_T} is equivalent to Conjecture \ref{conj:conjecture for another version of theta on T}.

\begin{example}
Let $T = \vcenter{\hbox{\scalebox{0.6}{\begin{ytableau}
    1&2&3\\ 4&4&5\\6&7&8
\end{ytableau}}}} \in \SSYT(3,[8])$. Then ${\bf g} = {\bf g}(T) = (0, -1, -2, 1, 0, 1, 1, 0)$ and $-{\bf g}^{\pi} = (0, -1, -1, 0, -1, 2, 1, 0)$. We have that $T_{-{\bf g}^{\pi}} = \vcenter{\hbox{\scalebox{0.6}{\begin{ytableau}
    1&2&3\\ 4&5&5\\6&7&8
\end{ytableau}}}} = {\rm eva}(T)$.
\end{example}  

\section{Fixed points of quasi-automorphisms} \label{sec:fixed points of quasi-automorphisms}
In this section, we study fixed points of quasi-homomorphisms of cluster algebras. We introduce the concept of stable and unstable fixed points of quas-automorphisms. Stable fixed points give a way to construct prime non-real modules. 

\subsection{Fixed points of quasi-automorphisms}
We recall the definition of a limit g-vector, see Section 11 in \cite{EL23}. For a vector $v = (v_1, \ldots, v_m)$ in $\RR^{m}$, denote its $l^2$-norm by $ \left\Vert v \right\Vert = \sqrt{\sum_{i=1}^m |v_i|^2}$. For a cluster algebra $\mathcal{A}$ of infinite type of rank $m$, we say that a sequence of (distinct) g-vectors ${\bf g}_1, {\bf g}_2, \ldots$ of $\mathcal{A}$ has a limit ${\bf g}$ if the greatest common factor of entries of ${\bf g}$ is $1$ and for every $\epsilon >0$, there is a positive integer $N$ such that for every $j \ge N$, there is some positive real number $c_j$ such that $\left\Vert c_j{\bf g}-{\bf g}_j \right\Vert < \epsilon$. 

\begin{definition} \label{def:stable and unstable fixed points} 
Let $\mathcal{A}$ be any cluster algebra of rank $n$ and $f$ is a quasi-automorphism on $\mathcal{A}$. We say that a g-vector ${\bf g} \in \ZZ^n$ is a fixed point of $f$ if $f({\bf g})={\bf g}$. For a g-vector ${\bf g}$ which is fixed by a quasi-automorphism $f$ on $\mathcal{A}$, we say that ${\bf g}$ is a stable fixed point if for every generic vector ${\bf g}'$ in $\ZZ^n$, the sequence $\sigma^j({\bf g}')$, $j=1,2,\ldots$, has a limit and the limit is $\frac{1}{r}{\bf g}$ for some $r \in \ZZ_{\ge 1}$. Otherwise, we say that ${\bf g}$ is an unstable fixed point.  
\end{definition}

Note that it is possible that ${\bf g}$ is a tropical point but not the g-vector of a cluster variable, see Remark \ref{remark:tropical points and g-vectors}. 

Let $d = \gcd(k,n)$ and $n=k+\ell+1$. We say that a $U_q(\widehat{\mathfrak{sl}_k})$-module $L(M)$ in $\mathcal{C}_{\ell}$ (resp. a tableau $T \in \SSYT(k,[n])$, $\ch(T)$ in $\CC[\Gr(k,n)]$) without frozen factors is a fixed point of some $\sigma \in {\rm Br}_d$ if the truncated g-vector of $L(M)$ (resp. $T$, $\ch(T)$) is fixed under $\sigma$. Moreover, we say that $L(M)$ (resp. $T$, $\ch(T)$) is a stable fixed point if the truncated g-vector of $L(M)$ (resp. $T$, $\ch(T)$) is a stable fixed point under some $\sigma \in {\rm Br}_d$. 

Recall that an element $\ch(T)$ in the dual canonical basis of $\CC[\Gr(k,n)]$ is prime if $\ch(T) \ne \ch(T')\ch(T'')$ for any non-empty tableaux $T',T''$. We have the following conjecture.
\begin{conjecture} \label{conj:stable fixed points and unstable fixed points prime elements and cluster variables}
Suppose that $d=\gcd(k,n)>1$ and $\CC[\Gr(k,n)]$ is of infinite type. Then for every $i \in [d-1]$, every stable fixed point of the generator $\sigma_i$ of ${\rm Br}_d$ in $\CC[\Gr(k,n)]$ is prime non-real. 
\end{conjecture}

% \begin{conjecture} \label{conj:stable fixed points and unstable fixed points prime elements and cluster variables}
% Suppose that $d=\gcd(k,n)>1$ and $\CC[\Gr(k,n)]$ is of infinite type. 
% \begin{enumerate}
% \item For any $i \in [d-1]$ and $(k,n)$ is $(3,9)$ or $(4,8)$, every unstable fixed point of the generator $\sigma_i$ of ${\rm Br}_d$ in $\CC[\Gr(k,n)]$ is a cluster monomial.  

% \item For every $i \in [d-1]$, every stable fixed point of the generator $\sigma_i$ of ${\rm Br}_d$ in $\CC[\Gr(k,n)]$ is prime non-real. 
% \end{enumerate}

% \end{conjecture}

% \begin{remark}
% It is a natural question to ask if (1) of Conjecture \ref{conj:stable fixed points and unstable fixed points prime elements and cluster variables} is true for any $(k,n)$. In Conjecture 3.2 of \cite{CDHHHL}, it is conjecture that if an element in $\CC[\Gr(k,n)]$ is a cluster variable, then it is a cluster variable in $\CC[\Gr(k,n')]$ for any $n' \ge n$. We conjecture that for any element in $\CC[\Gr(k,n)]$, we have that it is a cluster monomial $\CC[\Gr(k,n)]$ if and only if it is a cluster monomial in $\CC[\Gr(k,n')]$ for any $n' \ge n$. For any $k < n$, we choose a smallest $n'$ such that $n \le n'$ and $k \mid n'$. Suppose that $T \in \SSYT(k,[n]) \subset \SSYT(k,[n'])$ is an unstable fixed point of some $\sigma \in {\rm Br}_{k}$. If Conjecture 3.2 in \cite{CDHHHL} is true and (1) of Conjecture \ref{conj:stable fixed points and unstable fixed points prime elements and cluster variables} is true for any $(k,n)$, then $\ch(T)$ is a cluster monomial for $\CC[\Gr(k,n)]$.  

% \end{remark}

\section{\texorpdfstring{Unstable fixed points in $\CC[\Gr(k,n)]$}{Unstable fixed points in C[Gr(k,n)]}}
\label{sec:unstable fixed points and cluster monomials}

In this section, we study unstable fixed points of braid group actions on Grassmannian cluster algebras.  

\subsection{\texorpdfstring{Unstable fixed points in $\CC[\Gr(3,9)]$}{Unstable fixed points in C[Gr(3,9)]}}

Note that $\rho^{-1} \sigma_i \rho = \sigma_{i+1}$, $i \in [d-2]$. It suffices to study $\sigma_1$. By direct computation, there are $4$ unstable fixed points under $\sigma_1$ in $\SSYT(3,[9])$ with $1$ column:
\begin{align*}
\scalebox{0.7}{
\begin{ytableau}
1 \\
2 \\
6
\end{ytableau}, \ \begin{ytableau}
3 \\
6 \\
9
\end{ytableau}, \ \begin{ytableau}
3 \\
7 \\
8
\end{ytableau}, \ \begin{ytableau}
4 \\
5 \\
9
\end{ytableau},}
\end{align*}
All of them are prime and they correspond to Pl\"{u}cker coordinates and they are cluster variables. Their g-vectors are 
\begin{align*}
& (0, 0, 0, 0, 1, 0, 0, 0, 0, 0), \ 
(0, 0, 0, -1, 0, 0, -1, 1, 0, 0), \\
& 
(0, 0, 0, -1, 0, 0, 0, 0, 0, 1), \ 
(0, 0, 0, 0, -1, 0, 0, 0, 0, 0).
\end{align*}
respectively. 

The following are all fixed points of $\sigma_1$ in $\SSYT(3,[9])$ with $2$ columns:
\begin{align} \label{eq:unstable fixed point with 2 columns in ssyt under sigma1 Gr39}
\scalebox{0.7}{
\begin{ytableau}
1 & 1 \\
2 & 2 \\
6 & 6
\end{ytableau},  \ \begin{ytableau}
3 & 3 \\
6 & 6 \\
9 & 9
\end{ytableau},  \ \begin{ytableau}
3 & 3 \\
7 & 7 \\
8 & 8
\end{ytableau},  \ \begin{ytableau}
4 & 4 \\
5 & 5 \\
9 & 9
\end{ytableau}, \ \begin{ytableau}
1 & 3 \\
2 & 6 \\
6 & 9
\end{ytableau},   \ \begin{ytableau}
3 & 3 \\
6 & 7 \\
8 & 9
\end{ytableau}, \ \begin{ytableau}
3 & 4 \\
5 & 6 \\
9 & 9
\end{ytableau}, \ \begin{ytableau}
1 & 4 \\
2 & 6 \\
5 & 9
\end{ytableau}, \ \begin{ytableau}
1 & 4 \\
2 & 7 \\
5 & 8
\end{ytableau}, \ \begin{ytableau}
1 & 3 \\
2 & 7 \\
6 & 8
\end{ytableau}, \ \begin{ytableau}
3 & 4 \\
5 & 7 \\
8 & 9
\end{ytableau}.
}
\end{align}

All of these fixed points are unstable. The last $4$ tableaux in (\ref{eq:unstable fixed point with 2 columns in ssyt under sigma1 Gr39}) are prime. We also checked that these $4$ tableaux are cluster variables. The first $7$ tableaux in (\ref{eq:unstable fixed point with 2 columns in ssyt under sigma1 Gr39}) are non-prime and they correspond to the cluster monomials $P_{126}^2$, $P_{369}^2$, $P_{378}^2$, $P_{459}^2$, $P_{126}P_{369}$, $P_{369}P_{378}$, $P_{369}P_{459}$ respectively. 

The following are all fixed points of $\sigma_1$ in $\SSYT(3,[9])$ with $3$ columns:
\begin{align} \label{eq:unstable fixed point with 3 columns in ssyt under sigma1 Gr39}
\begin{split}
& \scalebox{0.7}{
\begin{ytableau}
1 & 1 & 1 \\
2 & 2 & 2 \\
6 & 6 & 6
\end{ytableau}, \ \begin{ytableau}
1 & 1 & 3 \\
2 & 2 & 6 \\
6 & 6 & 9
\end{ytableau}, \ \begin{ytableau}
1 & 1 & 3 \\
2 & 2 & 7 \\
6 & 6 & 8
\end{ytableau}, \ \begin{ytableau}
1 & 1 & 4 \\
2 & 2 & 6 \\
5 & 6 & 9
\end{ytableau}, \ \begin{ytableau}
1 & 1 & 4 \\
2 & 2 & 7 \\
5 & 6 & 8
\end{ytableau}, \ \begin{ytableau}
1 & 3 & 3 \\
2 & 6 & 6 \\
6 & 9 & 9
\end{ytableau}, \ \begin{ytableau}
1 & 3 & 3 \\
2 & 6 & 7 \\
6 & 8 & 9
\end{ytableau}, } \\
& \scalebox{0.7}{ \begin{ytableau}
1 & 3 & 3 \\
2 & 7 & 7 \\
6 & 8 & 8
\end{ytableau}, \
\begin{ytableau}
1 & 3 & 4 \\
2 & 6 & 6 \\
5 & 9 & 9
\end{ytableau}, \ \begin{ytableau}
1 & 3 & 4 \\
2 & 7 & 7 \\
5 & 8 & 8
\end{ytableau},  \ \begin{ytableau}
1 & 4 & 4 \\
2 & 5 & 6 \\
5 & 9 & 9
\end{ytableau}, \ \begin{ytableau}
1 & 4 & 4 \\
2 & 5 & 7 \\
5 & 8 & 9
\end{ytableau}, \ \begin{ytableau}
3 & 3 & 3 \\
6 & 6 & 6 \\
9 & 9 & 9
\end{ytableau}, \ \begin{ytableau}
3 & 3 & 3 \\
6 & 6 & 7 \\
8 & 9 & 9
\end{ytableau}, \ \begin{ytableau}
3 & 3 & 3 \\
6 & 7 & 7 \\
8 & 8 & 9
\end{ytableau}, } \\
& \scalebox{0.7}{ \begin{ytableau}
3 & 3 & 3 \\
7 & 7 & 7 \\
8 & 8 & 8
\end{ytableau}, \ \begin{ytableau}
3 & 3 & 4 \\
5 & 6 & 6 \\
9 & 9 & 9
\end{ytableau}, \ \begin{ytableau}
3 & 3 & 4 \\
5 & 6 & 7 \\
8 & 9 & 9
\end{ytableau}, \ \begin{ytableau}
3 & 3 & 4 \\
5 & 7 & 7 \\
8 & 8 & 9
\end{ytableau}, \ \begin{ytableau}
3 & 4 & 4 \\
5 & 5 & 6 \\
9 & 9 & 9
\end{ytableau}, \ \begin{ytableau}
3 & 4 & 4 \\
5 & 5 & 7 \\
8 & 9 & 9
\end{ytableau}, \ \begin{ytableau}
4 & 4 & 4 \\
5 & 5 & 5 \\
9 & 9 & 9
\end{ytableau}, \ \begin{ytableau}
1 & 3 & 4 \\
2 & 6 & 7 \\
5 & 8 & 9
\end{ytableau}.
}
\end{split}
\end{align}

The last tableau in (\ref{eq:unstable fixed point with 3 columns in ssyt under sigma1 Gr39}) is a stable fixed point. It is non-real and is one of the tableaux in Example 8.1 in \cite{CDFL}. The other tableaux in (\ref{eq:unstable fixed point with 3 columns in ssyt under sigma1 Gr39}) are non-prime and they correspond to the following cluster monomials respectively:
\begin{align*}
& P_{126}^3, \ P_{126}^2P_{369}, \ P_{126} \ch_{126,378}, \ P_{126}\ch_{125,469}, \ P_{126}\ch_{125,478}, \ P_{126}P_{369}^2, \ P_{369} \ch_{126,378}, \\
& P_{378}\ch_{126,378}, \ P_{369}\ch_{125,469}, \ P_{378}\ch_{125,478}, \ P_{459}\ch_{125,469}, \ P_{459}\ch_{125,478}, \ P_{369}^3, \ P_{369}^2P_{378}, \ P_{369}P_{378}^2, \\
& P_{378}^3, \ P_{369}^2P_{459}, \ P_{369}\ch_{358, 479}, \ P_{378}\ch_{358, 479}, \ P_{369}P_{459}^2, \ P_{459}\ch_{358, 479}, \ P_{459}^3,
\end{align*}
where $\ch_{126,378}$ is the cluster variable $P_{126}P_{378}-P_{123}P_{678}$ corresponding the tableau $\vcenter{\hbox{\scalebox{0.6}{\begin{ytableau}
    1&3\\ 2&7\\6&8
\end{ytableau}}}}$. 

\subsection{\texorpdfstring{Unstable fixed points in $\CC[\Gr(4,8)]$}{Unstable fixed points in C[Gr(4,8)]}}

Note that $\rho^{-1} \sigma_i \rho = \sigma_{i+1}$, $i \in [d-2]$. It suffices to study $\sigma_1$. The following are all fixed points under $\sigma_1$ in $\SSYT(4,[8])$ with one column:
\begin{align*}
& \scalemath{0.7}{
\begin{ytableau}
1 \\
2 \\
3 \\
7
\end{ytableau}, \ \begin{ytableau}
1 \\
2 \\
4 \\
7
\end{ytableau}, \ \begin{ytableau}
1 \\
2 \\
4 \\
8
\end{ytableau}, \ \begin{ytableau}
1 \\
2 \\
5 \\
6
\end{ytableau}, \ \begin{ytableau}
3 \\
4 \\
7 \\
8
\end{ytableau}, \ \begin{ytableau}
3 \\
5 \\
6 \\
7
\end{ytableau}, \ \begin{ytableau}
3 \\
5 \\
6 \\
8
\end{ytableau}, \ \begin{ytableau}
4 \\
5 \\
6 \\
8
\end{ytableau}. }
\end{align*}
They are all cluster variables. All of them are unstable fixed points. 

The following are all fixed points under $\sigma_1$ in $\SSYT(4,[8])$ with two columns:
\begin{align*} 
& 
\scalemath{0.7}{
\begin{ytableau}
1 & 1 \\
2 & 2 \\
3 & 3 \\
7 & 7
\end{ytableau}, \ \begin{ytableau}
1 & 1 \\
2 & 2 \\
3 & 5 \\
6 & 7
\end{ytableau}, \ \begin{ytableau}
1 & 1 \\
2 & 2 \\
3 & 4 \\
7 & 7
\end{ytableau}, \ \begin{ytableau}
1 & 1 \\
2 & 2 \\
4 & 4 \\
7 & 8
\end{ytableau}, \ \begin{ytableau}
1 & 1 \\
2 & 2 \\
4 & 5 \\
6 & 7
\end{ytableau}, \ \begin{ytableau}
1 & 1 \\
2 & 2 \\
4 & 4 \\
8 & 8
\end{ytableau}, \ \begin{ytableau}
1 & 1 \\
2 & 2 \\
4 & 4 \\
7 & 7
\end{ytableau}, \ \begin{ytableau}
1 & 1 \\
2 & 2 \\
5 & 5 \\
6 & 6
\end{ytableau}, \ \begin{ytableau}
1 & 1 \\
2 & 2 \\
4 & 5 \\
6 & 8
\end{ytableau}, \ \begin{ytableau}
1 & 3 \\
2 & 5 \\
3 & 6 \\
7 & 7
\end{ytableau}, \ \begin{ytableau}
1 & 3 \\
2 & 5 \\
3 & 7 \\
6 & 8
\end{ytableau}, } 
\end{align*}
\begin{align*} 
& 
\scalemath{0.7}{ \begin{ytableau}
1 & 3 \\
2 & 4 \\
3 & 7 \\
7 & 8
\end{ytableau}, \  \begin{ytableau}
1 & 3 \\
2 & 4 \\
4 & 7 \\
8 & 8
\end{ytableau}, \ \begin{ytableau}
1 & 3 \\
2 & 5 \\
4 & 6 \\
7 & 7
\end{ytableau}, \ \begin{ytableau}
1 & 3 \\
2 & 5 \\
4 & 6 \\
7 & 8
\end{ytableau}, \ \begin{ytableau}
1 & 3 \\
2 & 5 \\
4 & 6 \\
8 & 8
\end{ytableau}, \ \begin{ytableau}
1 & 3 \\
2 & 4 \\
4 & 7 \\
7 & 8
\end{ytableau}, \ \begin{ytableau}
1 & 4 \\
2 & 5 \\
3 & 7 \\
6 & 8
\end{ytableau}, \ \begin{ytableau}
1 & 3 \\
2 & 5 \\
5 & 6 \\
6 & 8
\end{ytableau}, \ \begin{ytableau}
1 & 4 \\
2 & 5 \\
4 & 6 \\
8 & 8
\end{ytableau}, 
\ \begin{ytableau}
1 & 4 \\
2 & 5 \\
4 & 7 \\
6 & 8
\end{ytableau}, \ \begin{ytableau}
1 & 3 \\
2 & 5 \\
5 & 6 \\
6 & 7
\end{ytableau}, } 
\end{align*}
\begin{align*} 
& 
\scalemath{0.7}{  \begin{ytableau}
1 & 4 \\
2 & 5 \\
5 & 6 \\
6 & 8
\end{ytableau}, \ \begin{ytableau}
3 & 3 \\
4 & 4 \\
7 & 7 \\
8 & 8
\end{ytableau}, \ \begin{ytableau}
3 & 3 \\
4 & 5 \\
6 & 7 \\
8 & 8
\end{ytableau}, \ \begin{ytableau}
3 & 3 \\
5 & 5 \\
6 & 6 \\
8 & 8
\end{ytableau}, \ \begin{ytableau}
3 & 3 \\
5 & 5 \\
6 & 6 \\
7 & 7
\end{ytableau}, \ \begin{ytableau}
3 & 3 \\
4 & 5 \\
6 & 7 \\
7 & 8
\end{ytableau}, \ \begin{ytableau}
3 & 3 \\
5 & 5 \\
6 & 6 \\
7 & 8
\end{ytableau}, \ \begin{ytableau}
3 & 4 \\
4 & 5 \\
6 & 7 \\
8 & 8
\end{ytableau}, \ \begin{ytableau}
3 & 4 \\
5 & 5 \\
6 & 6 \\
8 & 8
\end{ytableau}, \ \begin{ytableau}
4 & 4 \\
5 & 5 \\
6 & 6 \\
8 & 8
\end{ytableau}, \ \begin{ytableau}
1 & 3 \\
2 & 5 \\
4 & 7 \\
6 & 8
\end{ytableau}.
}
\end{align*}
All the above tableaux except the last one are unstable fixed points. The last one is a stable fixed point and it is prime and non-real. These unstable fixed points correspond to the following cluster monomials respectively:
\begin{align*}
& P_{1237}^2, \ P_{1256}P_{1237}, \ P_{1237}P_{1247}, \ P_{1247}P_{1248}, \ P_{1247}P_{1256}, \ P_{1248}^2, \ P_{1247}^2, \ P_{1256}^2, \ P_{1248}P_{1256}, \\
& P_{1237}P_{3567}, \ P_{1236}P_{3578}-P_{1235}P_{3678}, \ P_{1237}P_{3478}, \ P_{1248}P_{3478}, \ P_{1247}P_{3567}, \ P_{1247}P_{3568}, \\
& P_{1248}P_{3568}, \ P_{1247}P_{3478}, \ P_{1236}P_{4578}-P_{1235}P_{4678} + P_{1234}P_{5678}, \ P_{1256}P_{3568}, \ P_{1248}P_{4568}, \\
& P_{1246}P_{4578} - P_{1245}P_{4678}, \ P_{1256}P_{3567}, \ P_{1256}P_{4568}, \ P_{3478}^2, \ P_{3478}P_{3568}, \ P_{3568}^2, \ P_{3567}^2, \ P_{3567}P_{3478}, \\
& P_{3567}P_{3568}, \ P_{3478}P_{4568}, \ P_{3568}P_{4568}, \ P_{4568}^2.
\end{align*}  

The following are all fixed points under $\sigma_1$ in $\SSYT(4,[8])$ with three columns:
\begin{align*}
& \scalemath{0.7}{
\begin{ytableau}
1 & 1 & 1 \\
2 & 2 & 2 \\
4 & 4 & 4 \\
7 & 7 & 7
\end{ytableau}, \ \begin{ytableau}
1 & 1 & 1 \\
2 & 2 & 2 \\
4 & 4 & 4 \\
7 & 7 & 8
\end{ytableau}, \ \begin{ytableau}
1 & 1 & 1 \\
2 & 2 & 2 \\
4 & 4 & 4 \\
7 & 8 & 8
\end{ytableau}, \ \begin{ytableau}
1 & 1 & 1 \\
2 & 2 & 2 \\
4 & 4 & 4 \\
8 & 8 & 8
\end{ytableau}, \ \begin{ytableau}
1 & 1 & 1 \\
2 & 2 & 2 \\
5 & 5 & 5 \\
6 & 6 & 6
\end{ytableau}, \ \begin{ytableau}
1 & 1 & 1 \\
2 & 2 & 2 \\
3 & 3 & 5 \\
6 & 7 & 7
\end{ytableau}, \ \begin{ytableau}
1 & 1 & 1 \\
2 & 2 & 2 \\
3 & 4 & 5 \\
6 & 7 & 7
\end{ytableau}, \ \begin{ytableau}
1 & 1 & 1 \\
2 & 2 & 2 \\
3 & 4 & 4 \\
7 & 7 & 7
\end{ytableau}, \ \begin{ytableau}
1 & 1 & 1 \\
2 & 2 & 2 \\
4 & 5 & 5 \\
6 & 6 & 7
\end{ytableau}, \ \begin{ytableau}
1 & 1 & 1 \\
2 & 2 & 2 \\
3 & 3 & 4 \\
7 & 7 & 7
\end{ytableau}, } 
\end{align*}
\begin{align*} 
& 
\scalemath{0.7}{ \begin{ytableau}
1 & 1 & 1 \\
2 & 2 & 2 \\
4 & 5 & 5 \\
6 & 6 & 8
\end{ytableau}, \ \begin{ytableau}
1 & 1 & 1 \\
2 & 2 & 2 \\
4 & 4 & 5 \\
6 & 7 & 7
\end{ytableau}, \ \begin{ytableau}
1 & 1 & 1 \\
2 & 2 & 2 \\
3 & 3 & 3 \\
7 & 7 & 7
\end{ytableau}, \ \begin{ytableau}
1 & 1 & 1 \\
2 & 2 & 2 \\
4 & 4 & 5 \\
6 & 7 & 8
\end{ytableau}, \ \begin{ytableau}
1 & 1 & 1 \\
2 & 2 & 2 \\
3 & 5 & 5 \\
6 & 6 & 7
\end{ytableau}, \ \begin{ytableau}
1 & 1 & 3 \\
2 & 2 & 5 \\
3 & 3 & 6 \\
7 & 7 & 7
\end{ytableau}, \ \begin{ytableau}
1 & 1 & 1 \\
2 & 2 & 2 \\
4 & 4 & 5 \\
6 & 8 & 8
\end{ytableau}, \ \begin{ytableau}
1 & 1 & 3 \\
2 & 2 & 5 \\
3 & 4 & 7 \\
6 & 7 & 8
\end{ytableau}, \ \begin{ytableau}
1 & 1 & 3 \\
2 & 2 & 4 \\
3 & 3 & 7 \\
7 & 7 & 8
\end{ytableau}, \ \begin{ytableau}
1 & 1 & 3 \\
2 & 2 & 5 \\
3 & 5 & 6 \\
6 & 7 & 7
\end{ytableau}, } 
\end{align*}
\begin{align*} 
& 
\scalemath{0.7}{ \begin{ytableau}
1 & 1 & 3 \\
2 & 2 & 5 \\
3 & 3 & 7 \\
6 & 7 & 8
\end{ytableau}, \ \begin{ytableau}
1 & 1 & 4 \\
2 & 2 & 5 \\
3 & 3 & 7 \\
6 & 7 & 8
\end{ytableau}, \ \begin{ytableau}
1 & 1 & 3 \\
2 & 2 & 5 \\
3 & 5 & 7 \\
6 & 6 & 8
\end{ytableau}, \ \begin{ytableau}
1 & 1 & 4 \\
2 & 2 & 5 \\
4 & 5 & 6 \\
6 & 8 & 8
\end{ytableau}, \ \begin{ytableau}
1 & 1 & 3 \\
2 & 2 & 5 \\
4 & 4 & 6 \\
7 & 7 & 7
\end{ytableau}, \ \begin{ytableau}
1 & 1 & 3 \\
2 & 2 & 4 \\
4 & 4 & 7 \\
7 & 7 & 8
\end{ytableau}, \ \begin{ytableau}
1 & 1 & 3 \\
2 & 2 & 5 \\
4 & 5 & 6 \\
6 & 7 & 8
\end{ytableau}, \ \begin{ytableau}
1 & 1 & 3 \\
2 & 2 & 5 \\
5 & 5 & 6 \\
6 & 6 & 7
\end{ytableau}, \ \begin{ytableau}
1 & 1 & 3 \\
2 & 2 & 5 \\
3 & 4 & 6 \\
7 & 7 & 7
\end{ytableau}, \ \begin{ytableau}
1 & 1 & 3 \\
2 & 2 & 4 \\
3 & 4 & 7 \\
7 & 7 & 8
\end{ytableau}, } 
\end{align*}
\begin{align*} 
& 
\scalemath{0.7}{ \begin{ytableau}
1 & 1 & 3 \\
2 & 2 & 5 \\
4 & 5 & 7 \\
6 & 6 & 8
\end{ytableau}, \ \begin{ytableau}
1 & 1 & 3 \\
2 & 2 & 5 \\
4 & 4 & 6 \\
7 & 8 & 8
\end{ytableau}, \ \begin{ytableau}
1 & 1 & 4 \\
2 & 2 & 5 \\
5 & 5 & 6 \\
6 & 6 & 8
\end{ytableau}, \ \begin{ytableau}
1 & 1 & 3 \\
2 & 2 & 5 \\
5 & 5 & 6 \\
6 & 6 & 8
\end{ytableau}, \ \begin{ytableau}
1 & 1 & 3 \\
2 & 2 & 5 \\
4 & 4 & 7 \\
6 & 8 & 8
\end{ytableau}, \ \begin{ytableau}
1 & 1 & 3 \\
2 & 2 & 4 \\
4 & 4 & 7 \\
7 & 8 & 8
\end{ytableau}, \ \begin{ytableau}
1 & 1 & 3 \\
2 & 2 & 4 \\
4 & 4 & 7 \\
8 & 8 & 8
\end{ytableau}, \ \begin{ytableau}
1 & 1 & 4 \\
2 & 2 & 5 \\
3 & 5 & 7 \\
6 & 6 & 8
\end{ytableau}, \ \begin{ytableau}
1 & 1 & 4 \\
2 & 2 & 5 \\
4 & 4 & 7 \\
6 & 7 & 8
\end{ytableau}, \ \begin{ytableau}
1 & 1 & 4 \\
2 & 2 & 5 \\
4 & 4 & 7 \\
6 & 8 & 8
\end{ytableau}, } 
\end{align*}
\begin{align*} 
& 
\scalemath{0.7}{
\begin{ytableau}
1 & 1 & 4 \\
2 & 2 & 5 \\
4 & 4 & 6 \\
8 & 8 & 8
\end{ytableau}, \ \begin{ytableau}
1 & 1 & 3 \\
2 & 2 & 5 \\
4 & 4 & 6 \\
7 & 7 & 8
\end{ytableau}, \ \begin{ytableau}
1 & 1 & 3 \\
2 & 2 & 5 \\
4 & 4 & 6 \\
8 & 8 & 8
\end{ytableau}, \ \begin{ytableau}
1 & 1 & 3 \\
2 & 2 & 5 \\
4 & 4 & 7 \\
6 & 7 & 8
\end{ytableau}, \ \begin{ytableau}
1 & 1 & 3 \\
2 & 2 & 5 \\
4 & 5 & 6 \\
6 & 8 & 8
\end{ytableau}, \ \begin{ytableau}
1 & 1 & 3 \\
2 & 2 & 5 \\
4 & 5 & 6 \\
6 & 7 & 7
\end{ytableau}, \ \begin{ytableau}
1 & 1 & 4 \\
2 & 2 & 5 \\
3 & 4 & 7 \\
6 & 7 & 8
\end{ytableau}, \ \begin{ytableau}
1 & 1 & 4 \\
2 & 2 & 5 \\
4 & 5 & 7 \\
6 & 6 & 8
\end{ytableau}, \ \begin{ytableau}
1 & 3 & 3 \\
2 & 4 & 4 \\
3 & 7 & 7 \\
7 & 8 & 8
\end{ytableau}, \ \begin{ytableau}
1 & 3 & 3 \\
2 & 5 & 5 \\
3 & 6 & 7 \\
6 & 8 & 8
\end{ytableau}, } 
\end{align*}
\begin{align*} 
& 
\scalemath{0.7}{
\begin{ytableau}
1 & 3 & 3 \\
2 & 4 & 4 \\
4 & 7 & 7 \\
8 & 8 & 8
\end{ytableau}, \ \begin{ytableau}
1 & 3 & 3 \\
2 & 4 & 5 \\
4 & 6 & 7 \\
7 & 8 & 8
\end{ytableau}, \ \begin{ytableau}
1 & 3 & 3 \\
2 & 5 & 5 \\
4 & 6 & 6 \\
7 & 7 & 8
\end{ytableau}, \ \begin{ytableau}
1 & 3 & 3 \\
2 & 4 & 5 \\
4 & 7 & 7 \\
6 & 8 & 8
\end{ytableau}, \ \begin{ytableau}
1 & 3 & 3 \\
2 & 5 & 5 \\
4 & 6 & 7 \\
6 & 7 & 8
\end{ytableau}, \ \begin{ytableau}
1 & 3 & 3 \\
2 & 5 & 5 \\
4 & 6 & 6 \\
8 & 8 & 8
\end{ytableau}, \ \begin{ytableau}
1 & 3 & 4 \\
2 & 4 & 5 \\
3 & 7 & 7 \\
6 & 8 & 8
\end{ytableau}, \ \begin{ytableau}
1 & 3 & 4 \\
2 & 4 & 5 \\
4 & 6 & 7 \\
8 & 8 & 8
\end{ytableau}, \ \begin{ytableau}
1 & 3 & 3 \\
2 & 4 & 5 \\
3 & 6 & 7 \\
7 & 7 & 8
\end{ytableau}, \ \begin{ytableau}
1 & 4 & 4 \\
2 & 5 & 5 \\
4 & 6 & 7 \\
6 & 8 & 8
\end{ytableau}, } 
\end{align*}
\begin{align*} 
& 
\scalemath{0.7}{ 
\begin{ytableau}
1 & 3 & 4 \\
2 & 5 & 5 \\
4 & 6 & 7 \\
6 & 8 & 8
\end{ytableau}, \ \begin{ytableau}
1 & 4 & 4 \\
2 & 5 & 5 \\
4 & 6 & 6 \\
8 & 8 & 8
\end{ytableau}, \ \begin{ytableau}
1 & 3 & 3 \\
2 & 4 & 5 \\
4 & 6 & 7 \\
8 & 8 & 8
\end{ytableau}, \ \begin{ytableau}
1 & 3 & 3 \\
2 & 4 & 5 \\
3 & 7 & 7 \\
6 & 8 & 8
\end{ytableau}, \ \begin{ytableau}
1 & 4 & 4 \\
2 & 5 & 5 \\
3 & 6 & 7 \\
6 & 8 & 8
\end{ytableau}, \ \begin{ytableau}
1 & 3 & 3 \\
2 & 4 & 5 \\
4 & 6 & 7 \\
7 & 7 & 8
\end{ytableau}, \ \begin{ytableau}
1 & 3 & 3 \\
2 & 5 & 5 \\
4 & 6 & 6 \\
7 & 7 & 7
\end{ytableau}, \ \begin{ytableau}
1 & 3 & 4 \\
2 & 5 & 5 \\
5 & 6 & 6 \\
6 & 8 & 8
\end{ytableau}, \ \begin{ytableau}
1 & 3 & 4 \\
2 & 5 & 5 \\
4 & 6 & 6 \\
8 & 8 & 8
\end{ytableau}, \ \begin{ytableau}
1 & 3 & 3 \\
2 & 5 & 5 \\
4 & 6 & 6 \\
7 & 8 & 8
\end{ytableau}, } 
\end{align*}
\begin{align*} 
& 
\scalemath{0.7}{
\begin{ytableau}
1 & 3 & 3 \\
2 & 5 & 5 \\
5 & 6 & 6 \\
6 & 7 & 7
\end{ytableau}, \ \begin{ytableau}
1 & 3 & 3 \\
2 & 5 & 5 \\
4 & 6 & 7 \\
6 & 8 & 8
\end{ytableau}, \ \begin{ytableau}
1 & 3 & 3 \\
2 & 5 & 5 \\
3 & 6 & 6 \\
7 & 7 & 7
\end{ytableau}, \ \begin{ytableau}
1 & 3 & 3 \\
2 & 5 & 5 \\
5 & 6 & 6 \\
6 & 7 & 8
\end{ytableau}, \ \begin{ytableau}
1 & 3 & 4 \\
2 & 5 & 5 \\
3 & 6 & 7 \\
6 & 8 & 8
\end{ytableau}, \ \begin{ytableau}
1 & 3 & 3 \\
2 & 5 & 5 \\
5 & 6 & 6 \\
6 & 8 & 8
\end{ytableau}, \ \begin{ytableau}
1 & 3 & 3 \\
2 & 4 & 4 \\
4 & 7 & 7 \\
7 & 8 & 8
\end{ytableau}, \ \begin{ytableau}
1 & 3 & 4 \\
2 & 4 & 5 \\
4 & 7 & 7 \\
6 & 8 & 8
\end{ytableau}, \ \begin{ytableau}
1 & 3 & 3 \\
2 & 5 & 5 \\
3 & 6 & 7 \\
6 & 7 & 8
\end{ytableau}, \ \begin{ytableau}
3 & 3 & 3 \\
4 & 4 & 4 \\
7 & 7 & 7 \\
8 & 8 & 8
\end{ytableau}, } 
\end{align*}
\begin{align*} 
& 
\scalemath{0.7}{
\begin{ytableau}
1 & 4 & 4 \\
2 & 5 & 5 \\
5 & 6 & 6 \\
6 & 8 & 8
\end{ytableau}, \ \begin{ytableau}
4 & 4 & 4 \\
5 & 5 & 5 \\
6 & 6 & 6 \\
8 & 8 & 8
\end{ytableau}, \ \begin{ytableau}
3 & 3 & 3 \\
5 & 5 & 5 \\
6 & 6 & 6 \\
7 & 7 & 7
\end{ytableau}, \ \begin{ytableau}
3 & 3 & 3 \\
5 & 5 & 5 \\
6 & 6 & 6 \\
7 & 7 & 8
\end{ytableau}, \ \begin{ytableau}
3 & 3 & 3 \\
4 & 4 & 5 \\
6 & 7 & 7 \\
7 & 8 & 8
\end{ytableau}, \ \begin{ytableau}
3 & 3 & 3 \\
5 & 5 & 5 \\
6 & 6 & 6 \\
7 & 8 & 8
\end{ytableau}, \ \begin{ytableau}
3 & 3 & 4 \\
4 & 4 & 5 \\
6 & 7 & 7 \\
8 & 8 & 8
\end{ytableau}, \ \begin{ytableau}
3 & 3 & 3 \\
4 & 4 & 5 \\
6 & 7 & 7 \\
8 & 8 & 8
\end{ytableau}, \ \begin{ytableau}
3 & 3 & 3 \\
5 & 5 & 5 \\
6 & 6 & 6 \\
8 & 8 & 8
\end{ytableau}, \ \begin{ytableau}
3 & 4 & 4 \\
5 & 5 & 5 \\
6 & 6 & 6 \\
8 & 8 & 8
\end{ytableau}, } 
\end{align*}
\begin{align*} 
& 
\scalemath{0.7}{ \begin{ytableau}
3 & 3 & 4 \\
5 & 5 & 5 \\
6 & 6 & 6 \\
8 & 8 & 8
\end{ytableau}, \ \begin{ytableau}
3 & 3 & 3 \\
4 & 5 & 5 \\
6 & 6 & 7 \\
7 & 7 & 8
\end{ytableau}, \ \begin{ytableau}
3 & 4 & 4 \\
4 & 5 & 5 \\
6 & 6 & 7 \\
8 & 8 & 8
\end{ytableau}, \ \begin{ytableau}
3 & 3 & 3 \\
4 & 5 & 5 \\
6 & 6 & 7 \\
7 & 8 & 8
\end{ytableau}, \ \begin{ytableau}
3 & 3 & 3 \\
4 & 5 & 5 \\
6 & 6 & 7 \\
8 & 8 & 8
\end{ytableau}, \ \begin{ytableau}
3 & 3 & 4 \\
4 & 5 & 5 \\
6 & 6 & 7 \\
8 & 8 & 8
\end{ytableau}. }
\end{align*}
All of these fixed points are unstable. These tableaux also correspond to the cluster monomials. 

\begin{remark}
An unstable fixed point can also correspond to a non-real element. For example, $\vcenter{\hbox{\scalebox{0.6}{ \begin{ytableau}
    1 & 2  \\
3 & 4  \\
5 & 6 \\
7 & 8
\end{ytableau} }}}$ is non-real for $\Gr(4,12)$. It is a stable fixed point for $\Gr(4,8)$. But it is an unstable fixed point under $\sigma_2$ for $\Gr(4,12)$. 
\end{remark}

\section{\texorpdfstring{Stable fixed points in $\CC[\Gr(k,n)]$}{Stable fixed points in C[Gr(k,n)]}}
\label{sec:stable fixed points}
In this section, we study stable fixed points of $\CC[\Gr(k,n)]$.

\subsection{\texorpdfstring{Stable fixed points in $\CC[\Gr(k,n)]$}{Stable fixed points in C[Gr(k,n)]}}

Consider $\CC[\Gr(3,9)]$. Denote by $\sigma_3 = \rho \circ \sigma_2 \circ \rho^{-1}$. The stable fixed points of $\sigma_1$, $\sigma_2$, $\sigma_3$ in $\CC[\Gr(3,9)]$ are:
\begin{align} \label{eq:stable fixed points Gr39}
\vcenter{\hbox{\scalebox{0.6}{
\begin{ytableau}
1 & 3 & 4 \\
2 & 6 & 7 \\
5 & 8 & 9
\end{ytableau}}}},  \ \vcenter{\hbox{\scalebox{0.6}{ \begin{ytableau}
    1 & 2 & 5 \\
3 & 4 & 8 \\
6 & 7 & 9
\end{ytableau} }}}, \ \vcenter{\hbox{\scalebox{0.6}{ \begin{ytableau}
    1 & 2 & 3 \\
4 & 5 & 6 \\
7 & 8 & 9
\end{ytableau}}}},
\end{align}
respectively. These tableaux are prime and non-real and they appeared in Example 8.1 in \cite{CDFL}.

Consider $\CC[\Gr(4,8)]$. Denote by $\sigma_4 = \rho \circ \sigma_3 \circ \rho^{-1}$. The stable fixed points of $\sigma_1$, $\sigma_2$, $\sigma_3$, $\sigma_4$ in $\CC[\Gr(4,8)]$ are:
\begin{align} \label{eq:stable fixed points Gr48}
\vcenter{\hbox{\scalebox{0.6}{
\begin{ytableau}
1 & 3 \\
2 & 5 \\
4 & 7 \\
6 & 8
\end{ytableau}}}}, \ \vcenter{\hbox{\scalebox{0.6}{
\begin{ytableau}
1 & 2 \\
3 & 4 \\
5 & 6 \\
7 & 8
\end{ytableau}}}}, \ \vcenter{\hbox{\scalebox{0.6}{ \begin{ytableau}
1 & 3 \\
2 & 5 \\
4 & 7 \\
6 & 8
\end{ytableau}}}}, \ \vcenter{\hbox{\scalebox{0.6}{ \begin{ytableau}
1 & 2 \\
3 & 4 \\
5 & 6 \\
7 & 8
\end{ytableau}}}},
\end{align}
respectively. The first and the third tableaux are the same. The second and the fourth tableaux are the same. These tableaux are prime and non-real and they appeared in Example 8.1 in \cite{CDFL}. 

Consider $\CC[\Gr(4,12)]$. Denote by $\sigma_4 = \rho \circ \sigma_3 \circ \rho^{-1}$. The stable fixed points of $\sigma_1$, $\sigma_2$, $\sigma_3$, $\sigma_4$ in $\CC[\Gr(4,12)]$ are:
\begin{align*}
\vcenter{\hbox{\scalebox{0.6}{
\begin{ytableau}
1 & 3 & 5 \\
2 & 7 & 9 \\
4 & 8 & 11 \\
6 & 10 & 12
\end{ytableau}}}}, \ \vcenter{\hbox{\scalebox{0.6}{ \begin{ytableau}
1 & 2 & 6 \\
3 & 4 & 8 \\
5 & 9 & 10 \\
7 & 11 & 12
\end{ytableau}}}}, \ \vcenter{\hbox{\scalebox{0.6}{ \begin{ytableau}
1 & 3 & 7 \\
2 & 5 & 9 \\
4 & 6 & 11 \\
8 & 10 & 12
\end{ytableau}}}}, \ \vcenter{\hbox{\scalebox{0.6}{ \begin{ytableau}
1 & 2 & 4 \\
3 & 6 & 8 \\
5 & 7 & 10 \\
9 & 11 & 12
\end{ytableau}}}},
\end{align*}
respectively. These tableaux are easy to check to be prime by using Theorem 1.1 in \cite{BL22}. We expect that they are non-real. We can check the non-reality by checking if $\ch(T \cup T) \ne \ch(T)^2$ holds. For these examples, $\ch(T \cup T)$ is difficult to compute.

Consider $\CC[\Gr(5,10)]$. Denote by $\sigma_5 = \rho \circ \sigma_4 \circ \rho^{-1}$. The stable fixed points of $\sigma_1, \ldots, \sigma_5$ in $\CC[\Gr(5,10)]$ are:
\begin{align*}
\vcenter{\hbox{\scalebox{0.6}{
\begin{ytableau}
1 & 3 \\
2 & 6 \\
4 & 8 \\
5 & 9 \\
7 & 10
\end{ytableau}}}}, \ \vcenter{\hbox{\scalebox{0.6}{ \begin{ytableau}
1 & 2 \\
3 & 4 \\
5 & 7 \\
6 & 9 \\
8 & 10
\end{ytableau}}}}, \ \vcenter{\hbox{\scalebox{0.6}{ \begin{ytableau}
1 & 3 \\
2 & 5 \\
4 & 6 \\
7 & 8 \\
9 & 10
\end{ytableau}}}}, \ \vcenter{\hbox{\scalebox{0.6}{ \begin{ytableau}
1 & 4 \\
2 & 6 \\
3 & 7 \\
5 & 9 \\
8 & 10
\end{ytableau}}}}, \ \vcenter{\hbox{\scalebox{0.6}{ \begin{ytableau}
1 & 2 \\
3 & 5 \\
4 & 7 \\
6 & 8 \\
9 & 10
\end{ytableau}}}},
\end{align*}
respectively. These tableaux are prime by Theorem 1.1 in \cite{BL22}. We expect that these tableaux are non-real.

Consider $\CC[\Gr(4,10)]$. Denote by $\sigma_2 = \rho \circ \sigma_1 \circ \rho^{-1}$. The stable fixed points of $\sigma_1$ and $\sigma_2$ are:
\begin{align*}
\vcenter{\hbox{\scalebox{0.6}{
\begin{ytableau}
1 & 1 & 3 & 3 & 5 \\
2 & 2 & 5 & 7 & 7 \\
4 & 4 & 6 & 9 & 9 \\
6 & 8 & 8 & 10 & 10
\end{ytableau}}}}, \ \vcenter{\hbox{\scalebox{0.6}{ \begin{ytableau}
1 & 1 & 2 & 2 & 4 \\
3 & 3 & 4 & 6 & 6 \\
5 & 5 & 7 & 8 & 8 \\
7 & 9 & 9 & 10 & 10
\end{ytableau}}}},
\end{align*}
respectively. We expect that these tableaux are prime and non-real.

\subsection{Tableaux obtained by applying the braid group action to stable fixed points}

Recall that the Euler totient function $\varphi(n)$ counts the positive integers up to a given integer $n$ that are relatively prime to $n$.
Surprisingly, the numbers of tableaux of rank $r$ (we call the number of columns of a tableau its rank) in $\SSYT(3,[9])$ and $\SSYT(4,[8])$ obtained by braid group actions on stable fixed points are determined by Euler's totient function. The sequence of numbers of rank $r$ ($r \ge 1$) tableaux in $\SSYT(3,[9])$ obtained by the braid group action on the stable fixed points is 
\begin{align*}
& 0,0,3,0,0,3,0,0,6,0,0,6,0,0,12,0,0,6,0,0,18,0,0,12,0,0,12,0,0,12, \\
& 0,0,18,0,0,12,0,0,30,0,0,12,0,0,36,0,0,18,0,0,24,0,0,24, \ldots
\end{align*}
The sequence of numbers of rank $r$ ($r \ge 1$) tableaux in $\SSYT(4,[8])$ obtained by the braid group action on the stable fixed points is 
\begin{align*}
& 0,2,0,2,0,4,0,4,0,8,0,4,0,12,0,8,0,12,0,8,0,12,0,8,0, \\
& 20,0,8,0,24,0,12,0,16,0,16,0,32,0,12,0,36,0,16,0,24,0,20,0,44, \ldots
\end{align*} 

\begin{theorem} \label{thm:prime non real modules Gr39 and Gr48}
The number of rank $r$ ($r \ge 1$) tableaux in $\SSYT(3,[9])$ obtained by applying the braid group ${\rm Br}_{3}$ action to the stable fixed points is 
\begin{align*}
N_{3,9,r} = 
\begin{cases}
0, & r \pmod 3 = i, \ i \in \{1,2\}, \\
3 \phi(r/3), & r \pmod 3 = 0,
\end{cases}
\end{align*}
where $\phi(m)$ is the Euler totient function.

The number of rank $r$ ($r \ge 1$) tableaux in $\SSYT(4,[8])$ obtained by applying the braid group $\Br_4$ action to the stable fixed points is 
\begin{align*}
N_{4,8,r} = 
\begin{cases}
0, & r \pmod 2 = 1, \\
2 \phi(r/2), & r \pmod 2 = 0.
\end{cases}
\end{align*}
\end{theorem}

We will explain Theorem \ref{thm:prime non real modules Gr39 and Gr48} in the next two subsections. 

\begin{remark}
We expect that the braid group action on the sub-category $\mathcal{C}_{\ell}$ of the category of finite dimensional representations of $U_q(\widehat{\mathfrak{sl}_k})$ introduced in \cite{KKOP21} coincides with the braid group action on $\mathcal{C}_{\ell}$ induced by the braid group action on $\CC[\Gr(k,n)]$. By Proposition 6.2 in \cite{FHOO23}, Theorem 3.1 in \cite{KKOP21}, and Theorems 6.1.1, 6.2.1 in \cite{KQW22}, braid group action sends simple modules to simple modules. By Theorem 5.3 in \cite{Fra20}, any $\sigma \in {\rm Br}_d$ is an algebra homomorphism on $\CC[\Gr(k,n)]$, $d = \gcd(k,n)$. Therefore we expect that after removing frozen factors, the braid group action sends a prime (resp. non-prime) module to a prime (resp. non-prime) module, and sends a real (resp. non-real) module to a real (resp. non-real) module. 
\end{remark}

It is conjectured in \cite[Conjecture 11.5]{EL23} every prime non-real element in the dual canonical basis of $\CC[\Gr(3,9)]$ (resp. $\CC[\Gr(4,8)]$) is of the form $\ch(T^{(1)} \cup \cdots \cup T^{(r)})$ for some $r \in \ZZ_{\ge 1}$ and some $T^{(j)}$'s corresponding to limit g-vectors. We give the following conjecture.
\begin{conjecture} \label{conj: all imaginary prime tableaux in Gr39 and Gr48}
All limit g-vectors for $\CC[\Gr(3,9)]$ and $\CC[\Gr(4,8)]$ can be obtained by the braid group action and every prime non-real element in the dual canonical basis of $\CC[\Gr(3,9)]$ (resp. $\CC[\Gr(4,8)]$) is of the form $\ch(T^{(1)} \cup \cdots \cup T^{(r)})$ for some $r \in \ZZ_{\ge 1}$ and some $T^{(j)}$'s which are obtained from the stable fixed points in (\ref{eq:stable fixed points Gr39}) (resp. (\ref{eq:stable fixed points Gr48})) by the braid group action. 
\end{conjecture}

\begin{remark}
Consider $\CC[\Gr(4,8)]$, let 
$T = \vcenter{\hbox{\scalebox{0.6}{ \begin{ytableau}
    1 & 3 \\ 2 & 5 \\ 4 & 7 \\ 6 & 8
\end{ytableau}}}}$ and $A = \ch(T)$, $B= \ch(T)^2 - \ch(T\cup T)$. As explained in Remark 11.6 in \cite{EL23}, although for some $r$, $\ch(T^{\cup r})$ could be a non-prime polynomial in $A, B$, $\ch(T^{\cup r})$ could still be prime as an element in the dual canonical basis of $\CC[\Gr(k,n)]$ (meaning that $\ch(T^{\cup r}) \ne \ch(T')\ch(T'')$ for any non-trivial tableaux $T', T'' \in \SSYT(k,[n])$). 
\end{remark}

\subsection{\texorpdfstring{Proof of Theorem \ref{thm:prime non real modules Gr39 and Gr48} for $\Gr(4,8)$}{Proof of Theorem \ref{thm:prime non real modules Gr39 and Gr48} for Gr(4,8)}} \label{subsec:prime nonreal tableaux for Gr48}

Acting on the $\hat{y}$-variables of the initial seed with $\sigma_1^{-1}$ and $\sigma_1$ we find
\begin{align}
   \renewcommand{\arraystretch}{1.2}
\begin{bmatrix}
\chi_{11}\\ \chi_{12} \\ \chi_{13} \\ \chi_{21} \\ \chi_{22} \\ \chi_{23} \\ \chi_{31} \\ \chi_{32} \\ \chi_{33}
    \end{bmatrix}
\xmapsto{\sigma_1^{-1}}
  \begingroup
    \renewcommand{\arraystretch}{1.8}
\begin{bmatrix}
 \frac{1}{\chi_{12} (1 + \chi_{13})}\\
 \frac{\chi_{23} (1 + \chi_{33})(1 + \chi_{12} + \chi_{12} \chi_{13})}{
  1 + \chi_{13} + \chi_{13} \chi_{23} + \chi_{13} \chi_{23} \chi_{33}}\\
 \frac{1 + \chi_{13}}{\chi_{13} \chi_{23} (1 + \chi_{33})}\\
 {\scriptstyle \chi_{11} (1 + \chi_{21}) (1 + \chi_{12} + \chi_{12} \chi_{13})}\\
 \frac{\chi_{12} \chi_{21} \chi_{22} \chi_{33}(1 + \chi_{13})}{(1 + \chi_{21}) (1 + \chi_{33})(1 + \chi_{12} + \chi_{12} \chi_{13})}\\
 \frac{1 + \chi_{13} + \chi_{13} \chi_{23} + \chi_{13} \chi_{23} \chi_{33}}{\chi_{33}(1 + \chi_{13})}\\
 \frac{\chi_{21} \chi_{31}}{1 + \chi_{21}}\\
 \frac{1}{\chi_{21}}\\
 {\scriptstyle \chi_{32} (1 + \chi_{21}) (1 + \chi_{33})}
\end{bmatrix}
,
\endgroup
\quad
   \renewcommand{\arraystretch}{1.2}
\begin{bmatrix}
\chi_{11}\\ \chi_{12} \\ \chi_{13} \\ \chi_{21} \\ \chi_{22} \\ \chi_{23} \\ \chi_{31} \\ \chi_{32} \\ \chi_{33}
    \end{bmatrix}
\xmapsto{\sigma_1}
  \begingroup
    \renewcommand{\arraystretch}{1.8}
\begin{bmatrix}
 \frac{\chi_{11} \chi_{21} \chi_{32}}{(1 + \chi_{11}) (1 + \chi_{32})}\\
 \frac{\chi_{12} (1 + \chi_{13})}{1 + \chi_{11} + \chi_{11} \chi_{12} + \chi_{11} \chi_{12} \chi_{13}}\\
 \frac{1 + \chi_{11}}{\chi_{11} \chi_{12} (1 + \chi_{13})}\\
 \frac{1}{\chi_{32}}\\
 \frac{\chi_{22}(1 + \chi_{11})(1 + \chi_{32})(1 + \chi_{13} + \chi_{13} \chi_{23})}{1 + \chi_{13}}\\
 \frac{\chi_{23}(1 + \chi_{11} + \chi_{11} \chi_{12} + \chi_{11} \chi_{12} \chi_{13})}{(1 + \chi_{11}) (1 + \chi_{13} + 
     \chi_{13} \chi_{23})}\\
 {\scriptstyle \chi_{31} (1 + \chi_{32})}\\
 \frac{\chi_{13} \chi_{23} \chi_{32} \chi_{33}}{(1 + \chi_{32})(1 + \chi_{13} + \chi_{13} \chi_{23})}\\
 \frac{1 + \chi_{13}}{\chi_{13} \chi_{23}}
\end{bmatrix}
.
\endgroup
\notag
\end{align}
Tropicalising the first map above with `max' convention gives the piecewise linear map $Q^+_{\sigma_1}$ which acts on g-vectors as $Q^+_{\sigma_1} {\bf g}(x) = {\bf g} (\sigma_1(x))$. Tropicalising the second map above with `max' convention gives us the map $Q^+_{{\sigma_1}^{-1}} = (Q^+_{\sigma_1})^{-1}$ which acts on g-vectors as $Q^+_{\sigma_1^{-1}} {\bf g}(x) = {\bf g}(\sigma_1^{-1}(x))$.

Acting on the $\hat{y}$-variables of the initial seed with $\sigma_2^{-1}$ and $\sigma_2$ we find
\begin{align}
   \renewcommand{\arraystretch}{1.2}
\begin{bmatrix}
\chi_{11}\\ \chi_{12} \\ \chi_{13} \\ \chi_{21} \\ \chi_{22} \\ \chi_{23} \\ \chi_{31} \\ \chi_{32} \\ \chi_{33}
    \end{bmatrix}
\xmapsto{\sigma_2^{-1}}
  \begingroup
    \renewcommand{\arraystretch}{1.8}
\begin{bmatrix}
 \frac{\chi_{11} \chi_{22} (1 + \chi_{22} + \chi_{22} \chi_{23} + \chi_{22} \chi_{32} + \chi_{22} \chi_{23} \chi_{32})}{(1 + \chi_{22} + 
     \chi_{22} \chi_{23}) (1 + \chi_{22} + \chi_{22} \chi_{32})}\\
 \frac{\chi_{23}}{1 + \chi_{22} + \chi_{22} \chi_{23} + \chi_{22} \chi_{32} + \chi_{22} \chi_{23} \chi_{32}}\\
 {\scriptstyle \chi_{12} (1 + \chi_{22} + \chi_{22} \chi_{32})}\\
 \frac{\chi_{32}}{1 + \chi_{22} + \chi_{22} \chi_{23} + \chi_{22} \chi_{32} + \chi_{22} \chi_{23} \chi_{32}}\\
 \frac{(1 + \chi_{22} + \chi_{22} \chi_{23}) (1 + \chi_{22} + \chi_{22} \chi_{32})}{\chi_{22} \chi_{23} \chi_{32}}\\
 \frac{\chi_{13} (1 + \chi_{22} + \chi_{22} \chi_{23} + \chi_{22} \chi_{32} + \chi_{22} \chi_{23} \chi_{32})}{
  1 + \chi_{22} + \chi_{22} \chi_{32}}\\
 {\scriptstyle \chi_{21} (1 + \chi_{22} + \chi_{22} \chi_{23})}\\
 \frac{\chi_{31} (1 + \chi_{22} + \chi_{22} \chi_{23} + \chi_{22} \chi_{32} + \chi_{22} \chi_{23} \chi_{32})}{
  1 + \chi_{22} + \chi_{22} \chi_{23}}\\
 \frac{\chi_{22} \chi_{23} \chi_{32} \chi_{33}}{1 + \chi_{22} + \chi_{22} \chi_{23} + \chi_{22} \chi_{32} + \chi_{22} \chi_{23} \chi_{32}}
\end{bmatrix}
,
\endgroup
\quad
   \renewcommand{\arraystretch}{1.2}
\begin{bmatrix}
\chi_{11}\\ \chi_{12} \\ \chi_{13} \\ \chi_{21} \\ \chi_{22} \\ \chi_{23} \\ \chi_{31} \\ \chi_{32} \\ \chi_{33}
    \end{bmatrix}
\xmapsto{\sigma_2}
  \begingroup
    \renewcommand{\arraystretch}{1.8}
\begin{bmatrix}
 {\scriptstyle \chi_{11} (1 + \chi_{12} + \chi_{21} + \chi_{12} \chi_{21} + \chi_{12} \chi_{21} \chi_{22})}\\
 \frac{\chi_{12} \chi_{13} (1 + \chi_{21} + \chi_{21} \chi_{22})}{1 + \chi_{12} + \chi_{21} + \chi_{12} \chi_{21} + \chi_{12} \chi_{21} \chi_{22}}\\
 \frac{\chi_{21} \chi_{22} \chi_{23}}{1 + \chi_{21} + \chi_{21} \chi_{22}}\\
 \frac{\chi_{21} \chi_{31}(1 + \chi_{12} + \chi_{12} \chi_{22})}{1 + \chi_{12} + \chi_{21} + \chi_{12} \chi_{21} + \chi_{12} \chi_{21} \chi_{22}}\\
 \frac{\chi_{22}}{(1 + \chi_{12} + \chi_{12} \chi_{22}) (1 + \chi_{21} + \chi_{21} \chi_{22})}\\
 \frac{1 + \chi_{12} + \chi_{21} + \chi_{12} \chi_{21} + \chi_{12} \chi_{21} \chi_{22}}{\chi_{21} \chi_{22}}\\
 \frac{\chi_{12} \chi_{22} \chi_{32}}{1 + \chi_{12} + \chi_{12} \chi_{22}}\\
 \frac{1 + \chi_{12} + \chi_{21} + \chi_{12} \chi_{21} + \chi_{12} \chi_{21} \chi_{22}}{\chi_{12} \chi_{22}}\\
 \frac{\chi_{33}(1 + \chi_{12} + \chi_{12} \chi_{22}) (1 + \chi_{21} + \chi_{21} \chi_{22}) }{1 + \chi_{12} + \chi_{21} + \chi_{12} \chi_{21} + \chi_{12} \chi_{21} \chi_{22}}
\end{bmatrix}
.
\endgroup
\notag
\end{align}
Tropicalising the first map above with `max' convention gives the piecewise linear map $Q^+_{\sigma_2}$ which acts on g-vectors as $Q^+_{\sigma_2} {\bf g}(x) = {\bf g} (\sigma_2(x))$. Tropicalising the second map above with `max' convention gives us the map $Q^+_{{\sigma_2}^{-1}} = (Q^+_{\sigma_2})^{-1}$ which acts on g-vectors as $Q^+_{\sigma_2^{-1}} {\bf g}(x) = {\bf g}(\sigma_2^{-1}(x))$.

Acting on the $\hat{y}$-variables of the initial seed with $\sigma_3^{-1}$ and $\sigma_3$ we find
\begin{align}
   \renewcommand{\arraystretch}{1.2}
\begin{bmatrix}
\chi_{11}\\ \chi_{12} \\ \chi_{13} \\ \chi_{21} \\ \chi_{22} \\ \chi_{23} \\ \chi_{31} \\ \chi_{32} \\ \chi_{33}
    \end{bmatrix}
\xmapsto{\sigma_3^{-1}}
  \begingroup
    \renewcommand{\arraystretch}{1.8}
\begin{bmatrix}
 \frac{1}{\chi_{21} (1 + \chi_{31})}\\
 {\scriptstyle \chi_{11} (1 + \chi_{12}) (1 + \chi_{21} + \chi_{21} \chi_{31})}\\
 \frac{\chi_{12} \chi_{13}}{1 + \chi_{12}}\\
 \frac{\chi_{32} (1 + \chi_{33})(1 + \chi_{21} + \chi_{21} \chi_{31})}{1 + \chi_{31} + \chi_{31} \chi_{32} + \chi_{31} \chi_{32} \chi_{33}}\\
 \frac{\chi_{12} \chi_{21} \chi_{22}\chi_{33} (1 + \chi_{31})}{(1 + \chi_{12}) (1 + 
     \chi_{33}) (1 + \chi_{21} + \chi_{21} \chi_{31})}\\
 \frac{1}{\chi_{12}}\\
 \frac{1 + \chi_{31}}{\chi_{31} \chi_{32} (1 + \chi_{33})}\\
 \frac{1 + \chi_{31} + \chi_{31} \chi_{32} + \chi_{31} \chi_{32} \chi_{33}}{\chi_{33}(1 + \chi_{31})}\\
 {\scriptstyle \chi_{23}(1 + \chi_{12})(1 + \chi_{33})}
\end{bmatrix}
,
\endgroup
\quad
   \renewcommand{\arraystretch}{1.2}
\begin{bmatrix}
\chi_{11}\\ \chi_{12} \\ \chi_{13} \\ \chi_{21} \\ \chi_{22} \\ \chi_{23} \\ \chi_{31} \\ \chi_{32} \\ \chi_{33}
    \end{bmatrix}
\xmapsto{\sigma_3}
  \begingroup
    \renewcommand{\arraystretch}{1.8}
\begin{bmatrix}
 \frac{\chi_{11} \chi_{12} \chi_{23}}{(1 + \chi_{11}) (1 + \chi_{23})}\\
 \frac{1}{\chi_{23}}\\
 {\scriptstyle \chi_{13} (1 + \chi_{23})}\\
 \frac{\chi_{21} (1 + \chi_{31})}{1 + \chi_{11} + \chi_{11} \chi_{21} + \chi_{11} \chi_{21} \chi_{31}}\\
 \frac{\chi_{22}(1 + \chi_{11}) (1 + \chi_{23}) (1 + \chi_{31} + \chi_{31} \chi_{32})}{1 + \chi_{31}}\\
 \frac{\chi_{23} \chi_{31} \chi_{32} \chi_{33}}{(1 + \chi_{23}) (1 + \chi_{31} + \chi_{31} \chi_{32})}\\
 \frac{1 + \chi_{11}}{\chi_{11} \chi_{21} (1 + \chi_{31})}\\
 \frac{\chi_{32}(1 + \chi_{11} + \chi_{11} \chi_{21} + \chi_{11} \chi_{21} \chi_{31})}{(1 + \chi_{11}) (1 + \chi_{31} + \chi_{31} \chi_{32})}\\
 \frac{1 + \chi_{31}}{\chi_{31} \chi_{32}}
\end{bmatrix}
.
\endgroup
\notag
\end{align}
Tropicalising the first map above with `max' convention gives the piecewise linear map $Q^+_{\sigma_3}$ which acts on g-vectors as $Q^+_{\sigma_3} {\bf g}(x) = {\bf g} (\sigma_3(x))$. Tropicalising the second map above with `max' convention gives us the map $Q^+_{{\sigma_3}^{-1}} = (Q^+_{\sigma_3})^{-1}$ which acts on g-vectors as $Q^+_{\sigma_3^{-1}} {\bf g}(x) = {\bf g}(\sigma_3^{-1}(x))$.

When we write g-vectors, we use the following ordered initial cluster variables:  
\begin{align*}
& P_{1235},\ P_{1245},\ P_{1345},\ P_{1236},\ P_{1256},\ P_{1456},\ P_{1237},\ P_{1267},\ P_{1567}.  
\end{align*}   

Let ${\bf g}_1$ and ${\bf g}_2$ be the truncated g-vectors of the stable fixed points $ \vcenter{\hbox{\scalebox{0.6}{\begin{ytableau}
    1 & 3 \\ 2 & 5 \\ 4 & 7 \\ 6 & 8
\end{ytableau}}}}$, $\vcenter{\hbox{\scalebox{0.6}{\begin{ytableau}
    1 & 2 \\ 3 & 4 \\ 5 & 6 \\ 7 & 8
\end{ytableau}}}}$ of $\sigma_1$ and $\sigma_2$, respectively. Then 
\begin{align*}
& {\bf g}_1 = (-1, 1, 0, 1, 0, -1, 0, -1, 1), \quad {\bf g}_2 = (0, -1, 0, -1, 0, 1, 0, 1, 0).
\end{align*}
%We have that $g_1$ and $g_2$ are fixed points for the braid generators $\sigma_1$, $\sigma_2$ respectively. 

We have that
\begin{align*}
{\bf g}_3:=\sigma_2^{-1} ({\bf g}_1) = (1, 0, -1, 0, -2, 1, -1, 1, 1), \quad {\bf g}_4:=\sigma_1^{-1} ({\bf g}_2) = (-1, -1, 1, -1, 2, 0, 1, 0, -1),
\end{align*}
and ${\bf g}_2 = {\bf g}_3+{\bf g}_4$. The tableaux corresponding to ${\bf g}_3,{\bf g}_4$ are
\begin{align*}
T_3 = \vcenter{\hbox{\scalebox{0.6}{\begin{ytableau}
1 & 1 & 2 & 3 \\
2 & 4 & 4 & 5 \\
3 & 6 & 6 & 7 \\
5 & 7 & 8 & 8
\end{ytableau}}}}, \quad T_4 = \vcenter{\hbox{\scalebox{0.6}{ \begin{ytableau}
1 & 1 & 2 & 4 \\
2 & 3 & 3 & 6 \\
4 & 5 & 5 & 7 \\
6 & 7 & 8 & 8
\end{ytableau}}}},
\end{align*}
respectively. The vectors ${\bf g}_1,{\bf g}_3,{\bf g}_4$ are linearly independent. 

We have that
\begin{align*}
& \sigma_1^{-1}({\bf g}_1) = {\bf g}_1, \ \sigma_1({\bf g}_1) = {\bf g}_1, \ \sigma_1^{-1}({\bf g}_2) = {\bf g}_4, \ \sigma_1({\bf g}_2) = {\bf g}_3, \ \sigma_1^{-1}({\bf g}_3) = {\bf g}_3+{\bf g}_4 = {\bf g}_2, \\
& \sigma_1({\bf g}_3) = {\bf g}_1+{\bf g}_3, \ \sigma_1^{-1}({\bf g}_4) = {\bf g}_1+{\bf g}_4, \ \sigma_1({\bf g}_4) = {\bf g}_3+{\bf g}_4 = {\bf g}_2, \ \sigma_2^{-1}({\bf g}_1) = {\bf g}_3, \\
& \sigma_2({\bf g}_1) = {\bf g}_4, \ \sigma_2^{-1}({\bf g}_2)={\bf g}_2, \ \sigma_2({\bf g}_2)={\bf g}_2, \ \sigma_2^{-1}({\bf g}_3) = 2{\bf g}_3+{\bf g}_4 = {\bf g}_2+{\bf g}_3, \\
& \sigma_2({\bf g}_3) = {\bf g}_1, \ \sigma_2^{-1}({\bf g}_4) = {\bf g}_1, \ \sigma_2({\bf g}_4) = {\bf g}_3+2{\bf g}_4 = {\bf g}_2 + {\bf g}_4. 
\end{align*}
Moreover, for all $v \in \{{\bf g}_1, {\bf g}_2, {\bf g}_3, {\bf g}_4\}$, we have $\sigma_1 (v) = \sigma_3 (v)$ and $\sigma_1^{-1} (v) = \sigma_3^{-1} (v)$. 
According to the above computations, we consider g-vectors of the form $b{\bf g}+c{\bf g}'$, where $\{{\bf g},{\bf g}'\}$ is one of the following
\begin{align*}
    \{{\bf g}_1,{\bf g}_3\}, \ \{{\bf g}_1,{\bf g}_4\}, \ \{{\bf g}_3,{\bf g}_4\},
\end{align*}
and $b,c \in \ZZ_{\ge 0}$.

We first consider the g-vectors of the form $b{\bf g}_3+c{\bf g}_4$. The vector $b{\bf g}_3 + c{\bf g}_4$ is 
\begin{align*}
( b - c, -c, -b + c, -c, -2b + 2c, b, -b + c, b, b - c ). 
\end{align*}
The corresponding monomial is (see Section \ref{subsec:correspondence between modules and tableaux}) 
\begin{align*}
& Y_{1,-1}^{b-c} ( Y_{1,-1}Y_{1,-3} )^{-c} ( Y_{1,-1}Y_{1,-3}Y_{1,-5} )^{c-b} Y_{2,0}^{-c} ( Y_{2,0}Y_{2,-2} )^{-2b+2c} \times  \\
&  \qquad \qquad \times ( Y_{2,0}Y_{2,-2}Y_{2,-4} )^{b}  Y_{3,1}^{c-b} (Y_{3,1}Y_{3,-1})^{b} ( Y_{3,1}Y_{3,-1}Y_{3,-3} )^{b-c} \\
& = Y_{1,-1}^{-c} Y_{1,-3}^{-b} Y_{2,0}^{-b+c} Y_{2,-2}^{-b+2c} Y_{2,-4}^{b} Y_{3,1}^{b} Y_{3,-1}^{2b-c} Y_{3,-3}^{b-c}.
\end{align*}
First suppose that $b>c$. By Lemma \ref{lem: relation of g vectors and tableaux}, we multiply appropriate frozen variables to the above monomial and obtain a dominant monomial which does not have frozen factors:
\begin{align} \label{eq:dominant monomial after multiplying frozen variables Gr48 b greater than c}
Y_{1,-1}^{b-c} Y_{1,-5}^{c} Y_{1,-7}^{b} Y_{2,-2}^{c} Y_{2,-4}^{2b-c} Y_{2,-6}^{b-c} Y_{3,-1}^{2b-c} Y_{3,1}^{b} Y_{3,-3}^{b-c},
\end{align}
where the frozen variables we use are $b$ copies of $Y_{1,-1}Y_{1,-3}Y_{1,-5}Y_{1,-7}$, and $b-c$ copies of $Y_{2,0}Y_{2,-2}Y_{2,-4}Y_{2,-6}$. 

The dominant monomial (\ref{eq:dominant monomial after multiplying frozen variables Gr48 b greater than c}) decomposes into the following set of increasing dominant monomials with minimal size: $b$ copies of $Y_{1,-7}Y_{2,-4}Y_{3,-1}$, $c$ copies of $Y_{1,-5}Y_{2,-2}Y_{3,1}$, $b-c$ copies of $Y_{2,-6}Y_{3,-3}Y_{3,-1}Y_{3,1}$, $b-c$ copies of $Y_{1,-1}$, $b-c$ copies of $Y_{2,-4}$. Therefore for vectors of the form $b{\bf g}_3+c{\bf g}_4$ the degree of the corresponding polynomial in Pl\"{u}cker coordinates is $4b-2c$ if $b>c$. 

Similarly, for vectors of the form $b{\bf g}_3+c{\bf g}_4$, the degree of the corresponding polynomial in Pl\"{u}cker coordinates is $4c-2b$ if $b<c$. 

We plot the g-vectors on the space spanned by ${\bf g}_3, {\bf g}_4$ in Figure \ref{fig: prime nonreal g vectors obtained by braid group action bg3 plus cg4 Gr48}. The position $(b,c)$ corresponds to the g-vector $b{\bf g}_3+c{\bf g}_4$. The number on each position is the degree of the corresponding polynomial (degrees in Pl\"{u}cker coordinates). The bullets in Figure \ref{fig: prime nonreal g vectors obtained by braid group action bg3 plus cg4 Gr48} represent $(b,c)$ that are not coprime and they correspond to multiples of g-vectors which are already on the picture.

\begin{lemma} \label{lem:bg3 plus cg4 all obtained by braid group action}
For any coprime pairs $(b,c) \in \ZZ_{\ge 1}^2$, $b{\bf g}_3 + c{\bf g}_4$ can be obtained from braid group $\Br_4$ action on ${\bf g}_1, {\bf g}_2$. 
\end{lemma}

\begin{figure}
    \centering
    \adjustbox{scale=0.7,center}{%
    \begin{tikzcd}
	c \\
	10 & \bullet & 38 & \bullet & 34 & \bullet & \bullet & \bullet & 26 & \bullet & 22 & \bullet \\
	9 & \bullet & 34 & 32 & \bullet & 28 & 26 & \bullet & 22 & 20 & \bullet & 22 \\
	8 & \bullet & 30 & \bullet & 26 & \bullet & 22 & \bullet & 18 & \bullet & 20 & \bullet \\
	7 & \bullet & 26 & 24 & 22 & 20 & 18 & 16 & \bullet & 18 & 22 & 26 \\
	6 & \bullet & 22 & \bullet & \bullet & \bullet & 14 & \bullet & 16 & \bullet & \bullet & \bullet \\
	5 & \bullet & 18 & 16 & 14 & 12 & \bullet & 14 & 18 & 22 & 26 & \bullet \\
	4 & \bullet & 14 & \bullet & 10 & \bullet & 12 & \bullet & 20 & \bullet & 28 & \bullet \\
	3 & \bullet & 10 & 8 & \bullet & 10 & 14 & \bullet & 22 & 26 & \bullet & 34 \\
	2 & \bullet & 6 & \bullet & 8 & \bullet & 16 & \bullet & 24 & \bullet & 32 & \bullet \\
	1 & 4 & 2 & 6 & 10 & 14 & 18 & 22 & 26 & 30 & 34 & 38 \\
	0 & \bullet & 4 & \bullet & \bullet & \bullet & \bullet & \bullet & \bullet & \bullet & \bullet & \bullet \\
	& 0 & 1 & 2 & 3 & 4 & 5 & 6 & 7 & 8 & 9 & 10 & b
	\arrow["{\sigma_1^{-1}}"{description}, from=11-3, to=11-2]
	\arrow["{\sigma_1}"{description}, from=11-3, to=12-3]
	\arrow["{\sigma_2}"{description}, from=11-2, to=10-3]
	\arrow["{\sigma_2}"{description}, from=10-3, to=9-4]
	\arrow["{\sigma_2}"{description}, from=9-4, to=8-5]
	\arrow["{\sigma_2}"{description}, from=8-5, to=7-6]
	\arrow["{\sigma_2}"{description}, from=7-6, to=6-7]
	\arrow["{\sigma_2}"{description, pos=0.7}, from=6-7, to=5-8]
	\arrow["{\sigma_2}"{description}, from=5-8, to=4-9]
	\arrow["{\sigma_2}"{description}, from=4-9, to=3-10]
	\arrow["{\sigma_2}"{description}, from=3-10, to=2-11]
	\arrow["{\sigma_2^{-1}}"{description}, from=12-3, to=11-4]
	\arrow["{\sigma_2^{-1}}"{description}, from=11-4, to=10-5]
	\arrow["{\sigma_2^{-1}}"{description}, from=10-5, to=9-6]
	\arrow["{\sigma_2^{-1}}"{description}, from=9-6, to=8-7]
	\arrow["{\sigma_2^{-1}}"{description}, from=8-7, to=7-8]
	\arrow["{\sigma_2^{-1}}"{description}, from=7-8, to=6-9]
	\arrow["{\sigma_2^{-1}}"{description}, from=6-9, to=5-10]
	\arrow["{\sigma_2^{-1}}"{description}, from=5-10, to=4-11]
	\arrow["{\sigma_2^{-1}}"{description}, from=4-11, to=3-12]
	\arrow["{\sigma_1^{-1}}"{description}, from=11-4, to=10-3]
	\arrow["{\sigma_1}"{description}, from=9-4, to=11-5]
	\arrow["{\sigma_1}"{description, pos=0.3}, from=8-5, to=11-6]
	\arrow["{\sigma_1}"{description}, from=7-6, to=11-7]
	\arrow["{\sigma_1}"{description}, from=6-7, to=11-8]
	\arrow["{\sigma_1}"{description, pos=0.2}, from=5-8, to=11-9]
	\arrow["{\sigma_1}"{description}, from=4-9, to=11-10]
	\arrow["{\sigma_1}"{description, pos=0.7}, from=3-10, to=11-11]
	\arrow["{\sigma_1}"{description}, from=2-11, to=11-12]
	\arrow["{\sigma_1^{-1}}"{description}, from=10-5, to=9-3]
	\arrow["{\sigma_1^{-1}}"{description, pos=0.3}, from=9-6, to=8-3]
	\arrow["{\sigma_1^{-1}}"{description}, from=8-7, to=7-3]
	\arrow["{\sigma_1^{-1}}"{description}, from=7-8, to=6-3]
	\arrow["{\sigma_1^{-1}}"{description, pos=0.4}, from=6-9, to=5-3]
	\arrow["{\sigma_1^{-1}}"{description}, from=5-10, to=4-3]
	\arrow["{\sigma_1^{-1}}"{description, pos=0.6}, from=4-11, to=3-3]
	\arrow["{\sigma_1^{-1}}"{description}, from=3-12, to=2-3]
	\arrow["{\sigma_2}"{description}, bend left=20, from=9-3, to=7-5]
	\arrow["{\sigma_2}"{description, pos=0.3}, bend left=20, from=7-5, to=5-7]
	\arrow["{\sigma_2}"{description}, bend left=20, from=5-7, to=3-9]
	\arrow["{\sigma_2}"{description}, bend left=20, from=8-3, to=5-6]
	\arrow["{\sigma_2}"{description}, bend left=20, from=5-6, to=2-9]
	\arrow["{\sigma_2}"{description, pos=0.6}, bend left=20, from=7-4, to=4-7]
	\arrow["{\sigma_2}"{description}, bend left=20, from=7-3, to=3-7]
	\arrow["{\sigma_2^{-1}}"{description}, bend left=-20, from=11-5, to=9-7]
	\arrow["{\sigma_2^{-1}}"{description, pos=0.3}, bend left=-20, from=9-7, to=7-9]
	\arrow["{\sigma_2^{-1}}"{description, pos=0.4}, bend left=-20, from=7-9, to=5-11]
	\arrow["{\sigma_2^{-1}}"{description}, bend left=-20, from=11-6, to=8-9]
	\arrow["{\sigma_2^{-1}}"{description}, bend left=-20, from=8-9, to=5-12]
	\arrow["{\sigma_2^{-1}}"{description}, bend left=-20, from=11-7, to=7-11]
	\arrow["{\sigma_2^{-1}}"{description, pos=0.6}, bend left=-20, from=10-7, to=7-10]
	\arrow["{\sigma_1^{-1}}"{description}, from=9-7, to=7-4]
	\arrow["{\sigma_1^{-1}}"{description}, from=7-9, to=5-4]
	\arrow["{\sigma_1^{-1}}"{description, pos=0.4}, from=5-11, to=3-4]
	\arrow["{\sigma_1}"{description}, from=7-5, to=10-7]
	\arrow["{\sigma_1^{-1}}"{description}, from=8-9, to=5-5]
	\arrow["{\sigma_1^{-1}}"{description, pos=0.7}, from=7-10, to=4-5]
	\arrow["{\sigma_1^{-1}}"{description}, from=5-12, to=2-5]
	\arrow["{\sigma_1}"{description, pos=0.7}, from=5-6, to=9-9]
	\arrow["{\sigma_1^{-1}}"{description}, from=7-11, to=3-6]
	\arrow["{\sigma_1}"{description, pos=0.3}, from=5-7, to=10-9]
	\arrow["{\sigma_1}"{description, pos=0.2}, from=4-7, to=9-10]
	\arrow["{\sigma_1}"{description, pos=0.7}, from=3-7, to=8-11]
	\arrow["{\sigma_1}"{description, pos=0.3}, from=3-9, to=10-11]
	\arrow["{\sigma_1}"{description}, from=2-9, to=9-12]
\end{tikzcd} }
   \caption{The g-vectors of the form $b{\bf g}_3 + c{\bf g}_4$ obtained by the braid group action for $\Gr(4,8)$. We draw arrows which correspond to the braid group action. The underlying undirected graph of the quiver in the figure is symmetric.}
            \label{fig: prime nonreal g vectors obtained by braid group action bg3 plus cg4 Gr48 with maps sigmai}
\end{figure}
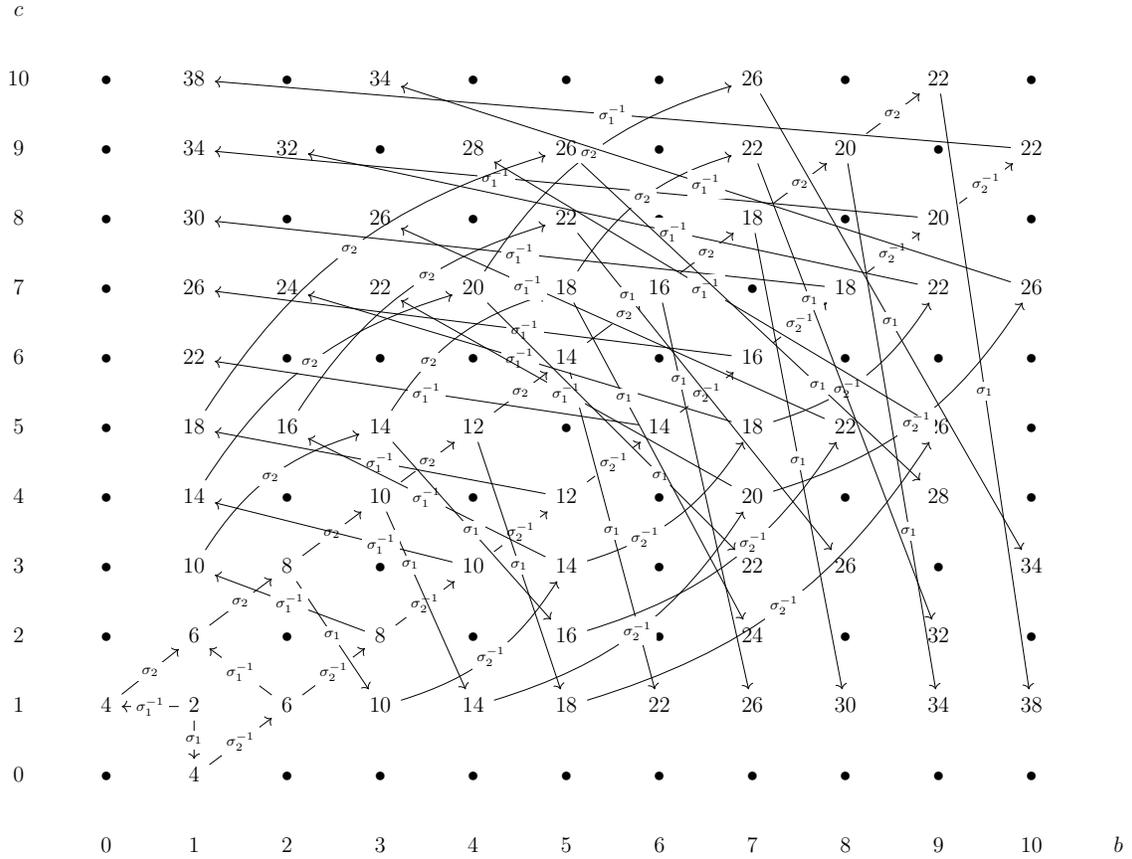

\begin{proof}
We draw the orbits of the braid group action in Figure 
\ref{fig: prime nonreal g vectors obtained by braid group action bg3 plus cg4 Gr48 with maps sigmai}.
Applying $\sigma_1^{-1}$ to ${\bf g}_2$ gives ${\bf g}_4$, and then applying $\sigma_2$ repeatedly to ${\bf g}_4$ yields
$b{\bf g}_3 + (b+1){\bf g}_4$, $b \in \mathbb{Z}_{\ge 1}$.
Similarly, applying $\sigma_1$ to ${\bf g}_2$ gives ${\bf g}_3$, and applying $\sigma_2^{-1}$ repeatedly to ${\bf g}_3$ gives
$(c+1){\bf g}_3 + c{\bf g}_4$, $c \in \mathbb{Z}_{\ge 1}$.
Next, applying $\sigma_1^{-1}$ to $(c+1){\bf g}_3 + c{\bf g}_4$, $c \ge 1$, produces
${\bf g}_3 + d{\bf g}_4$, $d \in \mathbb{Z}_{\ge 2}$,
which correspond to the g-vectors above the diagonal and in Column 1 of Figure 
\ref{fig: prime nonreal g vectors obtained by braid group action bg3 plus cg4 Gr48 with maps sigmai}.
Similarly, applying $\sigma_1$ to $b{\bf g}_3 + (b+1){\bf g}_4$, $b \ge 1$, gives
$d{\bf g}_3 + {\bf g}_4$, $d \in \mathbb{Z}_{\ge 2}$,
which are the g-vectors below the diagonal and in Row 1 of the figure.
Applying $\sigma_2^{-1}$ repeatedly to $3{\bf g}_3 + {\bf g}_4$ produces
$(3+2d){\bf g}_3 + (1+2d){\bf g}_4$, $d \in \mathbb{Z}_{\ge 1}$,
and applying $\sigma_1^{-1}$ to these vectors yields the g-vectors above the diagonal and in Column 2.
Similarly, applying $\sigma_2$ repeatedly to ${\bf g}_3 + 3{\bf g}_4$ gives
$(1+2d){\bf g}_3 + (3+2d){\bf g}_4$, $d \in \mathbb{Z}_{\ge 1}$,
and applying $\sigma_1$ produces the g-vectors below the diagonal and in Row 2.

Applying $\sigma_2^{-1}$ repeatedly to $4{\bf g}_3 + {\bf g}_4$ and $5{\bf g}_3 + 2{\bf g}_4$ produces
$(4+3d){\bf g}_3 + (1+3d){\bf g}_4$, $(5+3d){\bf g}_3 + (2+3d){\bf g}_4$, $d \in \mathbb{Z}_{\ge 1}$.
Applying $\sigma_1^{-1}$ to these vectors yields the g-vectors above the diagonal in Column 3.  
Similarly, applying $\sigma_2$ repeatedly to ${\bf g}_3 + 4{\bf g}_4$ and $2{\bf g}_3 + 5{\bf g}_4$ gives
$(1+3d){\bf g}_3 + (4+3d){\bf g}_4$, $(2+3d){\bf g}_3 + (5+3d){\bf g}_4, \quad d \in \mathbb{Z}_{\ge 1}$,
and applying $\sigma_1$ to these vectors yields the g-vectors below the diagonal in Row 3. Continuing this procedure, we obtain all g-vectors of the form
$b{\bf g}_3 + c{\bf g}_4$, $b,c \in \mathbb{Z}_{\ge 1}$ are coprime.
\end{proof}

We can also see highlighted the rays with a given degree in Figure \ref{fig: prime nonreal g vectors obtained by braid group action bg3 plus cg4 Gr48}. There are Euler totient $\phi(r)$ of them, for degree $2r$. Indeed, the number of coprime pairs $(b,c)$ ($b>c$) such that $4b-2c=2r$ ($r \ge 3$) is $\frac{1}{2}\phi(r)$. The number of coprime pairs $(b,c)$ ($b<c$) such that $4c-2b=2r$ ($r \ge 3$) is also $\frac{1}{2}\phi(r)$. 

Now we consider the g-vectors of the form $b{\bf g}_1+b'{\bf g}_3$, $c{\bf g}_1+c'{\bf g}_4$, $b,b',c,c' \in \ZZ_{\ge 0}$. The case of some of $b,b',c,c'$ are $0$ is clear. We consider the case that $b,b',c,c' \in \ZZ_{\ge 1}$. Similar to Lemma \ref{lem:bg3 plus cg4 all obtained by braid group action}, we have the following.
\begin{lemma} \label{lem:bg1 plus cg3, bg1 plus cg4 Gr48 all obtained by braid group action}
For any coprime pairs $(b,c) \in \ZZ_{\ge 1}^2$, $b{\bf g}_1 + c{\bf g}_3$ and $b{\bf g}_1+c{\bf g}_4$ can be obtained from braid group $\Br_4$ action on ${\bf g}_1, {\bf g}_2$. 
\end{lemma}

Note that $\sigma_1 (v) = \sigma_3 (v)$ and $\sigma_1^{-1} (v) = \sigma_3^{-1} (v)$ for $v \in \{{\bf g}_1, {\bf g}_2, {\bf g}_3, {\bf g}_4\}$, from the proofs of Lemmas \ref{lem:bg3 plus cg4 all obtained by braid group action} and \ref{lem:bg1 plus cg3, bg1 plus cg4 Gr48 all obtained by braid group action}, we also have the following. 
\begin{lemma} \label{lem:all g vectors obtained by braid group action are in the vector space Gr48}
The g-vectors obtained by the braid group action on ${\bf g}_1, {\bf g}_2$ are of the form $bg+cg'$, where $\{g,g'\}$ is one of the following
\begin{align*}
    \{{\bf g}_1,{\bf g}_3\}, \ \{{\bf g}_1,{\bf g}_4\}, \ \{{\bf g}_3,{\bf g}_4\},
\end{align*}
and $b,c \in \ZZ_{\ge 0}$ are coprime. 
\end{lemma}

The number of degree $2r$, $r \ge 3$, g-vectors of this form is $A_r + B_r$, where $A_r$ is the number of g-vectors with degree $2r$ of the form $b{\bf g}_1+b'{\bf g}_3$, and $B_r$ is the number of g-vectors with degree $2r$ of the form $c{\bf g}_1+c'{\bf g}_4$. Both of the numbers $A_r, B_r$ are equal to the number of pairs $(e,f)$ of coprime positive numbers such that $2e+4f=2r$, $r \ge 3$. Therefore $A_r+B_r$ is twice of the number of pairs $(e,f)$ of coprime positive numbers such that $e+2f=r$, $r \ge 3$. This is the value $\phi(r)$ of the Euler totient function, $r \ge 3$.  

\begin{remark}
Since ${\bf g}_2={\bf g}_3+{\bf g}_4$, the set $\{b{\bf g}_3+c{\bf g}_4: b, c \in \ZZ_{\ge 1}, b>c\}$ is the same as the set $\{c {\bf g}_2 + c' {\bf g}_3: c, c' \ge \ZZ_{\ge 1}\}$, and the set $\{b{\bf g}_3+c{\bf g}_4: b, c \in \ZZ_{\ge 1}, b<c\}$ is the same as the set $\{b {\bf g}_2 + b' {\bf g}_4: b, b' \ge \ZZ_{\ge 1}\}$.  
\end{remark}

\begin{figure}
\scalebox{0.7}{
\begin{tikzpicture}[scale=1] 
    % labels
    \foreach \i in {0,...,8}
      \path[blue] (\i,-1) node{\i} (-1,\i) node{\i};
    % loop over the lattice points
    % \foreach \i in {0,...,8}
    %   \foreach \j in {0,...,8}{
    %     \draw (\i,\j) circle(3pt);    
  %    }; 

     \node at (0,1) {$4$};

     \foreach \i in {0,2,3,4,5,6,7,8}
     \fill (0,\i) circle[radius=2pt]; 

     \foreach \i in {0,2,4,6,8}
     \fill (2,\i) circle[radius=2pt]; 

     \foreach \i in {0,3,6}
     \fill (3,\i) circle[radius=2pt]; 

     \foreach \i in {0,2,4,6,8}
     \fill (4,\i) circle[radius=2pt]; 

     \foreach \i in {0,5}
     \fill (5,\i) circle[radius=2pt]; 

     \foreach \i in {0,2,3,4,6,8}
     \fill (6,\i) circle[radius=2pt]; 

     \foreach \i in {0,7}
     \fill (7,\i) circle[radius=2pt]; 

      \foreach \i in {0,2,4,6,8}
     \fill (8,\i) circle[radius=2pt]; 

     % \node at (0,0) {$2$};
     % \node at (0,1) {$4$};
     % \node at (0,2) {$6$};
     % \node at (0,3) {$8$};
     % \node at (0,4) {$10$};
     % \node at (0,5) {$12$};
     % \node at (0,6) {$14$};
     % \node at (0,7) {$16$};
     % \node at (0,8) {$18$};

     \node at (1,0) {$4$};
     \node at (1,1) {$2$};
     \node at (1,2) {$6$};
     \node at (1,3) {$10$};
     \node at (1,4) {$14$};
     \node at (1,5) {$18$};
     \node at (1,6) {$22$};
     \node at (1,7) {$26$};
     \node at (1,8) {$30$};

     \node at (2,1) {$6$};
     \node at (2,3) {$8$};
     \node at (2,5) {$16$};
     \node at (2,7) {$24$};

     \node at (3,1) {$10$};
     \node at (3,2) {$8$};
     \node at (3,4) {$10$};
     \node at (3,5) {$14$};
     \node at (3,7) {$22$};
     \node at (3,8) {$26$};

     \node at (4,1) {$14$};
     \node at (4,3) {$10$};
     \node at (4,5) {$12$};
     \node at (4,7) {$20$};

     \node at (5,1) {$18$};
     \node at (5,2) {$16$};
     \node at (5,3) {$14$};
     \node at (5,4) {$12$}; 
     \node at (5,6) {$14$};
     \node at (5,7) {$18$};
     \node at (5,8) {$22$};

     \node at (6,1) {$22$};
     \node at (6,5) {$14$};
     \node at (6,7) {$16$};

    \node at (7,1) {$26$};
     \node at (7,2) {$24$};
     \node at (7,3) {$22$};
     \node at (7,4) {$20$};
     \node at (7,5) {$18$}; 
     \node at (7,6) {$16$}; 
     \node at (7,8) {$18$};

     \node at (8,1) {$30$};
     \node at (8,3) {$26$}; 
     \node at (8,5) {$22$}; 
     \node at (8,7) {$18$};

     \draw[-] (8,8)--(0,4);
     \draw[-] (8,8)--(4,0);

     \draw[-] (5,5)--(0,2.5);
     \draw[-] (5,5)--(2.5,0);

\end{tikzpicture} }
            \caption{The g-vectors of the form $b{\bf g}_3 + c{\bf g}_4$ obtained by the braid group action for $\Gr(4,8)$. In the case that the greatest common factor of $b,c$ is $1$, the number at coordinate $(b,c)$ is the degree of the polynomial in Pl\"{u}cker coordinates corresponding to the g-vector $b{\bf g}_3 + c{\bf g}_4$. If the greatest common factor of $b,c$ is greater than $1$, then we put a ``$\bullet$'' at the position $(b,c)$.}
            \label{fig: prime nonreal g vectors obtained by braid group action bg3 plus cg4 Gr48}
\end{figure}
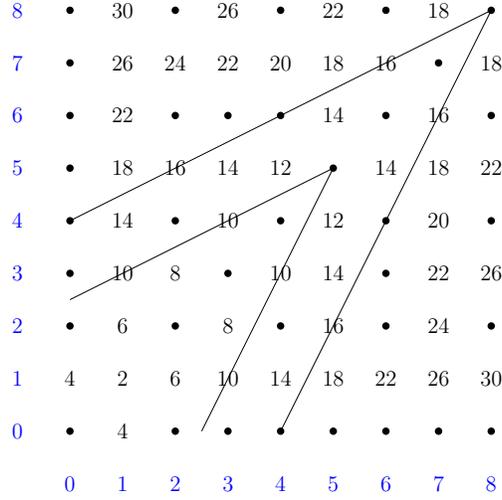

\subsection{\texorpdfstring{Proof of Theorem \ref{thm:prime non real modules Gr39 and Gr48} for $\Gr(3,9)$}{Proof of Theorem \ref{thm:prime non real modules Gr39 and Gr48} for Gr(3,9)}}
 \label{subsec:prime nonreal tableaux for Gr39}

Acting on the $\hat{y}$-variables of the initial seed with $\sigma_1^{-1}$ we find
\begin{align}
\sigma_1^{-1} :
   \renewcommand{\arraystretch}{1.2}
\begin{bmatrix}
\chi_{11}\\ \chi_{12} \\ \chi_{21} \\ \chi_{22} \\ \chi_{31} \\ \chi_{32} \\ \chi_{41} \\ \chi_{42} \\ \chi_{51} \\ \chi_{52}
    \end{bmatrix}
\mapsto
  \begingroup
    \renewcommand{\arraystretch}{1.8}
\begin{bmatrix}
 \frac{\chi_{22} (1 + \chi_{32} + \chi_{32} \chi_{42} + \chi_{32} \chi_{42} \chi_{52})}{
  1 + \chi_{12} + \chi_{12} \chi_{22} + \chi_{12} \chi_{22} \chi_{32} + \chi_{12} \chi_{22} \chi_{32} \chi_{42} + 
   \chi_{12} \chi_{22} \chi_{32} \chi_{42} \chi_{52}}\\
 \frac{1 + \chi_{12}}{\chi_{12} \chi_{22} (1 + \chi_{32} + \chi_{32} \chi_{42} + \chi_{32} \chi_{42} \chi_{52})}\\
 \frac{\chi_{11} (1 + \chi_{12}) (1 + \chi_{32} + \chi_{21} \chi_{32} + \chi_{32} \chi_{42} + \chi_{21} \chi_{32} \chi_{42} + 
     \chi_{32} \chi_{42} \chi_{52} + \chi_{21} \chi_{32} \chi_{42} \chi_{52})}{
  1 + \chi_{32} + \chi_{32} \chi_{42} + \chi_{32} \chi_{42} \chi_{52}}\\
 \frac{\chi_{21} (1 + \chi_{12} + \chi_{12} \chi_{22} + \chi_{12} \chi_{22} \chi_{32} + \chi_{12} \chi_{22} \chi_{32} \chi_{42} + 
     \chi_{12} \chi_{22} \chi_{32} \chi_{42} \chi_{52})}{(1 + \chi_{12}) (1 + \chi_{32} + \chi_{21} \chi_{32} + \chi_{32} \chi_{42} + 
     \chi_{21} \chi_{32} \chi_{42} + \chi_{32} \chi_{42} \chi_{52} + \chi_{21} \chi_{32} \chi_{42} \chi_{52})}\\
 \frac{\chi_{21} \chi_{31} \chi_{32} \chi_{42} (1 + \chi_{42} + \chi_{42} \chi_{52})}{(1 + \chi_{42}) (1 + \chi_{32} + 
     \chi_{21} \chi_{32} + \chi_{32} \chi_{42} + \chi_{21} \chi_{32} \chi_{42} + \chi_{32} \chi_{42} \chi_{52} + 
     \chi_{21} \chi_{32} \chi_{42} \chi_{52})}\\
 \frac{1 + \chi_{32} + \chi_{32} \chi_{42} + \chi_{32} \chi_{42} \chi_{52}}{\chi_{21} \chi_{32} (1 + \chi_{42} + \chi_{42} \chi_{52})}\\
 \frac{\chi_{52}}{1 + \chi_{42} + \chi_{42} \chi_{52}}\\
 \frac{(1 + \chi_{42}) (1 + \chi_{32} + \chi_{21} \chi_{32} + \chi_{32} \chi_{42} + \chi_{21} \chi_{32} \chi_{42} + 
     \chi_{32} \chi_{42} \chi_{52} + \chi_{21} \chi_{32} \chi_{42} \chi_{52})}{\chi_{42} \chi_{52}}\\
{\scriptstyle  \chi_{41} (1 + \chi_{42})}\\
 \frac{\chi_{51} (1 + \chi_{42} + \chi_{42} \chi_{52})}{1 + \chi_{42}}
\end{bmatrix}
.
\endgroup
\notag
\end{align}
Tropicalising the above map with `max' convention give the piecewise linear map $Q^+_{\sigma_1}$ which acts on g-vectors as $Q^+_{\sigma_1} {\bf g}(x) = {\bf g} (\sigma_1(x))$. 

Acting instead with $\sigma_1$ on the $\hat{y}$-variables of the initial seed we find
\begin{align}
\sigma_1 :
   \renewcommand{\arraystretch}{1.2}
\begin{bmatrix}
\chi_{11}\\ \chi_{12} \\ \chi_{21} \\ \chi_{22} \\ \chi_{31} \\ \chi_{32} \\ \chi_{41} \\ \chi_{42} \\ \chi_{51} \\ \chi_{52}
    \end{bmatrix}
\mapsto
  \begingroup
    \renewcommand{\arraystretch}{1.8}
\begin{bmatrix}
 \frac{\chi_{11} \chi_{21} \chi_{32} (1 + \chi_{12}) }{(1 + \chi_{32})(1 + \chi_{11} + \chi_{11} \chi_{12}) }\\
\frac{1}{\chi_{11} (1 + \chi_{12})}\\
 \frac{1 + \chi_{12} + \chi_{12} \chi_{22} + \chi_{12} \chi_{22} \chi_{32}}{\chi_{32}(1 + \chi_{12}) }\\
 \frac{\chi_{22} (1 + \chi_{32})(1 + \chi_{11} + \chi_{11} \chi_{12})}{
  1 + \chi_{12} + \chi_{12} \chi_{22} + \chi_{12} \chi_{22} \chi_{32}}\\
 \frac{\chi_{31} (1 + \chi_{32}) (1 + \chi_{12} + \chi_{12} \chi_{22} + \chi_{12} \chi_{22} \chi_{32} + \chi_{41} + \chi_{12} \chi_{41} + 
     \chi_{12} \chi_{22} \chi_{41} + \chi_{12} \chi_{22} \chi_{32} \chi_{41} + \chi_{12} \chi_{22} \chi_{32} \chi_{41} \chi_{42})}{
  1 + \chi_{12} + \chi_{12} \chi_{22} + \chi_{12} \chi_{22} \chi_{32}}\\
 \frac{\chi_{32} \chi_{41} \chi_{42}(1 + \chi_{12})}{(1 + \chi_{32}) (1 + \chi_{12} + \chi_{12} \chi_{22} + 
     \chi_{12} \chi_{22} \chi_{32} + \chi_{41} + \chi_{12} \chi_{41} + \chi_{12} \chi_{22} \chi_{41} + \chi_{12} \chi_{22} \chi_{32} \chi_{41} + 
     \chi_{12} \chi_{22} \chi_{32} \chi_{41} \chi_{42})}\\
 \frac{\chi_{41} \chi_{51}(1 + \chi_{12} + \chi_{12} \chi_{22} + \chi_{12} \chi_{22} \chi_{32} + \chi_{12} \chi_{22} \chi_{32} \chi_{42}) }{
  1 + \chi_{12} + \chi_{12} \chi_{22} + \chi_{12} \chi_{22} \chi_{32} + \chi_{41} + \chi_{12} \chi_{41} + \chi_{12} \chi_{22} \chi_{41} + 
   \chi_{12} \chi_{22} \chi_{32} \chi_{41} + \chi_{12} \chi_{22} \chi_{32} \chi_{41} \chi_{42}}\\
 \frac{1 + \chi_{12} + \chi_{12} \chi_{22} + \chi_{12} \chi_{22} \chi_{32}}{
  \chi_{41} (1 + \chi_{12} + \chi_{12} \chi_{22} + \chi_{12} \chi_{22} \chi_{32} + \chi_{12} \chi_{22} \chi_{32} \chi_{42})}\\
 \frac{\chi_{12} \chi_{22} \chi_{32} \chi_{42} \chi_{52}}{
  1 + \chi_{12} + \chi_{12} \chi_{22} + \chi_{12} \chi_{22} \chi_{32} + \chi_{12} \chi_{22} \chi_{32} \chi_{42}}\\
 \frac{1 + \chi_{12} + \chi_{12} \chi_{22} + \chi_{12} \chi_{22} \chi_{32} + \chi_{41} + \chi_{12} \chi_{41} + \chi_{12} \chi_{22} \chi_{41} + 
   \chi_{12} \chi_{22} \chi_{32} \chi_{41} + \chi_{12} \chi_{22} \chi_{32} \chi_{41} \chi_{42}}{\chi_{12} \chi_{22} \chi_{32} \chi_{42}}
\end{bmatrix}
.
\endgroup
\notag
\end{align}
Tropicalising the above map with `max' convention gives us the map $Q^+_{{\sigma_1}^{-1}} = (Q^+_{\sigma_1})^{-1}$ which acts on g-vectors as $Q^+_{\sigma_1^{-1}} {\bf g}(x) = {\bf g}(\sigma_1^{-1}(x))$. 

Acting on the $\hat{y}$-variables of the initial seed with $\sigma_2^{-1}$ we find
\begin{align}
\sigma_2^{-1} :
   \renewcommand{\arraystretch}{1.2}
\begin{bmatrix}
\chi_{11}\\ \chi_{12} \\ \chi_{21} \\ \chi_{22} \\ \chi_{31} \\ \chi_{32} \\ \chi_{41} \\ \chi_{42} \\ \chi_{51} \\ \chi_{52}
    \end{bmatrix}
\mapsto
  \begingroup
    \renewcommand{\arraystretch}{1.8}
\begin{bmatrix}
\frac{\chi_{22}}{1 + \chi_{21} + \chi_{21} \chi_{22} + \chi_{21} \chi_{31} + \chi_{21} \chi_{22} \chi_{31} + \chi_{21} \chi_{31} \chi_{41} + \chi_{21} \chi_{22} \chi_{31} \chi_{41} + \chi_{21} \chi_{31} \chi_{41} \chi_{51} + \chi_{21} \chi_{22} \chi_{31} \chi_{41} \chi_{51}}\\
{\scriptstyle \chi_{11} (1 + \chi_{21} + \chi_{21} \chi_{31} + \chi_{21} \chi_{31} \chi_{41} + \chi_{21} \chi_{31} \chi_{41} \chi_{51})}\\
\frac{1 + \chi_{21} + \chi_{21} \chi_{31} + \chi_{21} \chi_{31} \chi_{41} + \chi_{21} \chi_{31} \chi_{41} \chi_{51}}{\chi_{21} \chi_{22} (1 + \chi_{31} + \chi_{31} \chi_{41} + \chi_{31} \chi_{41} \chi_{51})}\\
\frac{\chi_{12} (1 + \chi_{21} + \chi_{21} \chi_{22} + \chi_{21} \chi_{31} + \chi_{21} \chi_{22} \chi_{31} + \chi_{21} \chi_{31} \chi_{41} + 
     \chi_{21} \chi_{22} \chi_{31} \chi_{41} + \chi_{21} \chi_{31} \chi_{41} \chi_{51} + \chi_{21} \chi_{22} \chi_{31} \chi_{41} \chi_{51})}{1 + \chi_{21} + \chi_{21} \chi_{31} + \chi_{21} \chi_{31} \chi_{41} + \chi_{21} \chi_{31} \chi_{41} \chi_{51}}\\
\frac{(1 + \chi_{31}) (1 + \chi_{21} + \chi_{21} \chi_{22} + \chi_{21} \chi_{31} + \chi_{21} \chi_{22} \chi_{31} + 
     \chi_{21} \chi_{31} \chi_{41} + \chi_{21} \chi_{22} \chi_{31} \chi_{41} + \chi_{21} \chi_{31} \chi_{41} \chi_{51} + 
     \chi_{21} \chi_{22} \chi_{31} \chi_{41} \chi_{51})}{\chi_{31} \chi_{41} (1 + \chi_{51})}\\
\frac{\chi_{21} \chi_{22} \chi_{31} \chi_{32} (1 + \chi_{31} + \chi_{31} \chi_{41} + \chi_{31} \chi_{41} \chi_{51})}{(1 + \chi_{31}) (1 +  \chi_{21} + \chi_{21} \chi_{22} + \chi_{21} \chi_{31} + \chi_{21} \chi_{22} \chi_{31} + \chi_{21} \chi_{31} \chi_{41} + \chi_{21} \chi_{22} \chi_{31} \chi_{41} + \chi_{21} \chi_{31} \chi_{41} \chi_{51} + \chi_{21} \chi_{22} \chi_{31} \chi_{41} \chi_{51})}\\
 \frac{\chi_{52}(1 + \chi_{31} + \chi_{31} \chi_{41} + \chi_{31} \chi_{41} \chi_{51}) }{(1 + \chi_{31}) (1 + \chi_{51} + \chi_{51} \chi_{52})}\\
 \frac{\chi_{41} (1 + \chi_{51})}{1 + \chi_{31} + \chi_{31} \chi_{41} + \chi_{31} \chi_{41} \chi_{51}}\\
 \frac{1 + \chi_{51}}{\chi_{51} \chi_{52}}\\
\frac{\chi_{42} (1 + \chi_{31})  (1 + \chi_{51} + \chi_{51} \chi_{52})}{1 + \chi_{51}}
\end{bmatrix}
.
\endgroup
\notag
\end{align}
Tropicalising the above map with `max' convention give the piecewise linear map $Q^+_{\sigma_2}$ which acts on g-vectors as $Q^+_{\sigma_2} {\bf g}(x) = {\bf g} (\sigma_2(x))$. 

Acting instead with $\sigma_2$ on the $\hat{y}$-variables of the initial seed we find
\begin{align}
\sigma_2 :
   \renewcommand{\arraystretch}{1.2}
\begin{bmatrix}
\chi_{11}\\ \chi_{12} \\ \chi_{21} \\ \chi_{22} \\ \chi_{31} \\ \chi_{32} \\ \chi_{41} \\ \chi_{42} \\ \chi_{51} \\ \chi_{52}
    \end{bmatrix}
\mapsto
  \begingroup
    \renewcommand{\arraystretch}{1.8}
\begin{bmatrix}
 \frac{\chi_{11} \chi_{12} (1 + \chi_{21})}{1 + \chi_{11} + \chi_{11} \chi_{21}}\\
 \frac{\chi_{21} \chi_{22}}{1 + \chi_{21}}\\
 \frac{\chi_{21} \chi_{31} \chi_{42}}{(1 + \chi_{21}) (1 + \chi_{11} + \chi_{11} \chi_{21} + \chi_{42} + \chi_{11} \chi_{42} +  \chi_{11} \chi_{21} \chi_{42} + \chi_{11} \chi_{21} \chi_{31} \chi_{42})}\\
 \frac{1 + \chi_{11} + \chi_{11} \chi_{21}}{\chi_{21}}\\
 \frac{1 + \chi_{11} + \chi_{11} \chi_{21}}{\chi_{42}(1 + \chi_{11} + \chi_{11} \chi_{21} + \chi_{11} \chi_{21} \chi_{31})}\\
 \frac{\chi_{32} (1 + \chi_{21}) (1 + \chi_{11} + \chi_{11} \chi_{21} + \chi_{42} + \chi_{11} \chi_{42} + \chi_{11} \chi_{21} \chi_{42} +  \chi_{11} \chi_{21} \chi_{31} \chi_{42})}{1 + \chi_{11} + \chi_{11} \chi_{21}}\\
 \frac{\chi_{41} (1 + \chi_{51})(1 + \chi_{11} + \chi_{11} \chi_{21} + \chi_{42} + \chi_{11} \chi_{42} + \chi_{11} \chi_{21} \chi_{42} + \chi_{11} \chi_{21} \chi_{31} \chi_{42})}{1 + \chi_{11} + \chi_{11} \chi_{21} + \chi_{11} \chi_{21} \chi_{31} + \chi_{11} \chi_{21} \chi_{31} \chi_{41} + \chi_{11} \chi_{21} \chi_{31} \chi_{41} \chi_{51}}\\
 \frac{\chi_{42} \chi_{51} \chi_{52} (1 + \chi_{11} + \chi_{11} \chi_{21} + \chi_{11} \chi_{21} \chi_{31}) }{(1 + \chi_{51}) (1 + \chi_{11} + \chi_{11} \chi_{21} + \chi_{42} + \chi_{11} \chi_{42} + \chi_{11} \chi_{21} \chi_{42} + \chi_{11} \chi_{21} \chi_{31} \chi_{42}) }\\
 \frac{1 + \chi_{11} + \chi_{11} \chi_{21} + \chi_{11} \chi_{21} \chi_{31}}{\chi_{11} \chi_{21} \chi_{31} \chi_{41} (1 + \chi_{51})},\\
 \frac{1 + \chi_{11} + \chi_{11} \chi_{21} + \chi_{11} \chi_{21} \chi_{31} + \chi_{11} \chi_{21} \chi_{31} \chi_{41} + 
   \chi_{11} \chi_{21} \chi_{31} \chi_{41} \chi_{51}}{\chi_{51} (1 + \chi_{11} + \chi_{11} \chi_{21} + \chi_{11} \chi_{21} \chi_{31}) }
\end{bmatrix}
.
\endgroup
\notag
\end{align}
Tropicalising the above map with `max' convention gives us the map $Q^+_{{\sigma_2}^{-1}} = (Q^+_{\sigma_2})^{-1}$ which acts on g-vectors as $Q^+_{\sigma_2^{-1}} {\bf g}(x) = {\bf g}(\sigma_2^{-1}(x))$.

When we write g-vectors, we use the following ordered initial cluster variables: 
\begin{align*}
P_{124},\ P_{134},\ P_{125},\ P_{145},\ P_{126},\ P_{156},\ P_{127},\ P_{167},\ P_{128},\ P_{178}.
\end{align*}

Let 
\begin{align*}
&   T_1 = \vcenter{\hbox{\scalebox{0.6}{ \begin{ytableau}
1 & 3 & 4 \\
2 & 6 & 7 \\
5 & 8 & 9
\end{ytableau}}}}, \
T_2 = \vcenter{\hbox{\scalebox{0.6}{ \begin{ytableau}
1 & 2 & 5 \\
3 & 4 & 8 \\
6 & 7 & 9
\end{ytableau}}}}, \
T_3 = \vcenter{\hbox{\scalebox{0.6}{ \begin{ytableau}
1 & 2 & 3 \\
4 & 5 & 6 \\
7 & 8 & 9
\end{ytableau}}}}, 
\end{align*}
\begin{align*}
&  T_4= \vcenter{\hbox{\scalebox{0.6}{
\begin{ytableau}
1 & 1 & 2 & 3 & 3 & 4 \\
2 & 5 & 5 & 6 & 6 & 7 \\
4 & 7 & 8 & 8 & 9 & 9
\end{ytableau}}}}, \
T_5 = \vcenter{\hbox{\scalebox{0.6}{ \begin{ytableau}
1 & 1 & 2 & 4 & 4 & 5 \\
2 & 3 & 3 & 7 & 7 & 8 \\
5 & 6 & 6 & 8 & 9 & 9
\end{ytableau}}}}, \
T_6 = \vcenter{\hbox{\scalebox{0.6}{ \begin{ytableau}
1 & 1 & 2 & 2 & 3 & 6 \\
3 & 4 & 4 & 5 & 5 & 8 \\
6 & 7 & 7 & 8 & 9 & 9
\end{ytableau}}}}.
\end{align*}

Their g-vectors are 
\begin{align*}
& {\bf g}_1=(0, 0, 1, -1, 0, -1, -1, 1, 0, 1), \quad {\bf g}_2=(-1, 0, -1, 1, 1, 0, 1, -1, 0, 0), \\
& {\bf g}_3=(0, -1, -1, 0, -1, 1, 0, 1, 1, 0), \quad {\bf g}_4=(1, -1, 0, -2, -2, 1, -1, 2, 1, 1), \\
& {\bf g}_5=(-2, 1, 1, 0, 2, -2, 0, -1, -1, 2), \quad {\bf g}_6=(-1, -1, -2, 1, -1, 2, 2, 0, 1, -1),
\end{align*}
respectively. The g-vectors ${\bf g}_1, \ldots, {\bf g}_6$ are linearly independent. 

We have that 
\begin{align*}
& \sigma_1^{-1} {\bf g}_1 = {\bf g}_1, \ \sigma_1^{-1} {\bf g}_2 = {\bf g}_5, \ \sigma_1^{-1} {\bf g}_3 = {\bf g}_2, \ \sigma_2^{-1} {\bf g}_1 = {\bf g}_3, \ \sigma_2^{-1} {\bf g}_2 = {\bf g}_2, \ \sigma_2^{-1} {\bf g}_3 = {\bf g}_6, \\
& \sigma_1 {\bf g}_1 = {\bf g}_1, \ \sigma_1 {\bf g}_2 = {\bf g}_3, \ \sigma_1 {\bf g}_3 = {\bf g}_4, \ \sigma_2 {\bf g}_1 = {\bf g}_5, \ \sigma_2 {\bf g}_2 = {\bf g}_2, \ \sigma_2 {\bf g}_3 = {\bf g}_1,
\end{align*}  
\begin{align*}
& \sigma_1^{-1} {\bf g}_4 = {\bf g}_3, \ \sigma_1^{-1} {\bf g}_5 = {\bf g}_1+{\bf g}_5, \ \sigma_1^{-1} {\bf g}_6 = {\bf g}_2+{\bf g}_5, \ \sigma_2^{-1} {\bf g}_4 = {\bf g}_3+{\bf g}_6, \ \sigma_2^{-1} {\bf g}_5 = {\bf g}_1, \ \sigma_2^{-1} {\bf g}_6 = {\bf g}_2+{\bf g}_6,  
\end{align*}
\begin{align*}
& \sigma_1 {\bf g}_4 = {\bf g}_1+{\bf g}_4, \ \sigma_1 {\bf g}_5 = {\bf g}_2, \ \sigma_1 {\bf g}_6 = {\bf g}_3+{\bf g}_4, \\
& \sigma_2 {\bf g}_4 =  {\bf g}_1+{\bf g}_5, \ \sigma_2 {\bf g}_5 = {\bf g}_2+{\bf g}_5, \ \sigma_2 {\bf g}_6 = {\bf g}_3.
\end{align*}
By the above computations, we consider g-vectors which are the form $c {\bf g} + c' {\bf g}'$, where $c, c' \in \ZZ_{\ge 0}$, $\{{\bf g}, {\bf g}'\}$ is one of the following 
\begin{align*}
\{{\bf g}_1, {\bf g}_4\}, \{{\bf g}_1, {\bf g}_5\}, \{{\bf g}_2, {\bf g}_5\}, \{{\bf g}_2, {\bf g}_6\}, \{{\bf g}_3, {\bf g}_4\}, \{{\bf g}_3, {\bf g}_6\}.
\end{align*}
We draw a picture of g-vectors of the form $c {\bf g}_1 + c' {\bf g}_4$ in Figure \ref{fig: prime nonreal g vectors obtained by braid group action bg1 plus cg4 Gr39}. 
The following lemma is similar to Lemmas \ref{lem:bg3 plus cg4 all obtained by braid group action} and \ref{lem:bg1 plus cg3, bg1 plus cg4 Gr48 all obtained by braid group action}.
\begin{lemma} \label{lem:bg1 plus cg4, bg1 plus cg5 all obtained by braid group action Gr39}
For any coprime pairs $(b,c) \in \ZZ_{\ge 1}^2$, $b{\bf g}_1 + c{\bf g}_4$, $b{\bf g}_1 + c{\bf g}_5$, $b{\bf g}_2 + c{\bf g}_5$, $b{\bf g}_2 + c{\bf g}_6$, $b{\bf g}_3 + c{\bf g}_4$, $b{\bf g}_3 + c{\bf g}_6$ can be obtained from braid group $\Br_3$ action on ${\bf g}_1, {\bf g}_2, {\bf g}_3$. 
\end{lemma}

\begin{proof}
Note that $\{T_2,T_5\}, \{T_3,T_6\}$ are promotions of $\{T_1, T_4\}$, and $\{T_2,T_6\}$, $\{T_3, T_4\}$ are promotions of $\{T_1, T_5\}$. It suffices to prove the cases of $b{\bf g}_1 + c{\bf g}_4$ and $b{\bf g}_1 + c{\bf g}_5$. We prove the case of $b{\bf g}_1 + c{\bf g}_4$. The case of $b{\bf g}_1 + c{\bf g}_5$ is similar. 

We draw orbits of the braid group action in Figure \ref{fig: prime nonreal g vectors obtained by braid group action bg1 plus cg4 Gr39}. Applying $\sigma_1$ repeatedly to ${\bf g}_4$ produces the g-vectors $b{\bf g}_1 + {\bf g}_4$, $b \in \mathbb{Z}_{\ge 1}$, which lie in Row 1. We have $\sigma_2^{-1}({\bf g}_3) = {\bf g}_6$, $\sigma_1({\bf g}_6) = {\bf g}_3 + {\bf g}_4$, and $\sigma_1({\bf g}_3 + {\bf g}_4) = {\bf g}_1 + 2{\bf g}_4$, which is the g-vector at position $(1,2)$. Applying $\sigma_1$ repeatedly to this vector produces the g-vectors in Row 2.
Similarly, we have $\sigma_2^{-1}({\bf g}_6) = {\bf g}_2 + {\bf g}_6$, $\sigma_1({\bf g}_2 + {\bf g}_6) = 2{\bf g}_3 + {\bf g}_4$, and $\sigma_1(2{\bf g}_3 + {\bf g}_4) = {\bf g}_1 + 3{\bf g}_4$, which is the g-vector at position $(1,3)$.
Next, we have $\sigma_1({\bf g}_3) = {\bf g}_4$, $\sigma_2^{-1}({\bf g}_4) = {\bf g}_3 + {\bf g}_6$, $\sigma_1({\bf g}_3 + {\bf g}_6) = {\bf g}_3 + 2{\bf g}_4$, and $\sigma_1({\bf g}_3 + 2{\bf g}_4) = 2{\bf g}_1 + 3{\bf g}_4$, which is the g-vector at position $(2,3)$.
Applying $\sigma_1$ repeatedly to the g-vectors at positions $(1,3)$ and $(2,3)$ produces the g-vectors in Row 3.

Continue this procedure, we obtain g-vectors of the form $b{\bf g}_1+c{\bf g}_4$, where $b,c \in \ZZ_{\ge 1}$ are coprime. 
\end{proof}

Similar to Lemma \ref{lem:all g vectors obtained by braid group action are in the vector space Gr48}, we have the following.
\begin{lemma} \label{lem:all g vectors obtained by braid group action are in the vector space Gr39}
The g-vectors obtained by the braid group action on ${\bf g}_1, {\bf g}_2, {\bf g}_3$ are of the form $bg+cg'$, where $\{g,g'\}$ is one of the following
\begin{align*}
\{{\bf g}_1, {\bf g}_4\}, \{{\bf g}_1, {\bf g}_5\}, \{{\bf g}_2, {\bf g}_5\}, \{{\bf g}_2, {\bf g}_6\}, \{{\bf g}_3, {\bf g}_4\}, \{{\bf g}_3, {\bf g}_6\},
\end{align*}
and $b,c \in \ZZ_{\ge 0}$ are coprime. 
\end{lemma}

Note that $\{T_2,T_5\}, \{T_3,T_6\}$ are promotions of $\{T_1, T_4\}$, and $\{T_2,T_6\}$, $\{T_3, T_4\}$ are promotions of $\{T_1, T_5\}$. The number of degree $3r$, $r \ge 3$, g-vectors is $3A_r + 3 B_r$, where $A_r$ is the number of g-vectors with degree $3r$ of the form $c{\bf g}_1+c'{\bf g}_4$, and $B_r$ is the number of g-vectors with degree $3r$ of the form $c{\bf g}_1+c'{\bf g}_5$. The number $A_r$ is equal to the number of pairs $(a,b)$ of coprime positive numbers such that $6a+3b=3r$, $r \ge 3$. The number $B_r$ is also equal to the number of coprime pairs $(a,b)$ of positive numbers such that $6a+3b=3r$, $r \ge 3$. Therefore $A_r+B_r$ is twice of the number of pairs $(a,b)$ of coprime positive numbers such that $2a+b=r$, $r \ge 3$. This is the value $\phi(r)$ of the Euler totient function, $r \ge 3$.

\begin{figure}
\scalebox{0.7}{
\begin{tikzpicture}[scale=1] 
    % labels
    \foreach \i in {0,...,8}
      \path[blue] (\i,-1) node{\i} (-1,\i) node{\i};
    % loop over the lattice points
    % \foreach \i in {0,...,8}
    %   \foreach \j in {0,...,8}{
    %     \draw (\i,\j) circle(3pt);    
  %    }; 

     \node at (0,1) {$6$}; 

     \foreach \i in {0,2,3,4,5,6,7,8}
     \fill (0,\i) circle[radius=2pt]; 

     \foreach \i in {0,2,4,6,8}
     \fill (2,\i) circle[radius=2pt]; 

     \foreach \i in {0,3,6}
     \fill (3,\i) circle[radius=2pt]; 

     \foreach \i in {0,2,4,6,8}
     \fill (4,\i) circle[radius=2pt]; 

     \foreach \i in {0,5}
     \fill (5,\i) circle[radius=2pt]; 

     \foreach \i in {0,2,3,4,6,8}
     \fill (6,\i) circle[radius=2pt]; 

     \foreach \i in {0,7}
     \fill (7,\i) circle[radius=2pt]; 

      \foreach \i in {0,2,4,6,8}
     \fill (8,\i) circle[radius=2pt];

     \node at (1,0) {$3$};
     \node at (1,1) {$9$};
     \node at (1,2) {$15$};
     \node at (1,3) {$21$};
     \node at (1,4) {$27$};
     \node at (1,5) {$33$};
     \node at (1,6) {$39$};
     \node at (1,7) {$45$};
     \node at (1,8) {$51$};

     \node at (2,1) {$12$};
     \node at (2,3) {$24$};
     \node at (2,5) {$36$};
     \node at (2,7) {$48$};
 
     \node at (3,1) {$15$};
     \node at (3,2) {$21$};
     \node at (3,4) {$33$};
     \node at (3,5) {$39$};
     \node at (3,7) {$51$};
     \node at (3,8) {$57$};
 
     \node at (4,1) {$18$};
     \node at (4,3) {$30$};
     \node at (4,5) {$42$};
     \node at (4,7) {$54$};
 
     \node at (5,1) {$21$};
     \node at (5,2) {$27$};
     \node at (5,3) {$33$};
     \node at (5,4) {$39$}; 
     \node at (5,6) {$51$};
     \node at (5,7) {$57$};
     \node at (5,8) {$63$};
 
     \node at (6,1) {$24$};
     \node at (6,5) {$48$};
     \node at (6,7) {$60$};
 
    \node at (7,1) {$27$};
     \node at (7,2) {$33$};
     \node at (7,3) {$39$};
     \node at (7,4) {$45$};
     \node at (7,5) {$51$}; 
     \node at (7,6) {$57$}; 
     \node at (7,8) {$69$};
 
     \node at (8,1) {$30$};
     \node at (8,3) {$42$}; 
     \node at (8,5) {$54$}; 
     \node at (8,7) {$66$}; 

     \draw[-] (1,4)--(7,1);
     \draw[-] (2,3)--(6,1);
     \draw[-] (1,3)--(5,1);
  
\end{tikzpicture} }
            \caption{The g-vectors of the form $b {\bf g}_1 + c{\bf g}_4$ obtained by the braid group action for $\Gr(3,9)$. In the case that the greatest common factor of $b,c$ is $1$, the number at coordinate $(b,c)$ is the degree of the polynomial in Pl\"{u}cker coordinates corresponding to the g-vector $b{\bf g}_1 + c{\bf g}_4$. If the greatest common factor of $b,c$ is greater than $1$, then we put a ``$\bullet$'' at the position $(b,c)$.}
            \label{fig: prime nonreal g vectors obtained by braid group action bg1 plus cg4 Gr39}
\end{figure}
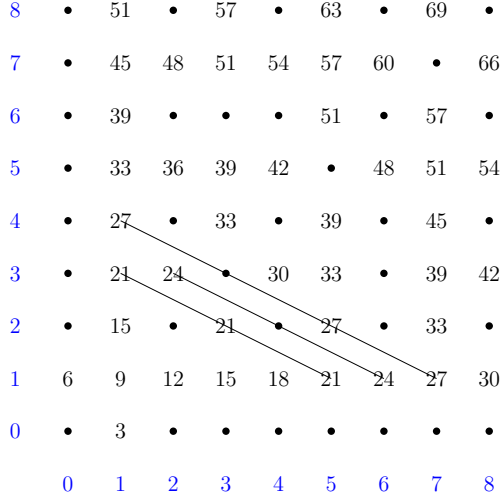

\begin{figure}
    \centering
    \adjustbox{scale=0.7,center}{%
    \begin{tikzcd}
	c \\
	8 & \bullet & 51 & \bullet & 57 & \bullet & 63 & \bullet & 69 & \bullet \\
	7 & \bullet & 45 & 48 & 51 & 54 & 57 & 60 & \bullet & 66 \\
	6 & \bullet & 39 & \bullet & \bullet & \bullet & 51 & \bullet & 57 & \bullet \\
	5 & \bullet & 33 & 36 & 39 & 42 & \bullet & 48 & 51 & 54 \\
	4 & \bullet & 27 & \bullet & 33 & \bullet & 39 & \bullet & 45 & \bullet \\
	3 & \bullet & 21 & 24 & \bullet & 30 & 33 & \bullet & 39 & 42 \\
	2 & \bullet & 15 & \bullet & 21 & \bullet & 27 & \bullet & 33 & \bullet \\
	1 & 6 & 9 & 12 & 15 & 18 & 21 & 24 & 27 & 30 \\
	0 & \bullet & 3 & \bullet & \bullet & \bullet & \bullet & \bullet & \bullet & \bullet \\
	& 0 & 1 & 2 & 3 & 4 & 5 & 6 & 7 & 8 & b
	\arrow["{\sigma_1}"{description}, from=9-2, to=9-3]
	\arrow["{\sigma_1}"{description}, from=9-3, to=9-4]
	\arrow["{\sigma_1}"{description}, from=9-4, to=9-5]
	\arrow["{\sigma_1}"{description}, from=9-5, to=9-6]
	\arrow["{\sigma_1}"{description}, from=9-6, to=9-7]
	\arrow["{\sigma_1}"{description}, from=9-7, to=9-8]
	\arrow["{\sigma_1}"{description}, from=9-8, to=9-9]
	\arrow["{\sigma_1}"{description}, from=9-9, to=9-10]
	\arrow["{\sigma_1}"{description}, bend left = 30, from=8-3, to=8-5]
	\arrow["{\sigma_1}"{description}, bend left = 30, from=8-5, to=8-7]
	\arrow["{\sigma_1}"{description}, bend left = 30, from=8-7, to=8-9]
	\arrow["{\sigma_1}"{description}, bend left = 15, from=7-3, to=7-6]
	\arrow["{\sigma_1}"{description}, bend left = 15, from=7-4, to=7-7]
	\arrow["{\sigma_1}"{description}, bend left = 15, from=7-7, to=7-10]
	\arrow["{\sigma_1}"{description}, bend left = 15, from=7-6, to=7-9]
	\arrow["{\sigma_1}"{description}, bend left = 15, from=6-3, to=6-7]
	\arrow["{\sigma_1}"{description}, bend left = 15, from=6-5, to=6-9]
	\arrow["{\sigma_1}"{description}, bend left = 15, from=5-3, to=5-8]
	\arrow["{\sigma_1}"{description}, bend left = 15, from=5-4, to=5-9]
	\arrow["{\sigma_1}"{description}, bend left = 15, from=5-5, to=5-10]
	\arrow["{\sigma_1}"{description}, bend left = 15, from=4-3, to=4-9]
	\arrow["{\sigma_1}"{description}, bend left = 15, from=3-3, to=3-10]
\end{tikzcd}
    }
   \caption{The g-vectors of the form $b{\bf g}_1 + c{\bf g}_4$ obtained by the braid group action for $\Gr(3,9)$. We draw arrows which correspond to the braid group action.}
            \label{fig: prime nonreal g vectors obtained by braid group action bg1 plus cg4 Gr39 with maps sigmai}
\end{figure}
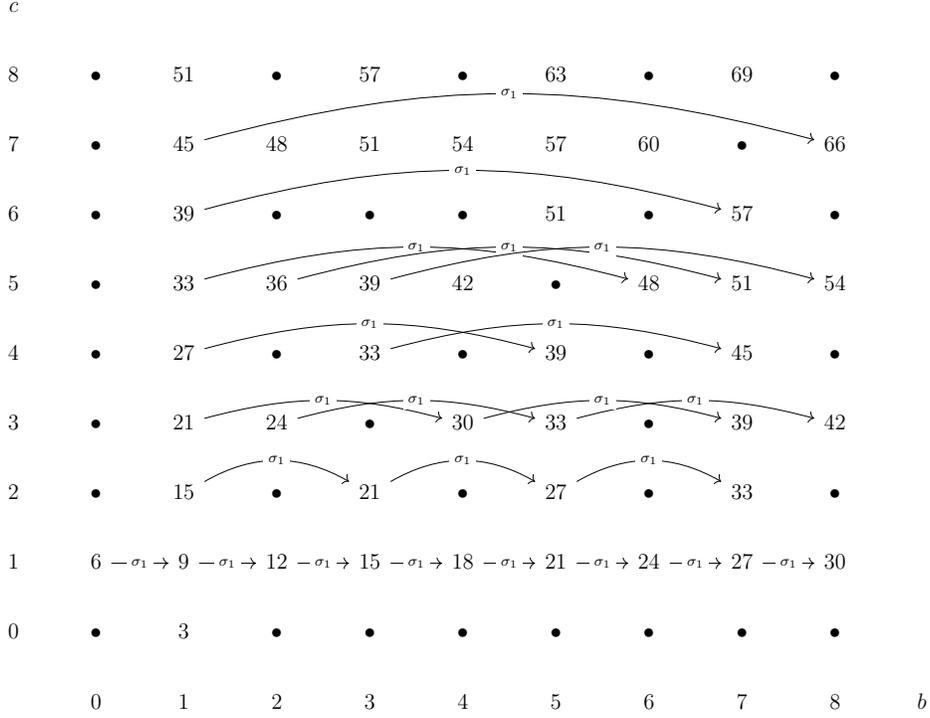

\section{Connection to scattering amplitudes} \label{sec:connection to scattering amplitudes}

In scattering amplitudes in physics, the singularities of planar loop amplitudes in $N = 4$ super Yang-Mills are related to Grassmannian cluster algebras $\CC[\Gr(4,n)]$. Cluster algebras are related to the branch cut singularities of amplitudes, see \cite{GGSVV, GPSV}. This connection relates the cluster variables of the cluster algebra with the rational symbol letters (potential logarithmic branch cuts) of the scattering amplitude. In the case when the cluster algebra is of infinite type, apart from rational symbol letters, there are algebraic letters. Algebraic letters are studied in \cite{ALS21, DFGK20, DFGK21, HP20, HP21}. In the case of $\Gr(4,8)$, algebraic letters are of the form
\begin{align*}
\frac{z_0 + B_z \sqrt{\Delta}}{z_0 - B_z \sqrt{\Delta}},
\end{align*}
where $B_z = \frac{2 z_1 - z_0 \ch(T)}{\Delta}$, $z_0, z_1$ are certain cluster variables, $T = \vcenter{\hbox{\scalebox{0.6}{ \begin{ytableau}
1 & 3 \\ 2 & 5 \\ 4 & 7 \\ 6 & 8
\end{ytableau}}}}$ or $\vcenter{\hbox{\scalebox{0.6}{ \begin{ytableau}
1 & 2 \\ 3 & 4 \\ 5 & 6 \\ 7 & 8
\end{ytableau}}}}$, $\Delta = A^2 - 4B$, $A = \ch(T)$, $B = \ch(T)^2 - \ch(T \cup T)$, see \cite{DFGK21, HP20}. The element $\sqrt{\Delta}$ is called a square root and it appears in the four-mass box integral \cite{Hodges79, HV79}. 

In \cite{DFGK21, HP20}, $\Delta$ is found using an infinite mutation sequence. In \cite{ALS21}, $\Delta$ is found using dual canonical basis of $\CC[\Gr(4,8)]$ and the formula of $\ch(T)$ in \cite{CDFL}. In this paper, we provide another method to find the tableau $\vcenter{\hbox{\scalebox{0.6}{ \begin{ytableau}
1 & 3 \\ 2 & 5 \\ 4 & 7 \\ 6 & 8
\end{ytableau}}}}$ (resp. $\vcenter{\hbox{\scalebox{0.6}{ \begin{ytableau}
1 & 2 \\ 3 & 4 \\ 5 & 6 \\ 7 & 8
\end{ytableau}}}}$) and obtain the corresponding square root $\sqrt{\Delta}$ (resp. the cyclic rotation of $\sqrt{\Delta}$) by using the formula of $\ch(T)$. That is, these tableaux are stable fixed points of the braid group action, see Section \ref{sec:stable fixed points}. Our method also produces stable fixed points for other cluster algebras $\CC[\Gr(k,n)]$. We expect that these stable fixed points will have applications in scattering amplitudes in physics. 
%These questions will be studied in a future publiction. 

\end{document}